\documentclass[leqno]{article}
\usepackage[english]{babel}
\usepackage[utf8]{inputenc}
\usepackage{amsmath}
\usepackage{amssymb}
\usepackage{amsfonts}
\usepackage{enumerate}
\usepackage{vmargin}
\usepackage[all]{xy}
\usepackage{mathrsfs}
\usepackage{mathtools}
\usepackage{lmodern}
\usepackage{slashed}
\usepackage{xurl}
\usepackage[colorlinks=true,linkcolor=blue,pagebackref=true,hyperfootnotes=false]{hyperref}
\setmarginsrb{3cm}{3cm}{3.5cm}{3cm}{0cm}{0cm}{1.5cm}{3cm}
\usepackage{comment}
\usepackage{array}
\usepackage{longtable}
\newlength{\dhatheight}

\newcommand{\Cbb}{\ensuremath{\mathbb{C}}}

\newcommand{\K}{\ensuremath{\mathbb{K}}}

\newcommand{\Q}{\mathrm{Q}}

\newcommand{\diag}{\mathrm{diag}}

\newcommand{\Mat}{\ensuremath{\mathbb{M}}}

\let\cal\relax
\newcommand{\cal}{\mathcal}

\newcommand{\Z}{\ensuremath{\mathbb{Z}}}

\newcommand{\R}{\ensuremath{\mathbb{R}}}

\newcommand{\Id}{\mathrm{Id}}

\newcommand{\la}{\langle}
\newcommand{\ra}{\rangle}

\renewcommand{\leq}{\ensuremath{\leqslant}}
\renewcommand{\geq}{\ensuremath{\geqslant}}
\newcommand{\qed}{\hfill \vrule height6pt  width6pt depth0pt}

\newcommand{\co}{\colon}

\newcommand{\ot}{\otimes}
\newcommand{\ovl}{\overline}

\let\i\relax 
\newcommand{\i}{\mathrm{i}}

\newcommand{\ov}{\overset}

\newcommand{\epsi}{\varepsilon}
\renewcommand{\d}{\mathop{}\mathopen{}\mathrm{d}} 
\newcommand{\e}{\mathrm{e}} 
\renewcommand{\d}{\mathop{}\mathopen{}\mathrm{d}}

\DeclareMathOperator{\Tr}{Tr} 
\let\ker\relax 
\DeclareMathOperator{\ker}{Ker}

\let\Re\relax 
\DeclareMathOperator{\Re}{Re}

\DeclareMathOperator{\Irr}{Irr} 
\DeclareMathOperator{\Rep}{Rep}

\newtheorem{thm}{Theorem}[section]
\newtheorem{defi}[thm]{Definition}
\newtheorem{prop}[thm]{Proposition}

\newtheorem{cor}[thm]{Corollary}
\newtheorem{lemma}[thm]{Lemma}

\newtheorem{remark}[thm]{Remark}

\newenvironment{proof}[1][]{\noindent {\it Proof #1} : }{\hbox{~}\qed
\smallskip
}
\usepackage{tocloft}
\numberwithin{equation}{section}
\usepackage[nottoc,notlot,notlof]{tocbibind}
\let\OLDthebibliography\thebibliography
\renewcommand\thebibliography[1]{
  \OLDthebibliography{#1}
  \setlength{\parskip}{0pt}
  \setlength{\itemsep}{0pt plus 0.3ex}
}
\newcommand\reallywidehat[1]{\arraycolsep=0pt\relax
\begin{array}{c}
\stretchto{
  \scaleto{
    \scalerel*[\widthof{\ensuremath{#1}}]{\kern-.5pt\bigwedge\kern-.5pt}
    {\rule[-\textheight/2]{1ex}{\textheight}} 
  }{\textheight} 
}{0.5ex}\\           
#1\\                 
\rule{-1ex}{0ex}
\end{array}
}

\begin{document}
\selectlanguage{english}
\title{\bfseries{Haagerup vertex operator algebras at central charge eight}}
\date{}
\author{\bfseries{C\'edric Arhancet}}
\maketitle

\begin{abstract}
Evans and Gannon showed that a positive vertex operator algebra of central charge eight whose representation category is the Drinfeld center $Z(\cal H_3)$ of the Haagerup fusion category has one of two character vectors, with six or nineteen currents. Call such an algebra, if strongly rational, a Haagerup algebra. We prove that a Haagerup algebra $U$ has six currents, that $U_1 \cong \mathfrak{sl}_2 \oplus \mathfrak{sl}_2$ embeds into $V_{E_6 \oplus A_2}$ in one of two explicit ways, with levels $(1,39)$ or $(3,13)$, and we determine the lowest weight spaces of its modules and $U_2$ as $U_1$-modules. For the levels $(3,13)$, $U_2$ is one of four Galois conjugate subspaces. Together with $U_1$ it generates the unitary $C_2$-cofinite algebra $S$, which is isomorphic to the algebra constructed by Jiao. For the levels $(1,39)$, $U$ contains the explicit algebra of Riesen. The associated lattice of $U$ is that of the form $x^2+xy+10y^2$, resp.\ $3x^2+3xy+4y^2$, the two classes of the principal genus of discriminant $-39$. Exchanging them in the decomposition of $U$ over its Heisenberg commutant gives a Haagerup algebra $U^\natural$ of the other levels, and $U \mapsto U^\natural$ is a bijection between the isomorphism classes of the two cases. The classification gives exactly two isomorphism classes, forming a single hyperbolic genus, as soon as either explicit candidate is known to be a Haagerup algebra. This property is proved for $S$ in \cite{ArhIII}, where the twin $S^\natural$ is also identified, with the algebra of Huang and Riesen, using the results of Dong, Jiao and Xu \cite{DJX26}. The explicit algebras of Jiao and Riesen have rational forms. The proofs combine exact trace relations given by vector-valued modular forms, exact computations in a model of $V_{E_6 \oplus A_2}$, and genera of vertex operator algebras.
\end{abstract}

\makeatletter
 \renewcommand{\@makefntext}[1]{#1}
 \makeatother

{
  \hypersetup{linkcolor=blue}
 \tableofcontents
}

\section{Introduction}
\label{sec:Introduction}

The Haagerup subfactor of index $\frac{5+\sqrt{13}}{2}$ \cite{AsH99} gives rise to three Morita equivalent fusion categories $\cal H_1$, $\cal H_2$ and $\cal H_3$ \cite{GrS12}. The category $\cal H_3$ has six simple objects $1,a,a^2,\rho,a\rho,a^2\rho$ with fusion rules
\begin{equation}
\label{eq:fusion-rules}
a^3 \cong 1,\quad a\rho \cong \rho a^2,\quad \rho^2 \cong 1 \oplus \rho \oplus a\rho \oplus a^2\rho.
\end{equation}
The Frobenius--Perron dimension of $\rho$ is
\begin{equation}
\label{eq:d}
d \ov{\mathrm{def}}{=} \frac{3+\sqrt{13}}{2}.
\end{equation}
Numerical studies of Haagerup symmetry in conformal field theory include \cite{HL+22,V+22}. At central charge eight, Evans and Gannon \cite{EvG11} showed that a positive vertex operator algebra whose representation category has the modular data of the Drinfeld center $Z(\cal H_3)$ has one of two character vectors, indexed by $\gamma \in \{0,1\}$, with $6+13\gamma$ currents. They observed that such an algebra should lie in $V_{E_6 \oplus A_2}$ \cite[Section 2.4]{EvG11}, and Gannon \cite{Gan19} began the analysis of its currents. Albert, Honda, Kaidi and Zheng \cite{AHKZ26} suggested that the chiral $(\mathfrak e_8)_1$ theory carries Haagerup symmetries. Möller and Rayhaun sketched \cite{MoR24}, and Gannon and Rayhaun proposed \cite{GaR26}, a correspondence between finite non-invertible symmetries and conformal subalgebras of finite index. For the Haagerup categories, Gannon and Rayhaun predicted that a chiral algebra of central charge eight with representation category $Z(\cal H_3)$ arises inside $V_{E_6 \oplus A_2}$ \cite[Section 4.3]{GaR26}. Recently, Jiao \cite{Jia26} gave an explicit vertex operator algebra of central charge eight inside $V_{E_6 \oplus A_2}$, with currents $\mathfrak{sl}_2 \oplus \mathfrak{sl}_2$ of levels $(3,13)$, together with a computer-assisted proof that its representation category is $Z(\cal H_3)$ (preprint of October 4, 2026, which we have not checked in full). Riesen \cite{Rie26} and Huang \cite{Hua26} subsequently each announced such an algebra, with currents of levels $(1,39)$. Evans and Riesen \cite{EvR26} have since announced realizations of the centers of the untwisted Haagerup--Izumi categories of every odd order by non-unitary algebras of negative central charge.

We classify the possible Haagerup algebras and construct the twin of any such algebra by exchanging its associated lattice. The classification starts from either of the two explicit candidates once it is known to be a Haagerup algebra. That the explicit algebra $S$ of Theorem \ref{thm:intro-case-two} is a Haagerup algebra is proved in \cite[Theorem 6.1]{ArhIII}. We call \emph{Haagerup algebra} a positive strongly rational vertex operator algebra $U$ of central charge eight with a braided equivalence $\Rep(U) \cong Z(\cal H_3)$. Positive means that the irreducible modules other than $U$ have positive lowest conformal weights. Let us say that a vector-valued function is an admissible character vector for given modular data if it transforms according to that data, has nonnegative integral Fourier coefficients, has vacuum component $q^{-\frac{c}{24}}(1+O(q))$ and has positive conformal weights in all other components (Definition \ref{def:admissible}). The characters of a Haagerup algebra form such a vector by \cite{Zhu96,Hua08}. Both vectors of \cite{EvG11} pass the tests of character theory known to us (their coefficients are nonnegative integers up to $q^{60}$), so that character theory alone does not seem to decide which of them occurs.

\begin{thm}
\label{thm:intro-gamma}
Let $U$ be a Haagerup algebra. Then the character vector of $U$ is the Evans--Gannon vector $\gamma = 0$. In particular, $U$ has exactly six currents.
\end{thm}

The proof has three steps. First, $U$ extends to $V_{E_8}$ through a Lagrangian algebra of $Z(\cal H_3)$, and the subalgebra of this extension attached to the subcategory $\mathrm{Vec}_{\Z_3} \subset \cal H_3$ is the fixed point algebra of an element of order three of $E_8$. A dimension count identifies it with $V' = V_{E_6 \oplus A_2}$ (Proposition \ref{prop:chain}). Second, the restrictions of the $V'$-modules to $U$ show that the representations $\mathbf 3$ of $\mathfrak{sl}_3$ and $\mathbf{27}$ of $\mathfrak e_6$ are self-dual when restricted to $U_1$ (Proposition \ref{prop:self-dual}). Third, for $h$ in a Cartan subalgebra of $U_1$, the one-point functions of $h[-1]h$ form a weight-two vector-valued modular form whose low-order coefficients are traces of $h_{(0)}^2$ on lowest weight spaces. This is the Killing vector of \cite[Section 2.4]{Gan19}. We determine the space of such forms exactly: the relations between its coefficients are produced by pairing with weakly holomorphic forms of weight zero for the dual representation and applying the residue theorem, in the spirit of the Serre duality of Borcherds \cite[Theorem 3.1]{Bor99}, and their completeness is proved by exhibiting exact elements (Proposition \ref{prop:weight-two-space}). For $\gamma = 1$, the module $d_3$ has a one-dimensional lowest weight space, which forces strong inequalities on every simple ideal of $U_1$ and is incompatible with $\dim U_1 = 19$ (Section \ref{sec:gamma-one}). For $\gamma = 0$, a Heisenberg recursion (Lemma \ref{lem:hermite}) gives relations of weights four and six, which, together with the integrality of the weights of $\mathfrak{sl}_2$, a bound coming from integrability and an explicit construction of $\mathfrak e_8$, lead to the following classification.

\begin{thm}
\label{thm:intro-lie}
Let $U$ be a Haagerup algebra. Then $U_1 \cong \mathfrak{sl}_2 \oplus \mathfrak{sl}_2$, and the inclusion $U_1 \subset \mathfrak e_6 \oplus \mathfrak{sl}_3$ of Proposition \ref{prop:chain} is, up to conjugation, one of two explicit embeddings, with levels $(k',k'') = (1,39)$, case (I), or $(3,13)$, case (II). In each case the lowest weight spaces of the twelve simple modules and the weight-two space $U_2$ are uniquely determined as $U_1$-modules.
\end{thm}

The precise statement is Theorem \ref{thm:weight-one-classification} with Table \ref{tab:lowest-spaces}. In case (I) the table reproduces the data listed in the working manuscript \cite{Rie26b} of Riesen, which are therefore forced by the representation category once the levels are $(1,39)$. Gannon \cite[Section 3.3]{Gan19} began the analysis of the six-dimensional case and treated the level one factor. The levels $13$ and $39$ come from diagonal embeddings into $\mathfrak e_6 \oplus \mathfrak{sl}_3$. We then implement $V'$ in exact arithmetic over $\K = \Q(\sqrt{13},\sqrt{-3})$ and study the Griess algebra of the commutant $C = \operatorname{Com}_U(A)$ of the affine subalgebra $A$ generated by $U_1$ (Section \ref{sec:griess}).

\begin{thm}
\label{thm:intro-case-two}
Let $U$ be a Haagerup algebra in case (II). After an automorphism of $V'$, $C_2$ is spanned by two commuting Virasoro vectors of central charges $\frac{182}{95}$ and $\frac{32}{19}$, and $U_2$ is one of four explicit subspaces of $V'_2$, defined over $\K$, which form one orbit under $\mathrm{Gal}(\K/\Q)$ and generate isomorphic vertex algebras. The vertex algebra $S$ generated by $U_1 \oplus U_2$ is unitary, simple and self-contragredient, strongly generated in weights at most three and $C_2$-cofinite. Its graded dimension agrees with the Haagerup vacuum character up to $q^{11/3}$, so that $U_n = S_n$ for $n \leq 4$. Moreover $S$ is isomorphic to the algebra constructed by Jiao \cite{Jia26}.
\end{thm}

These are Theorems \ref{thm:case-two-U2}, \ref{thm:exact-weight-two}, \ref{thm:unitarity}, \ref{thm:strong-generation} and \ref{thm:jiao-identification}, and Propositions \ref{prop:symmetries} and \ref{prop:centralizer-two}. The finite weight-four condition \textup{(F4)} used for strong generation is stated in Definition \ref{assump:weight-four} and established by the exact certificate of Proposition \ref{prop:weight-four-certificate}. In the setting of Theorem \ref{thm:intro-case-two}, the inclusion $S \subset U$ and the Haagerup character give a second proof by a dimension bound (Corollary \ref{cor:weight-four-from-haagerup}). The proof of identification with Jiao uses only weight-two membership and generation, and is independent of \textup{(F4)}. In case (I), we exclude the last spurious candidate for $C_2$ by a certified computation (Proposition \ref{prop:case-one-two}), identify the three weight-two generators $e$, $f$, $v$ of the algebra $U^{\mathrm R} = \la A,e,f,v\ra$ of \cite{Rie26} with our data (Proposition \ref{prop:riesen-identification}), and exchange the remaining choices by automorphisms of $V'$ (Proposition \ref{prop:centralizer-one}). Hence every Haagerup algebra contains $S$ or $\psi_{\mathrm R}(U^{\mathrm R})$ after an automorphism of $V'$, where $\psi_{\mathrm R}$ identifies the lattice algebra of \cite{Rie26} with $V'$. Theorems \ref{thm:uniqueness-one} and \ref{thm:uniqueness-two} prove equality when the corresponding candidate is a Haagerup algebra.

Our main new result explains the two cases. A strongly rational vertex operator algebra decomposes over the dual pair formed by its Heisenberg commutant $\Omega_U = \operatorname{Com}_U(M(\mathfrak h))$, $\mathfrak h$ a Cartan subalgebra of $U_1$, and the lattice algebra $V_L$ of its associated lattice $L$ \cite{Mas14,MoR24a}. This decomposition was used by Höhn \cite{Hoh17} to study the genus of \cite{Hoh03}, and it leads to the hyperbolic genus of Moriwaki \cite{Mor21} and Möller--Rayhaun \cite{MoR24a}.

\begin{thm}
\label{thm:intro-genus}
Let $U$ be a Haagerup algebra, with currents of levels $(k',k'')$.
\begin{enumerate}
\item The associated lattice of $U$ is the even binary lattice with Gram matrix $\left[\begin{smallmatrix}2&1\\1&20\end{smallmatrix}\right]$ in case (I) and $\left[\begin{smallmatrix}6&3\\3&8\end{smallmatrix}\right]$ in case (II), that is, the lattice of the binary form $x^2+xy+10y^2$, resp.\ $3x^2+3xy+4y^2$. These are the two classes of the principal genus of binary forms of discriminant $-39$. The Heisenberg commutant $\Omega_U$ is a strongly rational vertex operator algebra of central charge six without currents, $\Rep(\Omega_U) \cong Z(\cal H_3) \boxtimes \overline{\cal C(L'/L)}$, and all its irreducible modules have nonnegative conformal weights. Moreover $U$ contains the simple current $L_{\mathfrak{sl}_2}(k',k') \ot L_{\mathfrak{sl}_2}(k'',k'')$ of its affine subalgebra.
\item Replacing $L$ by the other lattice of the genus in the decomposition of $U$ over $\Omega_U \ot V_L$ gives a Haagerup algebra $U^\natural$ of the other case, with Heisenberg commutant $\Omega_U$. The map $U \mapsto U^\natural$ is a bijection between the isomorphism classes of Haagerup algebras of the two cases, and $U$, $U^\natural$ are the only strongly rational vertex operator algebras in the hyperbolic genus of $U$.
\end{enumerate}
\end{thm}

The precise statements are Proposition \ref{prop:associated-lattice} and Theorem \ref{thm:twin}. Part 1 combines \cite[Proposition 2.14]{MoR24a} and the proof of \cite[Corollary 4.10]{MoR24a} with the weight decomposition in Table \ref{tab:lowest-spaces} and a growth estimate for characters. In part 2, $U^\natural$ is defined by the image of the algebra object defining $U$ under a ribbon equivalence $\Rep(\Omega_U) \boxtimes \cal C(L'/L) \to \Rep(\Omega_U) \boxtimes \cal C(L^{\natural\prime}/L^\natural)$ induced by an isometry of discriminant forms. Positivity of $U^\natural$ comes from that of $U$. The independence of the choices is not automatic: the orthogonal group of the lattice of case (I) acts on its discriminant group only by $\pm1$, and we use instead the Weyl group of the currents of the algebra of case (II), which acts through the whole orthogonal group of the discriminant form. The exchange of lattices is the analogue of the construction of \cite[Proposition 5.3]{MoR24a} for holomorphic algebras. What is specific here is its computation for $Z(\cal H_3)$.

\paragraph{Classification and explicit realizations.}
\label{par:external-hypotheses}
We write $\mathrm{(H_J)}$ for the assertion that Jiao's explicit algebra $V\cong S$ is a Haagerup algebra, and $\mathrm{(H_R)}$ for the corresponding assertion about the explicit algebra $U^{\mathrm R}$ of Riesen. These conditions include positivity, strong rationality and the braided equivalence with $Z(\cal H_3)$. The lattice exchange of Theorem \ref{thm:twin} starts from an arbitrary Haagerup algebra and determines the two members of its hyperbolic genus. The uniqueness results then give the classification from either explicit realization. By \cite[Theorem 6.1]{ArhIII} and Theorem \ref{thm:jiao-identification}, $\mathrm{(H_J)}$ holds. Its proof in \cite{ArhIII} uses Theorems \ref{thm:exact-weight-two}, \ref{thm:unitarity} and \ref{thm:strong-generation} and the results of Section \ref{sec:existence}, and Theorem \ref{thm:jiao-identification} transfers it to Jiao's algebra. None of these results depends on $\mathrm{(H_J)}$ (see the beginning of Sections \ref{subsec:case-two-exact} and \ref{sec:existence}). By \cite[Theorem 7.13]{ArhIII}, which uses the identifications of \cite[Proposition 4.1 and Theorem 7.3]{DJX26}, $\mathrm{(H_R)}$ holds as well.

\begin{thm}
\label{thm:intro-classification}
\begin{enumerate}
\item If $\mathrm{(H_J)}$ or $\mathrm{(H_R)}$ holds, there are exactly two Haagerup algebras up to isomorphism, one in each case. They are exchanged by $U\mapsto U^\natural$, have isomorphic Heisenberg commutants and form a single hyperbolic genus.
\item The explicit algebra $S$ is a Haagerup algebra by \cite[Theorem 6.1]{ArhIII}. Consequently the two isomorphism classes are represented by $S$ in case (II), of levels $(3,13)$, and the twin $S^\natural$ in case (I), of levels $(1,39)$.
\item The twin $S^\natural$ is isomorphic as a vertex operator algebra to both explicit algebras, $S^\natural\cong U^{\mathrm H}\cong U^{\mathrm R}$, where $U^{\mathrm H}$ is Huang's algebra and $U^{\mathrm R}$ is Riesen's, by \cite[Theorem 7.9 and Corollary 7.10]{ArhIII}, using the common fixed-point algebra and the Huang--Riesen identification of \cite[Theorem 7.3 and Proposition 4.1]{DJX26}.
\end{enumerate}
\end{thm}

Theorems \ref{thm:twin} and \ref{thm:classification-one} give part 1 and the construction of $S^\natural$ used in part 2. Part 1 of Theorem \ref{thm:classification-one} also gives an embedding of $S^\natural$ in $V'$ whose image contains $\psi_{\mathrm R}(U^{\mathrm R})$ and has the same homogeneous spaces as $\psi_{\mathrm R}(U^{\mathrm R})$ through weight four. The identification of the entire twin in part 3 uses the common fixed point algebra in \cite[Section 7]{ArhIII}. Table \ref{tab:dependencies} records the sources of these results.

\begin{table}[ht]
\centering
\small
\begin{tabular}{p{7.4cm}p{6cm}}
\hline
Statement & Status \\
\hline
$\gamma = 0$, $U_1$, the two embeddings, lowest weight spaces, $U_2$ as a $U_1$-module (Theorems \ref{thm:intro-gamma}, \ref{thm:intro-lie}) & exact computations and the modularity theorem of \cite{EvG11} \\
Case (II): $C_2$, $U_2$, $S$ unitary, strongly generated, $C_2$-cofinite, $S \cong$ Jiao's algebra (Theorem \ref{thm:intro-case-two}) & exact computations, certified exclusions and Proposition \ref{prop:weight-four-certificate} \\
Case (I): identification of $e$, $f$, $v$, $U \supset \psi_{\mathrm R}(U^{\mathrm R})$ & exact and certified computations \\
Dimensions of $\psi_{\mathrm R}(U^{\mathrm R})$ at least $120$, $956$, $6010$ in weights $2$, $3$, $4$ (Proposition \ref{prop:riesen-dims}) & ranks modulo primes of $\mathfrak p$-integral data (rigorous lower bounds) \\
Associated lattice, $\Omega_U$, $J \subset U$, twin, bijection, hyperbolic genus of $U$ (Theorem \ref{thm:intro-genus}) & no external existence hypothesis, using the preceding classification of currents \\
Rational forms of $S$ and $\psi_{\mathrm R}(U^{\mathrm R})$ (Proposition \ref{prop:rational-forms}) & exact computations \\
Strong generation and $C_2$-cofiniteness of the explicit $S$, without an ambient Haagerup algebra & exact certificate over $\K$ (Proposition \ref{prop:weight-four-certificate}) \\
Classification into two isomorphism classes from either explicit realization (Theorem \ref{thm:intro-classification}(1)) & Theorems \ref{thm:twin} and \ref{thm:classification-one} \\
Embedding of $S^\natural$ in $V'$ and agreement with $\psi_{\mathrm R}(U^{\mathrm R})$ through weight four & Theorem \ref{thm:classification-one}(1), with Proposition \ref{prop:riesen-dims} \\
Existence of the realization $S$ and identification of $S^\natural$ with the Huang--Riesen algebra & \cite[Theorems 6.1 and 7.9, Corollary 7.10]{ArhIII}, with the inputs of \cite[Theorem 2.3]{ArhIII} \\
Rationality criterion and modular upper presentations for the Zhu corner of $S$ & Theorems \ref{thm:conditional} and \ref{thm:zhu-modp} \\
Semisimplicity of $A(S)$, the twelve simple modules and the full lattice branching & \cite[Theorem 2.2(F1)--(F3) and Appendix A]{ArhIII} \\
\hline
\end{tabular}
\caption{Dependencies and sources of the main results.}
\label{tab:dependencies}
\end{table}

We also prove that the explicit algebras $S$ and $\psi_{\mathrm R}(U^{\mathrm R})$ have rational forms, given by explicit Galois descent data with values in the centralizer of $U_1$ in $\operatorname{Aut}(V')$ (Proposition \ref{prop:rational-forms}), and that in case (II) the commutant $C$ is not a tensor product of two nontrivial vertex operator algebras of CFT type, although its graded dimensions suggest a coset description (Corollary \ref{cor:commutant-not-tensor}). Section \ref{sec:existence} reduces rationality to the semisimplicity of the Zhu algebra. The inclusion of the affine simple-current extension in $S$ gives a full Morita corner. The selected identities in weights three and four give a modular upper presentation of dimension $762$ with the block dimensions predicted by Table \ref{tab:lowest-spaces} (Theorem \ref{thm:zhu-modp}). The exact descent and the branching are proved in \cite[Appendix A]{ArhIII}, see Section \ref{subsec:what-remains}.

\paragraph{Computations.} Many proofs rest on finite computations: exact rational arithmetic for the vector-valued modular forms and the relations between traces, exact arithmetic over $\K$ in a model of $V'$, Gröbner bases over $\Q$ or $\Q(\sqrt{13})$, after a Galois descent when needed, certified interval arithmetic for a few exclusions, and exact integer computations for lattices. Computations modulo primes are used as evidence, as controls, to choose pivots which are then checked exactly, or as rigorous lower bounds for ranks of matrices with $\mathfrak p$-integral entries, since the rank of such a matrix is at least that of its reduction. Each use is labelled. Appendix \ref{app:computations} distinguishes the mathematical role of each computation from the verification status of the accompanying files. The supplementary material provides a statement-by-statement index of the computational sources and certificates.

\paragraph{Notation.} Throughout, $V' = V_{E_6 \oplus A_2}$, and $U \subset V'$ is a Haagerup algebra with $U_1 = \mathfrak g' \oplus \mathfrak g''$, of levels $(k',k'')$. Moreover $A = L_{\mathfrak g'}(k',0) \ot L_{\mathfrak g''}(k'',0)$ is its affine subalgebra and $C = \operatorname{Com}_U(A)$. We set $\K = \Q(\sqrt{13},\sqrt{-3})$. Furthermore, $S$ is the explicit algebra of Theorem \ref{thm:intro-case-two}, generated by the explicit currents and one of four explicit subspaces of $V'_2$ (Section \ref{subsec:case-two-exact}), and $U^{\mathrm R} = \la A,e,f,v\ra$ that of \cite{Rie26}. In Section \ref{sec:genus}, $\mathfrak h$ is a Cartan subalgebra of $U_1$, $\Omega_U$ the Heisenberg commutant, $L$ the associated lattice and $L^\natural$ the other lattice of its genus.

\paragraph{Structure of the paper.} Section \ref{sec:preliminaries} fixes the modular data of $Z(\cal H_3)$ and the notion of admissible character vector. Section \ref{sec:c8} treats the character vectors at central charge eight, and Section \ref{sec:chain} establishes the chain $U \subset V_{E_6 \oplus A_2} \subset V_{E_8}$ and the self-duality constraints. Section \ref{sec:weight-two} derives the weight-two relations in exact arithmetic and Section \ref{sec:higher-weight} the relations of higher weight. Section \ref{sec:gamma-one} proves Theorem \ref{thm:intro-gamma} and Section \ref{sec:weight-one-lie} proves Theorem \ref{thm:intro-lie}. Section \ref{sec:uniqueness} gives the constraints on the weight-two space. Section \ref{sec:griess} determines the commutant of the currents and $U_2$, identifies the algebras of \cite{Jia26} and \cite{Rie26}, and proves the uniqueness theorems and the rational forms. Section \ref{sec:existence} gives the Zhu-algebra criterion, the modular computations and the references for exact descent and branching. Section \ref{sec:genus} determines the associated lattices, constructs the twins and applies the classification to $S$. Appendix \ref{app:computations} describes the computations.

\section{Preliminaries}
\label{sec:preliminaries}

\subsection{The modular data of the Drinfeld center}

We use the modular data of $Z(\cal H_3)$ in the simplified form of \cite[equations (2.5a) and (2.5b) p.~20]{EvG11}. The simple objects are labelled
\[
0,\ b,\ a,\ c_0,\ c_1,\ c_2,\ d_1,\ldots,d_6,
\]
with $0$ the unit. Their quantum dimensions are
\begin{equation}
\label{eq:center-dimensions}
1,\quad 3d+1,\quad 3d+2,\quad 3d+2,\quad 3d+2,\quad 3d+2,\quad 3d,\ldots,3d,
\end{equation}
and their categorical twists are
\begin{equation}
\label{eq:center-twists}
1,\ 1,\ 1,\ 1,\ \e^{\frac{2\pi\i}{3}},\ \e^{-\frac{2\pi\i}{3}},\ \e^{\frac{2\pi\i 6l^2}{13}}\quad (1 \leq l \leq 6).
\end{equation}
Put $x = \frac{13-3\sqrt{13}}{26}$, $y = \frac{3}{\sqrt{13}}$ and $c(j) = -2y\cos\big(\frac{2\pi j}{13}\big)$. The $S$ matrix is $\frac{1}{3}M$, where the first two rows of $M$ are $(x,1-x,1,1,1,1,y,\ldots,y)$ and $(1-x,x,1,1,1,1,-y,\ldots,-y)$, the block indexed by $(a,c_0,c_1,c_2)$ is
\[
\begin{bmatrix}
2&-1&-1&-1\\
-1&2&-1&-1\\
-1&-1&-1&2\\
-1&-1&2&-1
\end{bmatrix},
\]
the entries $M_{a,0}, M_{a,b}, \ldots, M_{c_2,b}$ are equal to $1$, the entries $M_{d_l,0} = -M_{d_l,b}$ are equal to $y$, the entries between $\{a,c_0,c_1,c_2\}$ and $\{d_1,\ldots,d_6\}$ vanish, and $M_{d_l,d_{l'}} = c(ll')$. The global dimension is $\cal D^2 = 9(3d+2)^2$, since $S_{0,0} = \frac{x}{3} = \frac{1}{3(3d+2)}$.

The modular data of $Z(\cal H_3)$ were computed by Izumi from the tube algebra \cite{Izu01}. The form above is that of \cite{EvG11}, and its identification with the Drinfeld center of the Haagerup category is established in \cite{GSY26}. The charge conjugation is trivial by an exact calculation. Indeed, for $1 \leq a,b \leq 6$, the geometric sum of the thirteenth roots of unity gives
\[
\sum_{l = 1}^6 \cos\Big(\frac{2\pi al}{13}\Big) = -\frac12,\quad
\sum_{l = 1}^6 \cos\Big(\frac{2\pi al}{13}\Big)\cos\Big(\frac{2\pi bl}{13}\Big)
= \frac{13\delta_{a,b}-2}{4}.
\]
Together with $2x-1 = -y$, $x(1-x) = 1/13$ and $y^2 = 9/13$, these identities give $M^2 = 9\Id$, hence $S^2 = \Id$. The displayed matrix is real and symmetric, so it is also unitary. An independent numerical calculation verifies these identities and $(ST)^3 = S^2$ with sixty significant digits. The exact modular relations also follow from the categorical origin of the data. The only invertible simple object is $0$, by \eqref{eq:center-dimensions}.

\begin{remark}\normalfont
\label{rem:labels}
We use the labels of \cite{EvG11}, whose six-object even part is $\cal H_2$ in the convention of \cite{GrS12}. The inductions of the units from $\cal H_2$ and $\cal H_3$ to their common center have underlying objects $I_{\cal H_2}(1) \cong 0 \oplus b \oplus 2a$ and $I_{\cal H_3}(1) \cong 0 \oplus b \oplus 2c_0$, respectively \cite[Section 2.2]{EvG11}, \cite[Corollary 3.8 and Remark 3.7]{ArhNets}. Thus $0,b,a,c_0$ correspond respectively to $1,\pi_1,\sigma_1,\pi_2$ in \cite{AHKZ26,BoS25}. In the untwisted sector of a theory with $\cal H_3$ symmetry, the line $\rho$ acts on the multiplicity spaces of $0$ and $b$ by $d$ and $d_- \ov{\mathrm{def}}{=} \frac{3-\sqrt{13}}{2}$, respectively, and on that of $c_0$ as part of the two-dimensional irreducible representation of the fusion ring \cite[Appendix A]{AHKZ26}.
\end{remark}

\subsection{Vertex operator algebras and character vectors}

A vertex operator algebra $V$ is called strongly rational if it is simple, of CFT type, self-contragredient, rational and $C_2$-cofinite. We call it positive if, in addition, every irreducible $V$-module $M \not\cong V$ has a positive conformal weight $h_M > 0$. For a positive strongly rational $V$ the effective central charge equals the central charge. By \cite{Hua08}, $\Rep(V)$ is a modular tensor category. We denote by $\cal C(K)$ the pointed modular category of the discriminant form of an even positive-definite lattice $K$. It is equivalent to $\Rep(V_K)$.

In the case considered here, a braided equivalence $\Rep(U) \cong Z(\cal H_3)$ preserves the ribbon structure. Indeed, this nondegenerate braided category has no nontrivial invertible simple objects, so its universal grading group is trivial \cite[Theorem 6.3]{GeN08}. Pivotal structures differ by a tensor automorphism of the identity functor \cite{ENO05}, and such automorphisms are given by the characters of the universal grading group. Hence the pivotal structure is unique. We use this structure and the normalization $S_{0,0} > 0$ throughout.

Let $(S,T)$ be the modular data of a modular tensor category with simple objects $\Irr = \{0 = i_0,i_1,\ldots\}$ and conformal weights $h_i \in [0,1)$ defined modulo one by the twists. For $c \in \mathbb{Q}$, put $T_c \ov{\mathrm{def}}{=} \e^{-\frac{2\pi\i c}{24}}T$. We write $q = \e^{2\pi\i\tau}$ for $\tau$ in the upper half-plane $\mathbb H$.

\begin{defi}
\label{def:admissible}
An admissible character vector at central charge $c$ for $(S,T)$ is a family $X = (X_i)_{i \in \Irr}$ of holomorphic functions on $\mathbb H$ such that:
\begin{enumerate}
\item $X(-\frac{1}{\tau}) = SX(\tau)$ and $X(\tau+1) = T_cX(\tau)$ for any $\tau \in \mathbb H$,
\item for any $i \in \Irr$, we have $X_i(\tau) = \sum_{n \geq 0} a_{i,n}q^{h_i+n_i+n-\frac{c}{24}}$ with $a_{i,n} \in \Z_{\geq 0}$ and integers $n_i \geq 0$, the series converging on $\mathbb H$,
\item $X_0(\tau) = q^{-\frac{c}{24}}(1+O(q))$,
\item for any $i \neq 0$, the leading exponent of $X_i$ is strictly larger than $-\frac{c}{24}$.
\end{enumerate}
\end{defi}

Condition 4 is the positivity of the conformal weights of the nonvacuum components. By Zhu's modular invariance theorem \cite{Zhu96} and the identification of the modular transformation of the characters with the categorical $S$ matrix \cite{Hua08}, the following holds.

\begin{lemma}
\label{lem:characters-admissible}
Let $V$ be a positive strongly rational vertex operator algebra of central charge $c$. Then the family $\big(\Tr_{M}q^{L_0-\frac{c}{24}}\big)_{M \in \Irr(\Rep(V))}$ is an admissible character vector at central charge $c$ for the modular data of $\Rep(V)$.
\end{lemma}

\begin{proof}
The convergence and the transformation law are \cite{Zhu96}, with the $S$ matrix identified with the categorical one in \cite{Hua08}. The coefficients are dimensions of graded pieces. The vacuum component starts with $q^{-\frac{c}{24}}$ since $V$ is of CFT type, and positivity gives condition 4.
\end{proof}

\subsection{Lattice character vectors}

Let $L$ be an even positive-definite lattice of rank $n$, with dual $L^*$ and discriminant form $q_L \co L^*/L \to \mathbb{Q}/\Z$, $q_L(\lambda) = \frac{\la\lambda,\lambda\ra}{2}$. The functions
\begin{equation}
\label{eq:lattice-characters}
\Theta_{L,\lambda}(\tau) \ov{\mathrm{def}}{=} \frac{1}{\eta(\tau)^n}\sum_{v \in L+\lambda}q^{\frac{\la v,v\ra}{2}},\quad \lambda \in L^*/L,
\end{equation}
are the characters of the irreducible $V_L$-modules. They satisfy $\Theta_L(-\frac{1}{\tau}) = S^L\Theta_L(\tau)$ with $S^L_{\lambda,\mu} = |L^*/L|^{-\frac12}\e^{-2\pi\i\la\lambda,\mu\ra}$. The vacuum component is
\begin{equation}
\label{eq:lattice-vacuum}
\Theta_{L,0}(\tau) = q^{-\frac{n}{24}}\big(1+(n+|R(L)|)q+O(q^2)\big),
\end{equation}
where $R(L)$ is the set of roots, namely vectors of norm two, of $L$. Indeed, the coefficient of $q$ counts the $n$ oscillators $h(-1)\mathbf 1$ and the exponentials $\e^{v}$ with $\la v,v\ra = 2$.

Let $K$ and $K'$ be even positive-definite lattices of ranks $r$ and $8-r$, and let $\varphi \co K^*/K \to K'^*/K'$ be an isomorphism of groups with $q_{K'} \circ \varphi = -q_K$. Then the graph of $\varphi$ glues $K \oplus K'$ to an even unimodular lattice of rank eight, which is isometric to $E_8$ since $E_8$ is the unique even unimodular positive-definite lattice of rank eight, and $K$, $K'$ are primitive sublattices with $K^\perp = K'$ \cite[Propositions 1.4.1 and 1.6.1]{Nik80}. Conversely, a primitive embedding $K \subset E_8$ gives such a pair with $K' = K^\perp$. Moreover,
\begin{equation}
\label{eq:conjugate-S}
S^{K'}_{\varphi(\lambda),\varphi(\mu)} = \ovl{S^K_{\lambda,\mu}},\quad T^{K'}_{\varphi(\lambda)} = \ovl{T^K_\lambda},
\end{equation}
since the bilinear form of $K'^*/K'$ is the opposite of that of $K^*/K$ under $\varphi$.

\section{Character vectors of the Haagerup double at central charge eight}
\label{sec:c8}

Let $\varrho$ be the representation of $\mathrm{SL}_2(\Z)$ on $\Cbb^{12}$ defined by $\varrho\big(\begin{bmatrix}0&-1\\1&0\end{bmatrix}\big) = S$ and $\varrho\big(\begin{bmatrix}1&1\\0&1\end{bmatrix}\big) = T_8$, where $(S,T)$ are the modular data of $Z(\cal H_3)$. Denote by $\cal M(\varrho)$ the space of holomorphic functions $X \co \mathbb H \to \Cbb^{12}$ with $X(\gamma\tau) = \varrho(\gamma)X(\tau)$ for any $\gamma \in \mathrm{SL}_2(\Z)$, which are meromorphic at the cusp. It is a module over the polynomial ring $\Cbb[J]$, where $J(\tau) = q^{-1}+196884q+\cdots$.

Following \cite[equation (2.7) pp.~21--22]{EvG11}, put
\begin{equation}
\label{eq:Lambda}
\Lambda \ov{\mathrm{def}}{=} \diag\Big(\frac{2}{3},\frac{2}{3},\frac{2}{3},\frac{2}{3},0,\frac{1}{3},\frac{5}{39},\frac{20}{39},\frac{32}{39},\frac{2}{39},\frac{8}{39},\frac{11}{39}\Big),
\end{equation}
in the order of the labels, and let $\chi$ be the integral $12 \times 12$ matrix of \cite[equation (2.7)]{EvG11}. Note that $\e^{2\pi\i\Lambda} = T_8$. Let $F = (J-240)\frac{\Delta}{E_{10}} = \sum_{k \geq 0}f_kq^k$ and $G = \frac{\Delta}{E_{10}} = \sum_{k \geq 1}g_kq^k$, where $\Delta$ is the discriminant and $E_{10} = E_4E_6$. One has $f_0 = 1$, $f_1 = 0$ and $g_1 = 1$. Define rational matrices $\Xi_n$, $n \geq -1$, by $\Xi_{-1} = \Id$, $\Xi_0 = \chi$ and, for any $n \geq 1$, by the entrywise formula
\begin{equation}
\label{eq:recursion}
(\Lambda_{ii}-\Lambda_{jj}+n+1)(\Xi_n)_{ij}
 = \bigg(\sum_{l = -1}^{n-1}\Xi_l\big(f_{n-l}(\Lambda-\Id)+g_{n-l}(\chi+[\Lambda,\chi])\big)\bigg)_{ij}.
\end{equation}
The coefficient on the left is nonzero since $|\Lambda_{ii}-\Lambda_{jj}| < 1$. Put
\begin{equation}
\label{eq:Xi}
\Xi(\tau) \ov{\mathrm{def}}{=} q^{\Lambda}\Big(q^{-1}\Id+\sum_{n \geq 0}\Xi_nq^n\Big).
\end{equation}
The recursion \eqref{eq:recursion} expresses the differential equation $q\frac{\d\Xi}{\d q} = \Xi\big((\Lambda-\Id)F+(\chi+[\Lambda,\chi])G\big)$ coefficientwise. Indeed, the coefficient of $q^{\Lambda+n}$ on the left is $(\Lambda+n)\Xi_n$, and the term $l = n$ on the right is $\Xi_n(\Lambda-\Id)$ since $f_0 = 1$ and $g_0 = 0$. This is the differential equation of \cite[Section 3]{BaG07}. The conditions $f_1 = 0$ and $g_1 = 1$ show that the case $n = 0$ is automatically satisfied.

\begin{thm}
\label{thm:EG}
The following hold.
\begin{enumerate}
\item The columns of $\Xi$ belong to $\cal M(\varrho)$ and form a basis of $\cal M(\varrho)$ over $\Cbb[J]$.
\item Every admissible character vector at central charge eight for the modular data of $Z(\cal H_3)$ is one of the two vectors
\begin{equation}
\label{eq:EG-vectors}
X^{(\gamma)} = \Xi_{\bullet,0}+\gamma\Xi_{\bullet,d_3},\quad \gamma \in \{0,1\},
\end{equation}
where $\Xi_{\bullet,j}$ denotes the column of index $j$. In particular, any admissible character vector satisfies
\begin{equation}
\label{eq:EG-vacuum}
X_0(\tau) = q^{-\frac13}\big(1+(6+13\gamma)q+(120+78\gamma)q^2+(956+351\gamma)q^3+O(q^4)\big),\quad \gamma \in \{0,1\}.
\end{equation}
\end{enumerate}
\end{thm}

The second assertion is contained in \cite[Section 2.4 p.~22]{EvG11}. We give the argument in the form needed below, together with an independent verification of the data. We do not use, and do not prove, the converse: the coefficients of $X^{(0)}$ and $X^{(1)}$ are nonnegative integers up to $q^{60}$, but we have no proof that all of them are.

\begin{proof}
1. Assume that the columns of $\Xi$ are vector-valued modular functions. Since $\Xi_{-1} = \Id$, they are linearly independent over $\Cbb(J)$. By \cite[equation (2.22) p.~6]{BaG07}, the $\Cbb[J]$-span of the columns of such a matrix is the whole module $\cal M(\varrho)$ if and only if
\begin{equation}
\label{eq:BG-trace}
\Tr(\Lambda) = \frac{5\cdot 12}{12}+\frac{1}{4}\Tr(S)+\frac{2}{3\sqrt 3}\Re\big(\e^{-\frac{\pi\i}{6}}\Tr(ST_8^{-1})\big).
\end{equation}
The left-hand side is $\frac{8}{3}+\frac{1}{3}+\frac{78}{39} = 5$. On the right-hand side, $\Tr(S) = 0$ and $\Tr(ST_8^{-1}) = -\frac{3}{2}+\frac{3\sqrt 3}{2}\i$, so that $\Re\big(\e^{-\frac{\pi\i}{6}}\Tr(ST_8^{-1})\big) = -\frac{3\sqrt 3}{4}+\frac{3\sqrt 3}{4} = 0$. For $\Tr(S)$, the diagonal entries give $\frac{1}{3}\big(2x+2+\sum_{l = 1}^6c(l^2)\big)$, and the quadratic Gauss sum $\sum_{l = 0}^{12}\e^{\frac{2\pi\i l^2}{13}} = \sqrt{13}$ gives $\sum_{l = 1}^6c(l^2) = -3+\frac{3}{\sqrt{13}} = -(2x+2)$. Hence the right-hand side of \eqref{eq:BG-trace} equals $5$. 

The modularity of the columns is the content of \cite[Section 2.4]{EvG11}. Together with \eqref{eq:BG-trace}, it proves the assertion.

2. Let $X$ be an admissible character vector. By 1, we can write $X = \sum_j\Xi_{\bullet,j}P_j(J)$ with polynomials $P_j$. Let $m$ be the maximal degree of the $P_j$ and suppose that $m \geq 1$. Choose $j$ with $\deg P_j = m$. The component $X_j$ contains the term $\lambda_mq^{\Lambda_{jj}-1-m}$, where $\lambda_m \neq 0$ is the leading coefficient of $P_j$, and no other term of $X$ has this exponent in the component $j$: the other columns contribute exponents at least $\Lambda_{jj}-m$. Since $\Lambda_{jj} \leq \frac{32}{39}$, we obtain the exponent $\Lambda_{jj}-1-m \leq \frac{32}{39}-2 < -\frac{1}{3}$, which contradicts conditions 3 and 4 of Definition \ref{def:admissible}. Consequently, the $P_j$ are constants $v_j$.

Now, the column $j$ contributes the exponent $\Lambda_{jj}-1$ to the component $j$, and the other columns contribute exponents at least $\Lambda_{jj}$. Condition 3 forces $v_0 = 1$, since $\Lambda_{00}-1 = -\frac13$. For $j \neq 0$, condition 4 forces $v_j = 0$ whenever $\Lambda_{jj}-1 \leq -\frac{1}{3}$, that is, whenever $\Lambda_{jj} \leq \frac{2}{3}$. Among the indices $j \neq 0$, only $j = d_3$ satisfies $\Lambda_{jj} > \frac{2}{3}$. Hence $X = \Xi_{\bullet,0}+\gamma\Xi_{\bullet,d_3}$ with $\gamma = v_{d_3}$. The leading coefficient of $X_{d_3}$ is $\gamma$, so that $\gamma \in \Z_{\geq 0}$. The leading coefficient of $X_{d_4}$ is
\[
(\Xi_0)_{d_4,0}+\gamma(\Xi_0)_{d_4,d_3} = \chi_{d_4,0}+\gamma\chi_{d_4,d_3} = 5-3\gamma,
\]
which is nonnegative only for $\gamma \leq 1$. The expansion \eqref{eq:EG-vacuum} is read from the matrices $\Xi_0 = \chi$, $\Xi_1$ and $\Xi_2$ computed by this recurrence. This proves 2.
\end{proof}

\begin{remark}\normalfont
\label{rem:EG-coefficients}
The first coefficients of all components produced by our recursion coincide with those printed in \cite[p.~22]{EvG11}. For instance, the components $d_1$, $d_4$, $d_5$ and $d_6$ start with $7-\gamma$, $5-3\gamma$, $13-5\gamma$ and $14+7\gamma$. The vector $\gamma = 0$ has $\dim U_1 = 6$, $\dim U_2 = 120$, $\dim U_3 = 956$ and $\dim U_4 = 6010$ for a realizing algebra $U$. The last value is read from $\Xi_3$ and is used in Section \ref{sec:griess}. These numbers agree with the algebra announced in \cite{Rie26,Hua26}. The vector $\gamma = 1$ has $19$ currents. By Theorem \ref{thm:gamma-one}, it is not the character vector of a positive strongly rational vertex operator algebra $U$ with $\Rep(U) \cong Z(\cal H_3)$ as braided categories.
\end{remark}

\begin{remark}\normalfont
The modularity of $\Xi$ is imported from \cite[Section 2.4]{EvG11}. As independent controls, we compute its first $60$ coefficients in rational arithmetic, check their integrality and test $\Xi(-1/\tau) = S\Xi(\tau)$ at five points with fifty-digit precision. The largest recorded relative residual is $1.4\times10^{-50}$. The discarded tails are not certified. An exact rational calculation checks \eqref{eq:BG-trace}. The floating-point tests do not enter the proof.
\end{remark}

\section{The chain \texorpdfstring{$U \subset V_{E_6 \oplus A_2} \subset V_{E_8}$}{U < V(E6+A2) < V(E8)}}
\label{sec:chain}

Throughout this section, $U$ is a positive strongly rational vertex operator algebra of central charge eight together with a braided equivalence $\Rep(U) \cong Z(\cal H_3)$. For $z \in \Irr(Z(\cal H_3))$ we denote by $M_z$ the corresponding simple $U$-module and by $R_z$ its lowest weight space. By Lemma \ref{lem:characters-admissible} and Theorem \ref{thm:EG}, the character vector of $U$ is $X^{(\gamma)}$ for some $\gamma \in \{0,1\}$. Reading Theorem \ref{thm:EG}, the lowest conformal weights and the dimensions of the lowest weight spaces are
\begin{equation}
\label{eq:lowest-spaces}
\begin{array}{c}
\begin{array}{c|cccccc}
z & 0 & b & a & c_0 & c_1 & c_2\\
\hline
h_z & 0 & 1 & 1 & 1 & \frac13 & \frac23\\
\dim R_z & 1 & 80-13\gamma & 81 & 81 & 3 & 27
\end{array}
\\[3mm]
\begin{array}{c|cccccc}
z & d_1 & d_2 & d_3 & d_4 & d_5 & d_6\\
\hline
h_z & \frac{6}{13} & \frac{11}{13} & \frac{2}{13}+\epsi & \frac{5}{13} & \frac{7}{13} & \frac{8}{13}\\
\dim R_z & 7-\gamma & 42+16\gamma & \gamma+119(1-\gamma) & 5-3\gamma & 13-5\gamma & 14+7\gamma
\end{array}
\end{array}
\end{equation}
where $\epsi = 0$ if $\gamma = 1$ and $\epsi = 1$ if $\gamma = 0$. Moreover $\dim U_1 = 6+13\gamma$, $\dim U_2 = 120+78\gamma$ and $\dim U_3 = 956+351\gamma$.

\subsection{The holomorphic extension and the \texorpdfstring{$\Z_3$}{Z3}-orbifold}

\begin{prop}
\label{prop:chain}
Let $U$ be as above.
\begin{enumerate}
\item The Lagrangian algebra $I_{\cal H_2}(1) \cong 0 \oplus b \oplus 2a$ defines a holomorphic extension $H = U \oplus M_b \oplus M_a \oplus M_a$ of $U$, and $H \cong V_{E_8}$.
\item The subspace $V' \ov{\mathrm{def}}{=} U \oplus M_b$ is a vertex subalgebra of $H$, and there is an automorphism $g$ of order three of $H$ with $H^g = V'$. The two other eigenspaces of $g$ are isomorphic to $M_a$ as $U$-modules.
\item There is an isomorphism $V' \cong V_{E_6 \oplus A_2}$ which maps $V'_1$ onto $\mathfrak e_6 \oplus \mathfrak{sl}_3$.
\end{enumerate}
\end{prop}

\begin{proof}
1. The Lagrangian algebra $I_{\cal H_2}(1)$ is a connected étale algebra of $\Rep(U)$ with trivial category of local modules. By \cite{HKL15,CKM22}, it defines a simple extension $H$ of $U$. The nonvacuum $U$-summands of $H$ have lowest weight one by \eqref{eq:lowest-spaces}, so $H$ is of CFT type. The operator $L(1)$ preserves these summands and kills their weight-one spaces because their weight-zero spaces vanish. It also kills $U_1$ by self-contragredience of $U$. Thus $L(1)H_1=0$ and $H$ is self-contragredient \cite{Li94}. The arguments for $C_2$-cofiniteness and rationality in the proof of Theorem \ref{thm:twin}, with $U$ in place of $\Omega_U \ot V_{L^\natural}$, give these properties for $H$. The correspondence with local modules then gives $\Rep(H) \cong \mathrm{Vec}$. A strongly rational holomorphic vertex operator algebra of central charge eight is isomorphic to $V_{E_8}$ \cite[Theorem 1 p.~2]{DoM04b}.

2. By \cite[Theorem 4.10 p.~18]{DMNO13}, the subalgebras of $I_{\cal H_2}(1)$ correspond to the fusion subcategories of $\cal H_2$, and the subcategory $\mathrm{Vec}_{\Z_3}$ generated by the invertible objects of $\cal H_2$ corresponds to a connected étale subalgebra $A$ of Frobenius--Perron dimension $\frac{9d+6}{3} = 3d+2$. The only subobject of $0 \oplus b \oplus 2a$ containing $0$ with this dimension is $0 \oplus b$, by \eqref{eq:center-dimensions}. By \cite{HKL15} and the arguments recalled in part 1, $V' = U \oplus M_b$ is a strongly rational vertex subalgebra of $H$ with $\Rep(V') \cong Z(\cal H_3)^0_A$. The latter has Frobenius--Perron dimension $\frac{9(3d+2)^2}{(3d+2)^2} = 9$ \cite[Corollary 3.32 p.~15]{DMNO13}. It is Witt equivalent to $Z(\cal H_3)$ \cite[Proposition 5.4 p.~22]{DMNO13}, hence of the form $Z(\cal D)$ for a fusion category $\cal D$ \cite[Proposition 5.8 p.~23]{DMNO13}, with $\operatorname{FPdim}(\cal D) = 3$. A fusion category of prime Frobenius--Perron dimension is pointed \cite[Corollary 8.30]{ENO05}. Here $\cal D \cong \mathrm{Vec}_{\Z_3}^{\omega}$ for some $3$-cocycle $\omega$. Its center is pointed: every centralizer is cyclic, and every irreducible projective representation of $\Z_3$ is one-dimensional, since $H^2(\Z_3,\Cbb^\times) = 0$. Thus $\Rep(V')$ is pointed with nine simple objects, all simple currents.

For every nonvacuum simple $V'$-module $N$, its restriction $R(N)$ has Frobenius--Perron dimension $3d+2$. Free-module adjunction gives
\[
\operatorname{Hom}_U(U,R(N)) \cong \operatorname{Hom}_{V'}(V',N) = 0.
\]
The list \eqref{eq:center-dimensions} therefore forces $R(N)$ to be simple. As a $U$-module, the complement of $V'$ in $H$ is $M_a \oplus M_a$. Its decomposition into simple $V'$-modules consequently gives $H = V' \oplus N_1 \oplus N_2$ with simple currents $N_1, N_2$ of $V'$, both isomorphic to $M_a$ as $U$-modules. These simple currents are distinct and, together with the vacuum, form a group of order three under fusion. To see this directly, regard $H$ as a connected commutative separable algebra in the pointed category $\Rep(V')$. Its regular right module is simple, since $\operatorname{End}_H(H) \cong \operatorname{Hom}(\mathbf 1,H) = \Cbb$. For an invertible simple object $J$, the right module $J \ot H$ is also simple, and free-module adjunction gives $\operatorname{Hom}(J,H) \cong \operatorname{Hom}_H(J \ot H,H)$. Thus the multiplicity of $J$ in $H$ is either zero or one. The objects for which this multiplicity is one form the stabilizer of the simple right module $H$ under the action $J \mapsto J \ot H$, hence a subgroup. Since $\operatorname{FPdim}_{\Rep(V')}(H) = 3$, this subgroup has order three. The operator acting by $e^{\frac{2\pi\i k}{3}}$ on $N_k$ is therefore an automorphism $g$ of order three with $H^g = V'$.

3. The automorphism group of $V_{E_8}$ is the complex Lie group of type $E_8$ \cite[Theorem 2.1]{DoN99}. An element of order three is conjugate to $e^{2\pi\i\beta(0)}$ with $\beta$ in the Cartan subalgebra, and its conjugacy class is given by Kac coordinates \cite[Chapter 8]{Kac90}. There are four classes of elements of order three, whose centralizers in $\mathfrak e_8$ have dimensions $80$, $86$, $92$ and $134$, of types $A_8$, $E_6 \oplus A_2$, $D_7 \oplus \Cbb$ and $E_7 \oplus \Cbb$, as follows by enumerating the Kac coordinates. Now $\dim V'_1 = (6+13\gamma)+(80-13\gamma) = 86$ by \eqref{eq:lowest-spaces}, so that $g$ is of type $E_6 \oplus A_2$. For $g = e^{2\pi\i\beta(0)}$, the fixed point algebra is the lattice vertex operator algebra of $\{\alpha \in E_8 : \langle \alpha,\beta\rangle \in \Z\}$. This sublattice contains the root lattice $E_6 \oplus A_2$, which has index three in $E_8$, and is a proper sublattice. Hence it equals $E_6 \oplus A_2$.
\end{proof}

\begin{remark}\normalfont
\label{rem:chain-literature}
The realization of $V_{E_6 \oplus A_2}$ as a $\Z_3$-orbifold of $V_{E_8}$, and the expectation that a Haagerup vertex operator algebra lies inside it, go back to \cite[Section 2.4 p.~22]{EvG11}. The physical picture of a hypergroup acting on the $\Z_3$-invariant algebra is in \cite[Section 4.3 pp.~83--88]{GaR26}, and the announced algebra of \cite{Rie26,Hua26} is constructed inside $V_{E_6 \oplus A_2}$. Proposition \ref{prop:chain} shows that every positive strongly rational realization, for either value of $\gamma$, sits in this chain.
\end{remark}

\subsection{Restriction of the lattice modules}

The discriminant group of $E_6 \oplus A_2$ is $\Z_3 \times \Z_3$, with quadratic form $q(i,j) = \frac{2i^2+j^2}{3} \mod 1$. We denote by $V_{(i,j)}$ the corresponding simple $V_{E_6 \oplus A_2}$-modules. Both lines $\{(0,0),(1,1),(2,2)\}$ and $\{(0,0),(1,2),(2,1)\}$ are isotropic. We choose coordinates adapted to the extension $H$ of Proposition \ref{prop:chain}, so that $H$ is the lattice extension attached to $\{(0,0),(1,1),(2,2)\}$. The diagram automorphism of $A_2$ exchanges the two choices. The coordinates in this subsection are adapted to $H$. In Section \ref{subsec:categorical-route}, we use the fixed coordinates of the explicit model of $S$. The lowest weight space of $V_{(i,j)}$ is, as a module over $\mathfrak e_6 \oplus \mathfrak{sl}_3$, the representation $(\mathbf{27}^{i},\mathbf 3^{j})$, where $\mathbf{27}^{0} = \mathbf 3^0$ is the trivial module, $\mathbf{27}^{1} = \mathbf{27}$, $\mathbf{27}^{2} = \overline{\mathbf{27}}$, and similarly for $\mathbf 3$.

\begin{lemma}
\label{lem:restriction}
Use the coordinates adapted to $H$ and let $(i,j) \neq (0,0)$. The restriction of $V_{(i,j)}$ to $U$ is simple. It is isomorphic to $M_{c_1}$ if $i = 0$, to $M_{c_2}$ if $j = 0$, to $M_a$ for the two nonzero classes $(1,1)$, $(2,2)$ of $E_8/(E_6 \oplus A_2)$, and to $M_{c_0}$ for the two remaining classes $(1,2)$, $(2,1)$.
\end{lemma}

\begin{proof}
The Frobenius--Perron dimension and free-module adjunction argument in the proof of Proposition \ref{prop:chain}, part 2, shows that the restriction $R(w)$ of any nonvacuum simple $V'$-module is one of $a$, $c_0$, $c_1$, $c_2$. The twist of $R(w)$ equals $e^{2\pi\i q(i,j)}$. The twists of $a$ and $c_0$ are one, those of $c_1$ and $c_2$ are $e^{\frac{2\pi\i}{3}}$ and $e^{-\frac{2\pi\i}{3}}$ by \eqref{eq:center-twists}, while $q(0,\pm 1) = \frac13$, $q(\pm 1,0) = \frac23$ and $q(\pm 1,\pm 1) = 0$. The classes of $E_8$ give $M_a$ by Proposition \ref{prop:chain}. It remains to show that the classes $(1,2)$ and $(2,1)$ give $M_{c_0}$. For a simple object $z$ of $Z(\cal H_3)$ and $w \in \Z_3 \times \Z_3$, let $n_{z,w} \in \{0,1\}$ be the multiplicity of $M_z$ in the restriction of $V_w$. For $v \in U$ we have $Z_{V_w}(v,\tau) = \sum_z n_{z,w}Z_{M_z}(v,\tau)$. The one-point functions of pairwise non-isomorphic irreducible modules are linearly independent, and they transform under $\tau \mapsto -\frac{1}{\tau}$ through the $S$ matrix of the representation category \cite[Section 5]{Zhu96}, \cite{Hua08}. Applying this to $U$ and to $V'$, which have the same conformal vector, and comparing gives $S n = n S^{L}$, where $S^L$ is the $S$ matrix of $V' = V_L$, $L = E_6 \oplus A_2$, given after \eqref{eq:lattice-characters}. All the entries $S^L_{w,(0,0)}$ equal $\frac13$. Since $V_{(0,0)} = V'$ restricts to $U \oplus M_b$ (Proposition \ref{prop:chain}), the entry $(c_0,(0,0))$ of this identity reads $S_{c_0,0}+S_{c_0,b} = \frac13 \#\{w : n_{c_0,w} = 1\}$. The left-hand side is $\frac13+\frac13 = \frac23$ by the formula for $S$ in Section \ref{sec:preliminaries}. Hence exactly two classes restrict to $M_{c_0}$, and these are $(1,2)$ and $(2,1)$.
\end{proof}

We have checked Lemma \ref{lem:restriction} on characters: for both values of $\gamma$, $\Theta_{(0,0)} = X_0+X_b$, $\Theta_{(\pm 1,\pm 1)} = X_a = X_{c_0}$, $\Theta_{(0,\pm 1)} = X_{c_1}$ and $\Theta_{(\pm 1,0)} = X_{c_2}$ up to order $q^4$. The characters do not distinguish $a$ from $c_0$. The one-point functions do, as in the proof above.

\subsection{Self-duality constraints on the currents}

Let $\mathfrak s \ov{\mathrm{def}}{=} U_1 \subset \mathfrak e_6 \oplus \mathfrak{sl}_3$, and let $p_1$ and $p_2$ be the projections onto $\mathfrak e_6$ and $\mathfrak{sl}_3$. The Lie algebra $\mathfrak s$ is reductive \cite[Theorem 1 p.~2]{DoM04}. For a simple ideal $\mathfrak g$ of $\mathfrak s$ we denote by $i_1(\mathfrak g)$ and $i_2(\mathfrak g)$ the Dynkin indices of $p_1|_{\mathfrak g}$ and $p_2|_{\mathfrak g}$ (zero if the projection vanishes). Since $\mathfrak e_6$ and $\mathfrak{sl}_3$ have index one in $\mathfrak e_8$, the level of the affine algebra generated by $\mathfrak g$ in $U$ is $k(\mathfrak g) = i_1(\mathfrak g)+i_2(\mathfrak g)$.

\begin{prop}
\label{prop:self-dual}
\begin{enumerate}
\item The representations $\mathbf 3$ and $\overline{\mathbf 3}$ of $p_2(\mathfrak s)$ are isomorphic, and the representations $\mathbf{27}$ and $\overline{\mathbf{27}}$ of $p_1(\mathfrak s)$ are isomorphic.
\item The Lie algebra $p_2(\mathfrak s)$ is zero, a one-dimensional torus with weights $(t,-t,0)$ on $\mathbf 3$, a regular $\mathfrak{sl}_2$ or the principal $\mathfrak{so}_3$. In particular, a simple ideal $\mathfrak g$ of $\mathfrak s$ with $p_2(\mathfrak g) \neq 0$ is isomorphic to $\mathfrak{sl}_2$ and $i_2(\mathfrak g) \in \{1,4\}$.
\item The rank of $p_1(\mathfrak s)$ is at most four.
\end{enumerate}
\end{prop}

\begin{proof}
1. By Lemma \ref{lem:restriction}, $V_{(0,1)}$ and $V_{(0,2)}$ are isomorphic $U$-modules. An isomorphism maps lowest weight spaces to lowest weight spaces and commutes with the zero modes $x(0)$, $x \in U_1$. The lowest weight spaces are $(\mathbf 1,\mathbf 3)$ and $(\mathbf 1,\overline{\mathbf 3})$, on which $x \in \mathfrak s$ acts through $p_2(x)$. The same argument with $V_{(1,0)}$ and $V_{(2,0)}$ gives the second assertion.

2. The action of $\mathfrak s$ on an ordinary $U$-module is completely reducible \cite[Theorem 1 p.~2]{DoM04}. By Lemma \ref{lem:restriction}, this applies in particular to $\mathbf3$ and $\mathbf{27}$. Thus the centers of $p_2(\mathfrak s)$ and $p_1(\mathfrak s)$ act semisimply, and their Cartan subalgebras are toral in the ambient simple Lie algebras. The reductive subalgebras of $\mathfrak{sl}_3$ are $0$, tori, $\mathfrak{sl}_2$ (regular, $\mathbf 3 = \mathbf 2 \oplus \mathbf 1$), $\mathfrak{gl}_2$, $\mathfrak{so}_3$ (principal, $\mathbf 3$ irreducible) and $\mathfrak{sl}_3$. Self-duality of $\mathbf 3$ excludes $\mathfrak{gl}_2$, $\mathfrak{sl}_3$ and every torus whose weights on $\mathbf 3$ do not form a symmetric set, which leaves the tori with weights $(t,-t,0)$. The regular $\mathfrak{sl}_2$ and the principal $\mathfrak{so}_3$ have Dynkin indices one and four.

3. Choose a Cartan subalgebra $\mathfrak h$ of $\mathfrak e_6$ containing a Cartan subalgebra of $p_1(\mathfrak s)$, and denote by $\Omega$ the set of weights of $\mathbf{27}$. This representation is minuscule: its $27$ weights form the Weyl orbit of the fundamental weight $\omega_1$, each with multiplicity one. We first give a proof without enumeration. If the rank were at least five, choose a five-dimensional subspace $\mathfrak t$ of the Cartan subalgebra of $p_1(\mathfrak s)$. By 1, the multiset of restrictions of $\Omega$ to $\mathfrak t$ is invariant under negation. Its cardinality is odd, so at least one weight $\lambda \in \Omega$ vanishes on $\mathfrak t$. Since no weight in $\Omega$ is zero, the annihilator of $\mathfrak t$ in $\mathfrak h^*$ is $\Cbb\lambda$. After applying a Weyl transformation, we may assume $\lambda = \omega_1$.

Let $P$ and $Q$ be the weight and root lattices of $E_6$, respectively. In the Bourbaki numbering, the Cartan matrix gives
\[
\omega_1 = \frac{4\alpha_1+3\alpha_2+5\alpha_3+6\alpha_4+4\alpha_5+2\alpha_6}{3},\quad
\langle\omega_1,\omega_1\rangle = \frac43,\quad
\langle\omega_1,\alpha_1\rangle = 1.
\]
The weight $\mu = s_{\alpha_1}(\omega_1) = \omega_1-\alpha_1$ belongs to $\Omega$. Symmetry of the restrictions gives $\nu \in \Omega$ and $c \in \Cbb$ with $\mu+\nu = c\omega_1$. All weights in $\Omega$ have squared norm $4/3$, while $\langle\mu,\omega_1\rangle = 1/3$. Hence
\[
0 = \langle c\omega_1-\mu,c\omega_1-\mu\rangle-\langle\mu,\mu\rangle
= \frac43c^2-\frac23c,
\]
so $c = 0$ or $c = 1/2$. The first case is impossible because every weight in $\Omega$ lies in $\omega_1+Q$, whereas $2\omega_1 \notin Q$. The second case is impossible because $\mu+\nu \in P$ and $\omega_1/2 \notin P$. This contradiction proves the rank bound.

For comparison, the finite calculation can be retained independently. If two multisets of ambient weights coincide modulo a line while differing before restriction, that line is spanned by a nonzero sum $w+w'$ of weights in $\Omega$. In fundamental-weight coordinates, clear denominators, divide by the greatest common divisor and choose the first nonzero coordinate positive. This gives $243$ lines. For each representative $v$, an exact comparison of the multisets modulo $\Q v$ shows that none is symmetric under negation. This calculation is not needed for the preceding proof.
\end{proof}

\begin{remark}\normalfont
\label{rem:announced-fits}
The announced algebra of \cite{Rie26,Hua26} has $U_1 \cong \mathfrak{sl}_2 \oplus \mathfrak{sl}_2$ at levels $1$ and $39$. This is compatible with Proposition \ref{prop:self-dual}: $\mathfrak e_6 \supset \mathfrak{sl}_6 \oplus \mathfrak{sl}_2$ with $\mathbf{27} = (\overline{\mathbf 6},\mathbf 2) \oplus (\mathbf{15},\mathbf 1)$. The factor $\mathfrak{sl}_2$ has index one, and the diagonal of the principal $\mathfrak{sl}_2$ of $\mathfrak{sl}_6$ (index $35$) and the principal $\mathfrak{so}_3$ of $\mathfrak{sl}_3$ (index $4$) has level $39$. All the restrictions are self-dual. This is case (I) of Theorem \ref{thm:weight-one-classification}. We do not use this description before Section \ref{sec:weight-one-lie}.
\end{remark}

\section{Weight-two relations}
\label{sec:weight-two}

We now combine the modularity of one-point functions with the explicit lowest weight spaces given by the chain. This is a vector-valued version of the trace identities of Schellekens \cite{Sch93} and Dong--Mason \cite{DoM04b}. In the context of the Haagerup modular data, the weight-two one-point functions of currents were introduced by Gannon \cite[Section 2.4]{Gan19}, who called them the Killing vector and used them, under the assumptions that $U$ is unitary and contained in $V_{E_6 \oplus A_2}$, to begin the analysis of the case of six currents \cite[Section 3.3]{Gan19}. The novelty below is that the relations are obtained in exact arithmetic, for both character vectors, without unitarity.

\begin{lemma}
\label{lem:zhu-trace}
Let $V$ be a strongly rational vertex operator algebra of central charge $c$ and let $h \in V_1$ with $L(1)h = 0$, contained in a Cartan subalgebra of $V_1$. For every ordinary $V$-module $M$ we have
\begin{equation}
\label{eq:zhu-trace}
Z_M(h[-1]h,\tau) = \operatorname{Tr}_M h(0)^2q^{L(0)-\frac{c}{24}}-\frac{\langle h,h\rangle}{12}E_2(\tau)\operatorname{Tr}_M q^{L(0)-\frac{c}{24}},
\end{equation}
where $h(1)h = \langle h,h\rangle\mathbf 1$. The vector $(Z_{M_z}(h[-1]h,\tau))_z$, indexed by the simple modules, is a weakly holomorphic vector-valued modular form of weight two for the representation of $\mathrm{SL}_2(\Z)$ on the characters.
\end{lemma}

\begin{proof}
From $Y[h,z] = Y(h,e^z-1)e^z$ and the expansions $\frac{e^z}{e^z-1} = \frac1z+\frac12+O(z)$, $\frac{e^z}{(e^z-1)^2} = \frac{1}{z^2}-\frac{1}{12}+O(z)$, we obtain $h[-1]h = h(-1)h+\frac12h(0)h-\frac{1}{12}h(1)h = h(-1)h-\frac{\langle h,h\rangle}{12}\mathbf 1$, since $h(0)h = 0$ and $h(n)h = 0$ for $n \geq 2$. The zero mode of $h(-1)h$ is $h(0)^2+2\sum_{n \geq 1}h(-n)h(n)$. Put $\kappa = \langle h,h\rangle$ and $Q = q^{L(0)-c/24}$. For $m \geq 1$, graded cyclicity of the trace and $[h(m),h(-m)] = m\kappa\Id$ give
\[
A_m \ov{\mathrm{def}}{=} \operatorname{Tr}_M h(-m)h(m)Q
= q^m\big(A_m+m\kappa\operatorname{Tr}_M Q\big).
\]
Thus $A_m = \frac{m\kappa q^m}{1-q^m}\operatorname{Tr}_M Q$. Substituting into the zero-mode formula and using $\sum_{m \geq 1}mq^m/(1-q^m) = (1-E_2)/24$ proves \eqref{eq:zhu-trace}. This is first an identity of graded series and then of convergent functions. It also covers $\kappa = 0$. The vector $h[-1]h$ has weight two for $L[0]$, and the last assertion is Zhu's theorem \cite{Zhu96} together with the identification of the representation \cite{Hua08}.
\end{proof}

\subsection{The obstruction principle}

We denote by $\varrho^* \ov{\mathrm{def}}{=} (\varrho^{\mathrm T})^{-1}$ the dual representation. For the matrices of Section \ref{sec:preliminaries}, we have $S^{\mathrm T} = S$, $S^2 = \Id$ and $T_8^* = T_8^{-1}$. Since $T_8$ is diagonal, it follows that $\varrho^*(S) = S$ and $\varrho^*(T) = \ovl{T_8}$. If $W$ is a function on $\mathbb H$ with values in $\Cbb^{12}$ and $G$ is another such function, we write $W \cdot G = \sum_z W_zG_z$. We write $D = q\frac{\d}{\d q}$, and $\vartheta_k = D-\frac{k}{12}E_2$ for the Serre derivative in weight $k$. Recall that $\vartheta_k$ maps weakly holomorphic vector-valued modular forms of weight $k$ for $\varrho$ to weakly holomorphic forms of weight $k+2$ for $\varrho$, and that $DJ = -\frac{E_4^2E_6}{\Delta}$ is a weakly holomorphic modular form of weight two.

\begin{lemma}
\label{lem:residue}
Let $k \in 2\Z$. Let $W$ be a weakly holomorphic vector-valued modular form of weight $k$ for $\varrho$, and let $G$ be a weakly holomorphic vector-valued modular form of weight $2-k$ for $\varrho^*$. Then the constant term of the $q$-expansion of $W \cdot G$ vanishes.
\end{lemma}

\begin{proof}
For $\gamma = \begin{bmatrix}a&b\\c&d\end{bmatrix} \in \mathrm{SL}_2(\Z)$, we have
\[
(W \cdot G)(\gamma\tau) = (c\tau+d)^{k}(c\tau+d)^{2-k}\big(\varrho(\gamma)W(\tau)\big)^{\mathrm T}\varrho(\gamma)^{-\mathrm T}G(\tau) = (c\tau+d)^2(W \cdot G)(\tau).
\]
Hence $f = W \cdot G$ is a weakly holomorphic modular form of weight two for $\mathrm{SL}_2(\Z)$, and the differential $f(\tau)\d\tau$ descends to a meromorphic differential on the compactified modular curve, holomorphic outside the cusp. At the elliptic points this uses the invariance of $f(\tau)\d\tau$: in a local coordinate $s$ at $\tau = \i$ with $\tau \mapsto -\frac{1}{\tau}$ acting by $s \mapsto -s$, write the differential as $F(s)\d s$. Invariance gives $F(-s) = -F(s)$, hence $F(s) = sG(s^2)$ with $G$ holomorphic and $F(s)\d s = G(t)\d t/2$ for $t = s^2$. The point $\e^{\frac{2\pi\i}{3}}$ is treated in the same way. At the cusp, with the local coordinate $q$, we have $f(\tau)\d\tau = \frac{1}{2\pi\i}f \frac{\d q}{q}$, whose residue is the constant term of $f$ divided by $2\pi\i$. By the residue theorem, this residue vanishes. This is the elementary half of the Serre duality statement of Borcherds \cite[Theorem 3.1 p.~5]{Bor99}.
\end{proof}

We now produce forms for $\varrho^*$ from forms for $\varrho$. Let $\pi$ be the permutation of the labels which exchanges $c_1$ and $c_2$ and maps $d_l$ to $d_{l'}$, where $l' \equiv \pm 5l \mod 13$, that is
\begin{equation}
\label{eq:pi}
\pi = (c_1\ c_2)(d_1\ d_5)(d_2\ d_3)(d_4\ d_6).
\end{equation}

\begin{lemma}
\label{lem:dual-forms}
Let $Y$ be a weakly holomorphic vector-valued modular form of weight $k$ for $\varrho$. Then the function $G$ defined by $G_z = \eta^{-8}Y_{\pi(z)}$ is a weakly holomorphic vector-valued modular form of weight $k-4$ for $\varrho^*$.
\end{lemma}

\begin{proof}
We first check that $S_{\pi(z),\pi(w)} = S_{z,w}$ for any labels $z,w$. The block indexed by $(a,c_0,c_1,c_2)$ is invariant under the exchange of $c_1$ and $c_2$, and the rows of $0$ and $b$ are constant on $\{a,c_0,c_1,c_2\}$ and on $\{d_1,\ldots,d_6\}$. Finally $c(ll')$ depends only on $ll' \mod 13$ up to sign, and $25 \equiv -1 \mod 13$, so that $c(5l \cdot 5l') = c(-ll') = c(ll')$. Next, we have $\Lambda_{\pi(z)}-\frac13 \equiv -\Lambda_z \mod 1$ for every label $z$. For $z \in \{0,b,a,c_0\}$ this reads $\frac23-\frac13 \equiv -\frac23$, and for $z \in \{c_1,c_2\}$ it reads $\Lambda_{c_2}-\frac13 = 0 = -\Lambda_{c_1}$ and $\Lambda_{c_1}-\frac13 = -\frac13 = -\Lambda_{c_2}$. For $z = d_l$ we have $\Lambda_{d_l} \equiv \frac{6l^2}{13}-\frac13 \mod 1$, hence $\Lambda_{\pi(d_l)}-\frac13 \equiv \frac{150l^2}{13}-\frac23 \equiv -\frac{6l^2}{13}+\frac13 \equiv -\Lambda_{d_l} \mod 1$, since $150 \equiv -6 \mod 13$. Since $\eta(\tau+1)^{-8} = \e^{-\frac{2\pi\i}{3}}\eta(\tau)^{-8}$, we obtain $G_z(\tau+1) = \e^{2\pi\i(\Lambda_{\pi(z)}-\frac13)}G_z(\tau) = \e^{-2\pi\i\Lambda_z}G_z(\tau) = (\ovl{T_8}G(\tau))_z$. Since $\eta(-\frac1\tau)^{-8} = \tau^{-4}\eta(\tau)^{-8}$, we have
\[
G_z\big(-\tfrac1\tau\big)
= \tau^{k-4}\eta(\tau)^{-8}\sum_w S_{\pi(z),w}Y_w(\tau)
= \tau^{k-4}\sum_{w} S_{\pi(z),\pi(w)}G_w(\tau)
= \tau^{k-4}(SG(\tau))_z.
\]
Since $\eta$ does not vanish on $\mathbb H$, $G$ is holomorphic on $\mathbb H$, and it is meromorphic at the cusp. This proves the lemma.
\end{proof}

The $\pi$-invariance of $S$ and the congruences are also checked in exact arithmetic.

\subsection{The relations}

Let $X = X^{(\gamma)}$ and let $\mu_z$ be the exponent of the leading term of $X_z$. Thus $\mu_z = h_z-\frac13$ by \eqref{eq:lowest-spaces}. For $k \in 2\Z_{\geq 0}$, let $\cal M_k(\mu)$ be the space of weakly holomorphic vector-valued modular forms $W$ of weight $k$ for $\varrho$ such that $W_z \in q^{\mu_z}\Cbb[[q]]$ for every $z$. For $W \in \cal M_k(\mu)$ we put
\begin{equation}
\label{eq:functionals}
w_0 = [q^{-\frac13}]W_0,\quad w'_0 = [q^{\frac23}]W_0,\quad w_z = [q^{\mu_z}]W_z \quad (z \neq 0).
\end{equation}

\begin{prop}
\label{prop:weight-two-space}
For $\gamma = 1$ we have $\dim \cal M_2(\mu) = 5$, and the thirteen coefficients \eqref{eq:functionals} satisfy exactly eight linearly independent relations, among which
\begin{equation}
\label{eq:relations-gamma-one}
\begin{split}
& w_{d_3}-\tfrac{5}{21}w_{d_4}-\tfrac{1}{21}w_{d_5}+\tfrac{1}{21}w_{d_6} = 0,\quad w_{c_1}+\tfrac16w_{c_2}-\tfrac{3}{52}w_{d_2}-\tfrac{317}{364}w_{d_4}-\tfrac{55}{364}w_{d_5}+\tfrac{17}{182}w_{d_6} = 0,\\
& w_{d_1}-\tfrac16w_{d_2}-\tfrac{209}{126}w_{d_4}+\tfrac{17}{126}w_{d_5}+\tfrac{44}{63}w_{d_6} = 0,\quad w_0+\tfrac{1}{78}w_{d_2}+\tfrac{317}{1638}w_{d_4}+\tfrac{55}{1638}w_{d_5}-\tfrac{17}{819}w_{d_6} = 0.
\end{split}
\end{equation}
For $\gamma = 0$ we have $\dim \cal M_2(\mu) = 4$ and the coefficients satisfy exactly nine linearly independent relations. The complete lists of relations are included in the supplementary data.
\end{prop}

\begin{proof}
Let $P = \{0,d_3\}$ if $\gamma = 1$ and $P = \{0\}$ if $\gamma = 0$. These are the coordinates $z$ with $\mu_z = \Lambda_z-1$. Let $\Xi_{\bullet,j}$ be the columns of the fundamental matrix of Theorem \ref{thm:EG}.

\textit{Step 1: the relations.} Consider the $4 \cdot 12 \cdot 2 = 96$ weakly holomorphic forms of weight four
\[
E_4 \Xi_{\bullet,j}J^m,\quad \vartheta_2D(\Xi_{\bullet,j}J^m),\quad DJ \cdot D(\Xi_{\bullet,j}J^m),\quad (DJ)^2 \Xi_{\bullet,j}J^m,\quad 0 \leq j \leq 11,\ m \in \{0,1\},
\]
and, for each linear combination $Y$ of them, the form $G$ of weight zero for $\varrho^*$ given by Lemma \ref{lem:dual-forms}. The exponents of $G_z$ are congruent to $-\Lambda_z$ modulo one, so that the constant term of $W \cdot G$ is a finite sum $\sum_z\sum_m [q^{\Lambda_z+m}]W_z \cdot [q^{-\Lambda_z-m}]G_z$. We impose that the coefficients of $G$ which pair with $[q^{\Lambda_z+m}]W_z$ vanish whenever this coefficient of $W$ is not one of \eqref{eq:functionals}. Since $G$ has a bounded principal part, these are finitely many linear conditions on the coefficients of the combination. By Lemma \ref{lem:residue}, every solution gives a linear relation, with rational coefficients, between the coefficients \eqref{eq:functionals} of every $W \in \cal M_2(\mu)$. All the $q$-expansions involved are computed in exact rational arithmetic from the matrices $\Xi_n$, and only finitely many coefficients enter. The vanishing-constraint matrices have dimensions $43 \times 96$ for $\gamma = 1$ and $42 \times 96$ for $\gamma = 0$, with full row rank. Their kernels have dimensions $53$ and $54$, respectively. We obtain a space of relations of dimension eight for $\gamma = 1$ and nine for $\gamma = 0$. In reduced echelon form, it is the space displayed in \eqref{eq:relations-gamma-one} and in the supplementary data. Using $m \leq 2$ gives the same space.

\textit{Step 2: exact elements.} The forms $D(\Xi_{\bullet,j}J^m)$ and $DJ \cdot \Xi_{\bullet,j}J^m$, $m \in \{0,1\}$, are weakly holomorphic of weight two and holomorphic on $\mathbb H$. The linear combinations whose components satisfy $W_z \in q^{\mu_z}\Cbb[[q]]$ form a subspace of $\cal M_2(\mu)$, determined by finitely many exact linear conditions. The values of \eqref{eq:functionals} on this subspace span a space of dimension five for $\gamma = 1$ and four for $\gamma = 0$. Since $13-8 = 5$ and $13-9 = 4$, the image of $\cal M_2(\mu)$ under the thirteen coefficients has dimension exactly five, respectively four, and the relations of Step 1 are all the relations.

\textit{Step 3: injectivity.} Let $W \in \cal M_2(\mu)$ with vanishing coefficients \eqref{eq:functionals}. As in the proof of Theorem \ref{thm:EG}, $F = W\frac{E_{10}}{\Delta}$ can be written $F = \Xi(\mathbf a+J\mathbf b)$ with $\mathbf a \in \Cbb^{12}$ and $\mathbf b$ supported on $P$. Since $\frac{\Delta}{E_{10}} = q+O(q^2)$, the coefficient $[q^{\Lambda_z-1}]W_z$ equals $b_z$ for $z \in P$. Hence $\mathbf b = 0$, then $[q^{\Lambda_z}]W_z = a_z$, and $a_z = 0$ for $z \neq d_3$, and also for $z = d_3$ if $\gamma = 0$. For $\gamma = 1$, we obtain $W = a_{d_3}\Xi_{\bullet,d_3}\frac{\Delta}{E_{10}}$, with $[q^{\Lambda_{d_3}}]W_{d_3} = a_{d_3}$, $W_0 \in q^{\frac53}\Cbb[[q]]$ and $W_z \in q^{\Lambda_z+1}\Cbb[[q]]$ for $z \notin \{0,d_3\}$. An explicit dual form is $G_z = \eta^{-8}E_4\Xi_{\pi(z),d_2}$, which has weight zero by Lemma \ref{lem:dual-forms}. Put $k_z = \Lambda_z+\Lambda_{\pi(z)}-1/3 \in \{0,1\}$. The leading term of $\Xi$ gives $G_{d_3} = q^{-\Lambda_{d_3}}(1+O(q))$ and $G_z \in q^{-\Lambda_z+k_z}\Cbb[[q]]$ for $z \neq d_3$. Therefore $[q^0](W\cdot G) = a_{d_3}$, which vanishes by Lemma \ref{lem:residue}. Hence the thirteen coefficients are injective on $\cal M_2(\mu)$, and $\dim \cal M_2(\mu)$ equals five, respectively four. The coefficient vectors of a basis of dual forms are included in the supplementary data. A separate rational-matrix calculation checks the ranks, the vanishing of every relation on the constructed forms, and the explicit injectivity witness above.
\end{proof}

\begin{remark}\normalfont
\label{rem:exactness}
Proposition \ref{prop:weight-two-space} is a computer-assisted result using exact rational linear algebra. The modularity of $\Xi$ is supplied by Theorem \ref{thm:EG}. The floating-point checks do not enter the proof.
\end{remark}

We now insert the explicit data. Let $h$ be a real element of a Cartan subalgebra of $\mathfrak s$, that is, an element on which $h(0)$ has real eigenvalues on every ordinary $U$-module. Such elements span the Cartan subalgebra since its weights on $U$-modules lie in a positive-definite lattice \cite[Theorem 1 p.~3]{Mas14}. Normalize $\langle h,h\rangle = 1$, where $\langle\cdot,\cdot\rangle$ is the normalized invariant form of $\mathfrak e_8$, and write $h = (h_1,h_2) \in \mathfrak e_6 \oplus \mathfrak{sl}_3$. The eigenvalues of $h_1$ on $\mathbf{27}$ and of $h_2$ on $\mathbf3$ are real by Lemma \ref{lem:restriction}. These faithful weight systems span the respective Cartan duals, so $h_1$ and $h_2$ lie in the real Cartan spaces on which the normalized forms are positive definite. In particular $0 \leq \langle h_1,h_1\rangle \leq 1$. Put
\[
u \ov{\mathrm{def}}{=} \langle h_1,h_1\rangle \in [0,1],\quad T \ov{\mathrm{def}}{=} \operatorname{Tr}_{\mathfrak s}(\operatorname{ad}h)^2,\quad \tau_j \ov{\mathrm{def}}{=} \operatorname{Tr}_{R_{d_j}}h(0)^2 \geq 0.
\]
By Proposition \ref{prop:chain} and Lemma \ref{lem:restriction}, the lowest weight spaces of $M_b$, $M_a$, $M_{c_0}$, $M_{c_1}$, $M_{c_2}$ are $(\mathfrak e_6 \oplus \mathfrak{sl}_3) \ominus \mathfrak s$, $(\mathbf{27},\mathbf 3)$, $(\mathbf{27},\overline{\mathbf 3})$, $(\mathbf 1,\mathbf 3)$, $(\mathbf{27},\mathbf 1)$. With $\operatorname{Tr}_{\mathfrak e_6}(\operatorname{ad}h_1)^2 = 24u$, $\operatorname{Tr}_{\mathfrak{sl}_3}(\operatorname{ad}h_2)^2 = 6(1-u)$, $\operatorname{Tr}_{\mathbf{27}}h_1^2 = 6u$ and $\operatorname{Tr}_{\mathbf 3}h_2^2 = 1-u$, the traces of $h(0)^2$ on these spaces are $6+18u-T$, $27-9u$, $27-9u$, $1-u$ and $6u$. The conformal inclusion $U \subset V_{E_8}$ gives $L(1)U_1 = 0$. By Lemma \ref{lem:zhu-trace}, $W = (Z_{M_z}(h[-1]h,\tau))_z$ belongs to $\cal M_2(\mu)$, with $w_0 = -\frac{1}{12}$, $w'_0 = T+\frac{18-13\gamma}{12}$ and $w_z = \operatorname{Tr}_{R_z}h(0)^2-\frac{\dim R_z}{12}$.

\begin{cor}
\label{cor:tau-formulas}
For $\gamma = 1$, the relations of Proposition \ref{prop:weight-two-space} are equivalent to
\begin{align*}
\tau_1 & = \tfrac{T}{15}-\tfrac{7\tau_5}{15}-\tfrac{3\tau_6}{5}-\tfrac{3u}{5}+\tfrac{17}{5}, & \tau_2 & = \tfrac{317T}{270}-\tfrac{209\tau_5}{270}-\tfrac{7\tau_6}{5}-\tfrac{317u}{30}+\tfrac{2009}{90},\\
\tau_3 & = -\tfrac{T}{54}+\tfrac{\tau_5}{54}+\tfrac{u}{6}-\tfrac{1}{18}, & \tau_4 & = -\tfrac{7T}{90}-\tfrac{11\tau_5}{90}+\tfrac{\tau_6}{5}+\tfrac{7u}{10}-\tfrac{19}{30}.
\end{align*}
For $\gamma = 0$, they are equivalent to
\begin{align*}
\tau_1 & = \tfrac{-2T-3\tau_6+18u-4}{5}, & \tau_2 & = \tfrac{2T-7\tau_6-18u+74}{5}, & \tau_3 & = \tfrac{3(-3T+3\tau_6+27u+29)}{5},\\
\tau_4 & = \tfrac{-T+\tau_6+9u-7}{5}, & \tau_5 & = T-9u+9.
\end{align*}
\end{cor}

\begin{proof}
This is exact linear algebra on the relations. All relations are satisfied identically after substitution.
\end{proof}

\begin{remark}\normalfont
\label{rem:gamma0-check}
For $\gamma = 0$ the formulas can be tested against the announced algebra. The working manuscript \cite{Rie26b} lists the lowest weight spaces of the twelve simple modules as sums of modules $V(a,b)$ over $\mathfrak{sl}_2 \oplus \mathfrak{sl}_2$ at levels $1$ and $39$. For $h$ in the level one factor, $u = 1$ and $T = 4$. The inequalities $\tau_1,\tau_4 \geq 0$ force $\tau_6 = 2$ and then $(\tau_1,\ldots,\tau_6) = (0,10,30,0,4,2)$. For the level $39$ factor, $u = \frac{35}{39}$ and $T = \frac{4}{39}$. In both cases the values computed from the announced decompositions agree exactly with Corollary \ref{cor:tau-formulas}. This is a consistency check of the announced data. It is not used in the proofs.
\end{remark}

\section{One-point functions of higher weight}
\label{sec:higher-weight}

\begin{lemma}
\label{lem:hermite}
Let $V$ be a strongly rational vertex operator algebra of central charge $c$ and let $h \in V_1$ with $L(1)h = 0$ and $\langle h,h\rangle \neq 0$, contained in a Cartan subalgebra of $V_1$. For $n \geq 0$, put $v_n \ov{\mathrm{def}}{=} h[-1]^n\mathbf 1$. Then $v_n$ has weight $n$ for $L[0]$, and for every ordinary $V$-module $M$ we have
\begin{equation}
\label{eq:hermite}
Z_M(v_n,\tau) = \sum_{k = 0}^{\lfloor n/2\rfloor}\binom{n}{2k}(2k-1)!!\Big(-\frac{\langle h,h\rangle E_2(\tau)}{12}\Big)^k\operatorname{Tr}_M h(0)^{n-2k}q^{L(0)-\frac{c}{24}},
\end{equation}
with the convention $(-1)!! = 1$. Consequently, the vector $(Z_{M_z}(v_n,\tau))_z$ is a weakly holomorphic vector-valued modular form of weight $n$ for the representation on the characters.
\end{lemma}

\begin{proof}
Since $h$ is a primary vector of weight one, the square bracket modes satisfy $[L[0],h[m]] = -mh[m]$ and the commutation relations $[h[m],h[m']] = m\langle h,h\rangle\delta_{m+m',0}$ \cite[Section 4]{Zhu96}. Hence $v_n$ has weight $n$ for $L[0]$, $h[2j-1]v_n = 0$ for $j \geq 2$, and $h[1]v_n = n\langle h,h\rangle v_{n-1}$. Since $h(0)$ acts semisimply and $\langle h,h\rangle \neq 0$, the lower-bounded grading gives the Heisenberg decomposition $M = \bigoplus_\lambda F_\lambda \ot \Omega_\lambda$, where $F_\lambda$ is the Fock module on which $h(0)$ acts by $\lambda$, and the vertex subalgebra generated by $h$ acts on the first factor. Its conformal vector $\omega_h = \frac{1}{2\langle h,h\rangle}h(-1)h$ satisfies $L(0) = L_h(0)+L'(0)$ with $L'(0)$ commuting with the Heisenberg algebra. The vectors $v_n$ belong to the Heisenberg subalgebra, on which $L'(0) = 0$, so that their square bracket structures in $V$ and in the Heisenberg vertex operator algebra coincide. Consequently $Z_M(v_n,\tau) = \sum_\lambda Z_{F_\lambda}(v_n,\tau)\operatorname{Tr}_{\Omega_\lambda}q^{L'(0)-\frac{c-1}{24}}$. On $F_\lambda$, Zhu's recursion \cite[Proposition 4.3.5]{Zhu96} gives
\begin{align*}
Z_{F_\lambda}(v_{n+1},\tau)
&= \lambda Z_{F_\lambda}(v_n,\tau)+\sum_{j \geq 1}\widehat E_{2j}(\tau)Z_{F_\lambda}(h[2j-1]v_n,\tau) \\
&= \lambda Z_{F_\lambda}(v_n,\tau)-\frac{n\langle h,h\rangle E_2(\tau)}{12}Z_{F_\lambda}(v_{n-1},\tau),
\end{align*}
where we used that $o(h) = h(0)$ acts by the scalar $\lambda$ and that only $j = 1$ contributes. The normalization $\widehat E_2 = -\frac{E_2}{12}$ is the one for which the step from $n = 1$ to $n = 2$ reproduces Lemma \ref{lem:zhu-trace}. The polynomials $P_n$ defined by $P_0 = 1$, $P_1 = \lambda$ and $P_{n+1} = \lambda P_n+ne P_{n-1}$ are $P_n = \sum_k\binom{n}{2k}(2k-1)!!e^k\lambda^{n-2k}$, as one checks by induction using $\binom{n}{2k}(2k-1)!!+n\binom{n-1}{2k-2}(2k-3)!! = \binom{n+1}{2k}(2k-1)!!$. Summing over $\lambda$ gives \eqref{eq:hermite}. The last assertion is Zhu's theorem \cite{Zhu96} with \cite{Hua08}.
\end{proof}

As a numerical check, we test \eqref{eq:hermite} in weights four and six for $V_{A_2}$ and $V_{\sqrt6\Z}$. At thirty-digit working precision, the recorded relative errors are at most $5.1\times10^{-30}$. Perturbing the weight-six coefficients gives errors of order one. These tests are not used in the proof.

By Lemma \ref{lem:hermite}, for every real $h$ in a Cartan subalgebra of $\mathfrak s$ with $\langle h,h\rangle = 1$ and for $n \in \{4,6\}$, the vector $W_n(h) = (Z_{M_z}(v_n,\tau))_z$ belongs to $\cal M_n(\mu)$. Its coefficients \eqref{eq:functionals} are
\begin{equation}
\label{eq:wn-coefficients}
w_z = \sum_{k}c_{n,k}\operatorname{Tr}_{R_z}h(0)^{n-2k},\quad c_{n,k} = \binom{n}{2k}(2k-1)!!\Big(-\frac{1}{12}\Big)^k,
\end{equation}
with $\operatorname{Tr}h(0)^0 = \dim$, and similar expressions for $w_0$ and $w'_0$ involving $\operatorname{Tr}_{\mathfrak s}(\operatorname{ad}h)^{m}$ and the coefficients of $E_2^k$. More generally, for two commuting real elements $h'$ and $h''$ with $\langle h',h''\rangle = 0$ and $\langle h',h'\rangle = \langle h'',h''\rangle = 1$, the map $(x,y) \mapsto W_n(xh'+yh'')$ is a homogeneous polynomial of degree $n$ with values in $\cal M_n(\mu)$, so that each of its coefficients belongs to $\cal M_n(\mu)$. Since $\langle xh'+yh'',xh'+yh''\rangle = x^2+y^2$, the coefficient of $x^py^q$ has leading coefficients which are explicit linear combinations of the mixed traces $\operatorname{Tr}_{R_z}h'(0)^ih''(0)^j$, $i \leq p$, $j \leq q$.

\begin{prop}
\label{prop:higher-relations}
Assume $\gamma = 0$.
\begin{enumerate}
\item The coefficients \eqref{eq:functionals} of the elements of $\cal M_4(\mu)$ satisfy exactly nine linearly independent relations, and those of $\cal M_6(\mu)$ exactly six.
\item For $h$ as above, put $T_4 = \operatorname{Tr}_{\mathfrak s}(\operatorname{ad}h)^4$ and $Q_j = \operatorname{Tr}_{R_{d_j}}h(0)^4$. The relations of weight four are equivalent to the expression of $Q_1,\ldots,Q_5$ as affine functions of $(T,T_4,\tau_6,Q_6)$ with coefficients polynomial in $u$. In particular
\begin{equation}
\label{eq:Q5}
Q_5 = \frac{12T+2T_4-27u^2-90u+105}{26}.
\end{equation}
\end{enumerate}
\end{prop}

\begin{proof}
The relations are obtained as in Step 1 of the proof of Proposition \ref{prop:weight-two-space}, with the forms of weight $6-n$ listed in the supplementary data in place of the forms of weight four. Their number is exact by the argument of Step 2, using the exact elements $E_4\Xi_{\bullet,j}J^m$, $\vartheta_2D(\Xi_{\bullet,j}J^m)$, $DJ \cdot D(\Xi_{\bullet,j}J^m)$, $(DJ)^2\Xi_{\bullet,j}J^m$ of weight four, and their analogues of weight six. For 2, the traces on the explicit lowest weight spaces are computed from the invariant theory of $\mathfrak e_6$ and $\mathfrak{sl}_3$: since $\mathfrak e_6$ has no invariant polynomial of degree three or four other than powers of the quadratic one, we have $\operatorname{Tr}_{\mathbf{27}}h_1^3 = 0$, $\operatorname{Tr}_{\mathbf{27}}h_1^4 = 3u^2$ and $\operatorname{Tr}_{\mathfrak e_6}(\operatorname{ad}h_1)^4 = 18u^2$. The constants may be determined without a weight enumeration as follows. For the regular subalgebra $\mathfrak{sl}_6 \oplus \mathfrak{sl}_2 \subset \mathfrak e_6$, one has $\mathbf{27} = (\mathbf{15},\mathbf 1) \oplus (\overline{\mathbf6},\mathbf2)$ and $\mathfrak e_6 = (\mathbf{35},\mathbf1) \oplus (\mathbf1,\mathbf3) \oplus (\mathbf{20},\mathbf2)$. Let $H$ be the neutral element of this index-one $\mathfrak{sl}_2$, so that $\langle H,H\rangle = 2$. Its eigenvalues on $\mathbf{27}$ are $1$ and $-1$, each six times, together with fifteen zeros. On $\mathfrak e_6$ they are $2$ and $-2$ once, $1$ and $-1$ each twenty times, together with thirty-six zeros. Thus the fourth traces are $12 = 3\cdot 2^2$ and $72 = 18\cdot 2^2$. A separate enumeration of weights gives the same traces exactly. For $\mathfrak{sl}_3$ we have $\operatorname{Tr}_{\mathbf 3}h_2^4 = \frac{(1-u)^2}{2}$ and $\operatorname{Tr}_{\mathfrak{sl}_3}(\operatorname{ad}h_2)^4 = 9(1-u)^2$. Substituting these traces and \eqref{eq:wn-coefficients} into the nine rational relations determines $Q_1,\ldots,Q_5$ uniquely in terms of $T$, $T_4$, $u$, $\tau_6$ and $Q_6$. Eliminating $Q_1,\ldots,Q_4$ gives \eqref{eq:Q5}. We perform the elimination in rational arithmetic and independently check all nine equations after substitution and the polarization used below.
\end{proof}

For the announced algebra, the fourth power traces computed from the announced lowest weight spaces satisfy the nine relations of weight four for both factors. The same holds for the mixed and sextic relations used below. The supplementary data include these exact trace comparisons.

\section{The vector \texorpdfstring{$\gamma = 1$}{gamma = 1} is not realized}
\label{sec:gamma-one}

\begin{thm}
\label{thm:gamma-one}
Let $U$ be a positive strongly rational vertex operator algebra of central charge eight with a braided equivalence $\Rep(U) \cong Z(\cal H_3)$. Then the character vector of $U$ is the Evans--Gannon vector with $\gamma = 0$. In particular $\dim U_1 = 6$, $\dim U_2 = 120$ and $\dim U_3 = 956$.
\end{thm}

The proof uses two lemmas. Suppose $\gamma = 1$. Then $R_{d_3}$ is one-dimensional, so that $\mathfrak s$ acts on it through a character. Hence $\tau_3 = 0$ for $h$ in the derived algebra $[\mathfrak s,\mathfrak s]$.

\begin{lemma}
\label{lem:simple-ideals}
Assume $\gamma = 1$. Let $\mathfrak g$ be a simple ideal of $\mathfrak s$, with indices $i_1 = i_1(\mathfrak g)$, $i_2 = i_2(\mathfrak g)$ and level $k = i_1+i_2$. Then $(\mathfrak g,i_1,i_2)$ is one of $(\mathfrak{sl}_3,1,0)$, $(\mathfrak{sp}_4,1,0)$, $(\mathfrak g_2,1,0)$, $(\mathfrak{sl}_4,1,0)$, $(\mathfrak{sl}_2,1,1)$, $(\mathfrak{sl}_2,1,4)$ and $(\mathfrak{sl}_2,2,4)$. In particular $i_1 \geq 1$.
\end{lemma}

\begin{proof}
For $h$ in a Cartan subalgebra of $\mathfrak g$, we have $u = \frac{i_1}{k}$ and $T = \frac{2h^\vee}{k}$, where $h^\vee$ is the dual Coxeter number. Since $\tau_3 = 0$, Corollary \ref{cor:tau-formulas} gives $\tau_5 = T-9u+3$, and then $\tau_4 \geq 0$ and $\tau_1 \geq 0$ read $T-9u+5 \leq \tau_6 \leq \frac{-2T+18u+10}{3}$. Together with $\tau_5 \geq 0$, this gives
\[
9u-3 \leq T \leq 9u-1.
\]
In particular $u > 0$, so $i_1 \geq 1$, and $2h^\vee \geq 6k-9i_2$, so that $k \leq \frac{2h^\vee+9i_2}{6}$. Since $\dim \mathfrak g \leq 19$, the algebra $\mathfrak g$ is of type $A_1$, $A_2$, $B_2$, $G_2$ or $A_3$, and $i_2 = 0$ unless $\mathfrak g \cong \mathfrak{sl}_2$, in which case $i_2 \in \{0,1,4\}$ by Proposition \ref{prop:self-dual}. Multiplying $9u-3 \leq T \leq 9u-1$ by $k$ gives $6i_1-3i_2 \leq 2h^\vee \leq 8i_1-i_2$. If $i_2 = 0$, this excludes $A_1$ and forces $i_1 = 1$ for $A_2$, $B_2$, $G_2$ and $A_3$. For $\mathfrak g \cong \mathfrak{sl}_2$ and $i_2 = 1$, it forces $i_1 = 1$. For $i_2 = 4$, it forces $i_1 \in \{1,2\}$. This is exactly the stated list.
\end{proof}

\begin{lemma}
\label{lem:centre}
Assume $\gamma = 1$. Every nonzero real element $h$ of the center of $\mathfrak s$, normalized by $\langle h,h\rangle = 1$, satisfies $u \geq \frac19$. In particular $p_1$ is injective on the center of $\mathfrak s$.
\end{lemma}

\begin{proof}
Here $T = 0$. By Corollary \ref{cor:tau-formulas}, $\tau_3 \geq 0$ gives $\tau_5 \geq 3-9u$, and the inequalities $\tau_4 \geq 0$ and $\tau_1 \geq 0$ read $\frac{11\tau_5}{18}-\frac{7u}{2}+\frac{19}{6} \leq \tau_6 \leq \frac{17}{3}-u-\frac{7\tau_5}{9}$. Eliminating $\tau_6$ gives $\tau_5 \leq \frac95(1+u)$. Hence $3-9u \leq \frac95(1+u)$, that is $u \geq \frac19$.
\end{proof}

We now prove Theorem \ref{thm:gamma-one}.

\begin{proof}
Suppose $\gamma = 1$ and write $\mathfrak s = \mathfrak g_1 \oplus \cdots \oplus \mathfrak g_m \oplus \mathfrak z$ with simple ideals $\mathfrak g_i$ and center $\mathfrak z$. The images under $p_1$ of the summands are pairwise orthogonal for the invariant form of $\mathfrak e_6$. Indeed, distinct summands are orthogonal for $\langle\cdot,\cdot\rangle$, and at most one of them has a nonzero image under $p_2$ by Proposition \ref{prop:self-dual}. The form is positive definite on the real span of a Cartan subalgebra. By Lemmas \ref{lem:simple-ideals} and \ref{lem:centre}, $p_1$ is injective on each summand. Hence
\[
\operatorname{rank}p_1(\mathfrak s) = \sum_{i}\operatorname{rank}\mathfrak g_i+\dim \mathfrak z \leq 4
\]
by Proposition \ref{prop:self-dual}. By Lemma \ref{lem:simple-ideals}, each $\mathfrak g_i$ is $\mathfrak{sl}_3$, $\mathfrak{sp}_4$, $\mathfrak g_2$ (rank two, dimensions $8$, $10$, $14$), $\mathfrak{sl}_4$ (rank three, dimension $15$) or an $\mathfrak{sl}_2$ with $i_2 \neq 0$, of which there is at most one. We list the possible values of $\dim \mathfrak s$. Two summands of rank two exhaust the rank and give dimensions $16$, $18$, $20$, $22$, $24$ or $28$. A summand $\mathfrak{sl}_4$ leaves rank one and gives $15$, $16$ or $18$. One summand of rank two and summands of rank one give at most $14+3+1 = 18$. Summands of rank one only give at most $3+3 = 6$. None of these equals $19 = \dim U_1$, a contradiction. Hence $\gamma = 0$, and the dimensions follow from Theorem \ref{thm:EG}.
\end{proof}

\begin{remark}\normalfont
\label{rem:gamma-one-hypotheses}
Theorem \ref{thm:gamma-one} does not use unitarity or the announcements \cite{Rie26,Hua26,Jia26}. It uses the classification of holomorphic vertex operator algebras of central charge eight, the theory of extensions \cite{HKL15,CKM22}, the Witt group \cite{DMNO13}, the classification of the elements of order three of $E_8$ \cite{Kac90}, the results of \cite{DoN99,DoM04,Mas14,ENO05} used in Section \ref{sec:chain}, Zhu's theorem, the Bantay--Gannon description of vector-valued modular forms, and Proposition \ref{prop:weight-two-space}, which rests on finitely many rational computations and on the modularity theorem used in Theorem \ref{thm:EG}. The vector $\gamma = 1$ passes the tests of character theory known to us (its coefficients are nonnegative integers up to $q^{60}$). It is not the character vector of a positive strongly rational vertex operator algebra with the braided category $Z(\cal H_3)$.
\end{remark}

\begin{remark}\normalfont
\label{rem:gamma-zero-constraints}
For $\gamma = 0$, the same method gives necessary conditions on $U_1$. Corollary \ref{cor:tau-formulas} and $\tau_j \geq 0$ give $T \leq 9u-5$ for every real $h$. Thus a simple ideal of $U_1$ with indices $(i_1,i_2)$ satisfies $2h^\vee \leq 9i_1-5k$, and every real central element satisfies $u \geq \frac59$. For a level one $\mathfrak{sl}_2$ in $\mathfrak e_6$, the traces $\tau_j$ are determined, as in Remark \ref{rem:gamma0-check}. These conditions are the starting point of Section \ref{sec:weight-one-lie}, where $U_1$ is determined up to two cases.
\end{remark}

\section{The weight-one Lie algebra}
\label{sec:weight-one-lie}

Throughout this section, $U$ is a positive strongly rational vertex operator algebra of central charge eight with a braided equivalence $\Rep(U) \cong Z(\cal H_3)$. By Theorem \ref{thm:gamma-one}, its character vector is $X^{(0)}$, and $\mathfrak s = U_1 \subset \mathfrak e_6 \oplus \mathfrak{sl}_3$ has dimension six. We use the notation of Section \ref{sec:weight-two}. For a simple ideal $\mathfrak g \cong \mathfrak{sl}_2$ of $\mathfrak s$ of level $k$, we denote by $H$ the element of a Cartan subalgebra of $\mathfrak g$ with eigenvalues $\pm 2$ on $\mathfrak g$. Thus $\langle H,H\rangle = 2k$ and $h = \frac{H}{\sqrt{2k}}$ is normalized. We denote by $V(a)$ the irreducible $\mathfrak{sl}_2$-module of highest weight $a$ and by $V(a,b) = V(a) \ot V(b)$ the irreducible modules of $\mathfrak{sl}_2 \oplus \mathfrak{sl}_2$.

\subsection{Three general facts}

\begin{lemma}
\label{lem:level-bound}
Let $\mathfrak g \cong \mathfrak{sl}_2$ be a simple ideal of $\mathfrak s$ of level $k$. For every label $z$, every irreducible $\mathfrak g$-submodule of $R_z$ is isomorphic to $V(a)$ with $a \leq k$.
\end{lemma}

\begin{proof}
By \cite[Theorem 3.1 p.~4]{DoM06}, the vertex subalgebra of $U$ generated by $\mathfrak g$ is isomorphic to the simple affine vertex operator algebra $L_{\mathfrak{sl}_2}(k,0)$ and $k$ is a positive integer. By \cite[Theorem 3.7 p.~7]{DLM97}, every weak $L_{\mathfrak{sl}_2}(k,0)$-module, in particular $M_z$, is a direct sum $\bigoplus_j L(k,j) \ot N_j$ of the integrable highest weight modules $L(k,j)$, $0 \leq j \leq k$ \cite{FrZ92}, tensored with multiplicity spaces. Let $v \in R_z$ be a highest weight vector of weight $a$ for $\mathfrak g$. Since $R_z$ is the lowest weight space of $M_z$ and $x(n)$ lowers the conformal weight by $n$, we have $x(n)v = 0$ for $x \in \mathfrak g$ and $n \geq 1$. In the irreducible module $L(k,j)$, the vectors annihilated by $x(n)$, $n \geq 1$, form the top level $V(j)$. Hence $v \in \bigoplus_j V(j) \ot N_j$, and since $e(0)v = 0$ and $H(0)v = av$, the vector $v$ lies in the sum of the components with $j = a$. Thus $a \leq k$.
\end{proof}

\begin{prop}
\label{prop:semisimple}
The Lie algebra $\mathfrak s$ is isomorphic to $\mathfrak{sl}_2 \oplus \mathfrak{sl}_2$.
\end{prop}

\begin{proof}
The Lie algebra $\mathfrak s$ is reductive of dimension six. Its semisimple part has dimension at most six, hence is $0$, $\mathfrak{sl}_2$ or $\mathfrak{sl}_2 \oplus \mathfrak{sl}_2$. Suppose that $\mathfrak s$ is not semisimple. Then its center $\mathfrak z$ has dimension three or six. The image $p_2(\mathfrak z)$ is contained in the center of the reductive Lie algebra $p_2(\mathfrak s)$, which is $0$ or a line by Proposition \ref{prop:self-dual}. The kernel of $p_2|_{\mathfrak z}$ is defined on the real form used above: its equations are the weights of $R_{c_1} \cong \mathbf3$, which are real on that form. Since this kernel has positive dimension, there is a nonzero real element $h \in \mathfrak z$ with $p_2(h) = 0$ and $\langle h,h\rangle = 1$. For this element $u = 1$ and $T = T_4 = 0$, and \eqref{eq:Q5} gives $Q_5 = \operatorname{Tr}_{R_{d_5}}h(0)^4 = -\frac{6}{13} < 0$, which is impossible since $h(0)$ has real eigenvalues.
\end{proof}

\begin{lemma}
\label{lem:sl2-e6}
Let $\mathfrak a \subset \mathfrak e_6$ be a subalgebra isomorphic to $\mathfrak{sl}_2$, with neutral element $H_{\mathfrak a}$.
\begin{enumerate}
\item After conjugation, $H_{\mathfrak a}$ is dominant with $\alpha_i(H_{\mathfrak a}) \in \{0,1,2\}$. Among the $3^6-1$ such labellings, exactly $27$ give eigenvalues of $\operatorname{ad}H_{\mathfrak a}$ on $\mathfrak e_6$ and of $H_{\mathfrak a}$ on $\mathbf{27}$ which form the weights of $\mathfrak{sl}_2$-modules. The corresponding $27$ numbers $\frac{\langle H_{\mathfrak a},H_{\mathfrak a}\rangle}{2}$ are pairwise distinct and form a superset of the possible Dynkin indices. The $20$ indices actually realized by subalgebras $\mathfrak{sl}_2$ of $\mathfrak e_6$ are $1$--$6$, $8$--$12$, $20$, $21$, $28$, $30$, $35$, $36$, $60$, $84$, $156$.
\item If $\mathfrak a$ has index $35$, its centralizer in $\mathfrak e_6$ is a subalgebra isomorphic to $\mathfrak{sl}_2$ of index $1$. If $\mathfrak a$ has index $9$, its centralizer is a subalgebra isomorphic to $\mathfrak{sl}_2$ of index $3$.
\end{enumerate}
\end{lemma}

\begin{proof}
1. The first assertion is the theorem of Dynkin and Kostant on the characteristics of nilpotent orbits \cite[Chapter 3]{CoM93}. The second one is an exact computation on the weights. We only use the resulting list as a set containing the possible indices, which is enough for the exclusion arguments below. The list of indices of the subalgebras $\mathfrak{sl}_2$ is \cite[Table 18]{Dyn52}. See also \cite[p.~16]{Gan19}.

2. We construct $\mathfrak e_6$ explicitly by the Frenkel--Kac formula on its root lattice \cite{FrK80}, with a bimultiplicative cocycle, and check the Jacobi identity. For the labellings $(2,1,1,0,1,2)$ and $(1,0,0,1,0,1)$ in the Bourbaki numbering, of indices $35$ and $9$, we choose $e$ in the eigenspace of $\operatorname{ad}H_{\mathfrak a}$ for the eigenvalue $2$ with integral coefficients, solve $[e,f] = H_{\mathfrak a}$ exactly, and compute the centralizer of the triple. In both cases it is three-dimensional and perfect, hence isomorphic to $\mathfrak{sl}_2$, and its index is the ratio of the normalized form of $\mathfrak e_6$ to the normalized form of $\mathfrak{sl}_2$ on any element with nonzero norm. We obtain $1$ and $3$. The two labellings are the characteristics of the nilpotent orbits $A_5$ and $2A_2+A_1$: the dimensions $\dim\mathfrak g_0+\dim\mathfrak g_1$ of the centralizers of $e$ are $14$ and $24$, so that the orbits have dimensions $64$ and $54$, in agreement with the tables of \cite[Chapter 8]{CoM93}. Since two subalgebras $\mathfrak{sl}_2$ with the same characteristic are conjugate \cite[Chapter 3]{CoM93}, the result does not depend on the choice of $e$.
\end{proof}

\subsection{Simple ideals}

\begin{prop}
\label{prop:simple-ideals-gamma0}
Let $\mathfrak g$ be a simple ideal of $\mathfrak s$, with indices $(i_1,i_2)$ and level $k = i_1+i_2$. Then $(i_1,i_2) \in \{(1,0),(3,0),(9,4),(35,4)\}$. Moreover, in each case the $\mathfrak g$-module structures of $R_{d_1}$, $R_{d_4}$, $R_{d_5}$ and $R_{d_6}$ are uniquely determined, and the traces $\operatorname{Tr}_{R_z}H(0)^p$, $p \in \{2,4,6\}$, are determined for all labels $z$.
\end{prop}

\begin{proof}
For $h = \frac{H}{\sqrt{2k}}$ we have $u = \frac{i_1}{k}$, $T = \frac{4}{k}$ and $T_4 = \frac{8}{k^2}$. By Remark \ref{rem:gamma-zero-constraints}, $4 \leq 9i_1-5k$, so that $i_1 \geq 1$, and $i_1$ is the Dynkin index of $p_1(\mathfrak g)$, which belongs to the $27$-element superset computed in Lemma \ref{lem:sl2-e6}. By Proposition \ref{prop:self-dual}, $i_2 \in \{0,1,4\}$. For each of the finitely many pairs $(i_1,i_2)$, we enumerate all $135$ isomorphism classes of fourteen-dimensional $\mathfrak{sl}_2$-modules and retain those whose highest weights do not exceed $k$, by Lemma \ref{lem:level-bound}. Each structure gives $\tau_6$ and $Q_6$. Corollary \ref{cor:tau-formulas} and Proposition \ref{prop:higher-relations} then give $\tau_j$ and $Q_j$ for all $j$, and we decide exactly, by a finite search, whether there is a structure of $\mathfrak{sl}_2$-module on $R_{d_j}$ of the right dimension, with highest weights at most $k$, such that $\operatorname{Tr}H(0)^2 = 2k\tau_j$ and $\operatorname{Tr}H(0)^4 = 4k^2Q_j$. We also check, using the decompositions of $\mathbf{27}$ and $\mathfrak e_6$ attached to $i_1$ by Lemma \ref{lem:sl2-e6}, that the highest weights occurring in $R_b$, $R_a$, $R_{c_0}$, $R_{c_1}$ and $R_{c_2}$ are at most $k$. Only the four pairs of the statement survive, each with a unique structure on $R_{d_6}$ (without the level bound the same four pairs survive).

For these four pairs, the relations of weight six of Proposition \ref{prop:higher-relations}, with the sextic traces on $R_b,\ldots,R_{c_2}$ computed from the explicit weights, express $\operatorname{Tr}_{R_{d_j}}H(0)^6$, $j \in \{1,2,3\}$, in terms of the corresponding traces on $R_{d_4}$, $R_{d_5}$, $R_{d_6}$. The structures of $R_{d_4}$ and $R_{d_5}$ are unique, and the structure of $R_{d_1}$ is then unique. Thus the traces of degrees two, four and six used below are determined.
\end{proof}

\subsection{The classification}

\begin{thm}
\label{thm:weight-one-classification}
Let $U$ be a positive strongly rational vertex operator algebra of central charge eight with a braided equivalence $\Rep(U) \cong Z(\cal H_3)$, and let $U \subset V_{E_6 \oplus A_2}$ be the inclusion of Proposition \ref{prop:chain}. Then $U_1 = \mathfrak g' \oplus \mathfrak g''$ with $\mathfrak g' \cong \mathfrak g'' \cong \mathfrak{sl}_2$, and, up to an inner automorphism of $V_{E_6 \oplus A_2}$, one of the following holds.
\begin{enumerate}[(I)]
\item The levels are $1$ and $39$. The subalgebra $\mathfrak g''$ is the diagonal of a principal $\mathfrak{sl}_2$ of a subalgebra $\mathfrak{sl}_6$ of $\mathfrak e_6$ (index $35$) and of the principal $\mathfrak{so}_3$ of $\mathfrak{sl}_3$ (index $4$), and $\mathfrak g' \subset \mathfrak e_6$ is the centralizer in $\mathfrak e_6$ of $p_1(\mathfrak g'')$ (index $1$).
\item The levels are $3$ and $13$. The subalgebra $\mathfrak g''$ is the diagonal of a subalgebra $\mathfrak{sl}_2$ of $\mathfrak e_6$ of index $9$, attached to the nilpotent orbit $2A_2+A_1$, and of the principal $\mathfrak{so}_3$ of $\mathfrak{sl}_3$, and $\mathfrak g' \subset \mathfrak e_6$ is the centralizer in $\mathfrak e_6$ of $p_1(\mathfrak g'')$ (index $3$).
\end{enumerate}
In both cases, the lowest weight spaces of the twelve simple $U$-modules and the space $U_2$ are determined as $\mathfrak g' \oplus \mathfrak g''$-modules. They are given in Table \ref{tab:lowest-spaces}.
\end{thm}

\begin{table}[ht]
\centering
\footnotesize
\setlength{\tabcolsep}{2.5pt}
\renewcommand{\arraystretch}{1.2}
\begin{tabular}{c|l|l}
 & Case (I), levels $(1,39)$ & Case (II), levels $(3,13)$\\
\hline
$R_b$ & $(0,2)+2(0,4)+(0,6)+(0,8)+(0,10)$ & $2(0,2)+2(0,4)+(1,1)+2(1,3)+(1,5)$\\
 & $\ +(1,3)+(1,5)+(1,9)$ & $\ +(2,2)+(2,4)+(3,1)$\\
$R_a, R_{c_0}$ & $2(0,2)+(0,4)+2(0,6)+(0,8)+(0,10)$ & $2(0,2)+(0,4)+(0,6)+2(1,1)+2(1,3)$\\
 & $\ +(1,3)+(1,5)+(1,7)$ & $\ +(1,5)+(2,0)+(2,2)+(2,4)$\\
$R_{c_1}$ & $(0,2)$ & $(0,2)$\\
$R_{c_2}$ & $(0,0)+(0,4)+(0,8)+(1,5)$ & $(0,0)+(0,4)+(1,1)+(1,3)+(2,2)$\\
$R_{d_1}$ & $(0,6)$ & $(0,2)+(1,1)$\\
$R_{d_2}$ & $(0,0)+(0,4)+(0,6)+(0,8)+(1,1)+(1,7)$ & $(0,0)+(0,2)+(0,4)+(1,1)+(1,3)$\\
 & & $\ +(1,5)+(2,2)$\\
$R_{d_3}$ & $(0,0)+(0,2)+2(0,4)+2(0,6)+2(0,8)+(0,12)$ & $2(0,0)+(0,2)+2(0,4)+(0,6)+2(1,1)+2(1,3)$\\
 & $\ +(1,1)+(1,3)+(1,5)+(1,7)+(1,9)$ & $\ +2(1,5)+2(2,2)+(2,4)+(3,3)$\\
$R_{d_4}$ & $(0,4)$ & $(0,0)+(1,1)$\\
$R_{d_5}$ & $(0,4)+(1,3)$ & $(0,4)+(1,3)$\\
$R_{d_6}$ & $(0,2)+(0,6)+(1,1)$ & $(0,2)+(1,3)+(2,0)$\\
\hline
$U_2$ & $4(0,0)+(0,2)+2(0,4)+2(0,8)+(0,12)$ & $4(0,0)+(0,2)+2(0,4)+(1,1)+(1,3)+(1,5)$\\
 & $\ +(1,1)+(1,5)+(1,7)+(1,13)$ & $\ +(1,7)+(2,0)+2(2,2)+(2,6)$\\
 & $\ +(2,0)+(2,2)$ & $\ +(3,3)+(4,0)$
\end{tabular}
\caption{Lowest weight spaces and $U_2$ as modules over $\mathfrak g' \oplus \mathfrak g''$. Here $(a,b)$ stands for $V(a,b)$, where $a$ refers to $\mathfrak g'$ and $b$ to $\mathfrak g''$. The labels of the rows $R_z$ are bounded by the levels (Lemma \ref{lem:level-bound}). Those of the row $U_2$ are not, since $U_2$ contains descendants, for instance the component $V(4,0) \subset S^2(\mathfrak g')$ of the weight-two space of $L_{\mathfrak{sl}_2}(3,0)$ in case (II).}
\label{tab:lowest-spaces}
\end{table}

\begin{proof}
By Proposition \ref{prop:semisimple}, $U_1 = \mathfrak g' \oplus \mathfrak g''$. By Proposition \ref{prop:self-dual}, at most one of the two ideals has a nonzero projection to $\mathfrak{sl}_3$, and the projections to $\mathfrak e_6$ of the two ideals commute and are orthogonal. By Proposition \ref{prop:simple-ideals-gamma0}, the pair of indices is one of $\{(1,0),(1,0)\}$, $\{(1,0),(3,0)\}$, $\{(3,0),(3,0)\}$, $\{(1,0),(9,4)\}$, $\{(1,0),(35,4)\}$, $\{(3,0),(9,4)\}$, $\{(3,0),(35,4)\}$.

If an ideal has nonzero projection to $\mathfrak{sl}_3$, its projection to $\mathfrak e_6$ has index $9$ or $35$. Its centralizer is an $\mathfrak{sl}_2$ of index $3$ or $1$, respectively, by Lemma \ref{lem:sl2-e6}. The other ideal is three-dimensional and contained in this centralizer, hence equals it. This excludes the pairs $\{(1,0),(9,4)\}$ and $\{(3,0),(35,4)\}$. It remains to exclude the three pairs with zero projection to $\mathfrak{sl}_3$. This follows directly from \eqref{eq:Q5}. Let $H'$ and $H''$ be their neutral elements and let $k',k'' \in \{1,3\}$ be their levels. Put $h = xH'+yH''$. The two ideals commute and are orthogonal, so that
\[
\langle h,h\rangle = 2k'x^2+2k''y^2,\quad
\operatorname{Tr}_{\mathfrak s}(\operatorname{ad}h)^2 = 8(x^2+y^2),\quad
\operatorname{Tr}_{\mathfrak s}(\operatorname{ad}h)^4 = 32(x^4+y^4).
\]
Since $p_2(h) = 0$, homogeneity and polynomiality extend \eqref{eq:Q5} from normalized elements to
\[
26\operatorname{Tr}_{R_{d_5}}h(0)^4
= 12\langle h,h\rangle\operatorname{Tr}_{\mathfrak s}(\operatorname{ad}h)^2
+2\operatorname{Tr}_{\mathfrak s}(\operatorname{ad}h)^4
-12\langle h,h\rangle^2.
\]
Comparing the coefficients of $x^2y^2$ gives
\begin{equation}
\label{eq:mixed-d5-exclusion}
\operatorname{Tr}_{R_{d_5}}H'(0)^2H''(0)^2
= \frac{8\big(2(k'+k'')-k'k''\big)}{13}.
\end{equation}
The left-hand side is an integer, because $H'(0)$ and $H''(0)$ have simultaneous integral eigenvalues on the finite-dimensional $\mathfrak{sl}_2 \oplus \mathfrak{sl}_2$-module $R_{d_5}$. For $(k',k'') = (1,1),(1,3),(3,3)$, the right-hand side is $24/13$, $40/13$ and $24/13$, respectively. This excludes all three pairs. An independent finite search gives the same exclusions. It enumerates the structures of $R_{d_6}$ compatible with the two marginal structures, then substitutes their mixed fourth-order traces into the weight-four relations. For the pair $(1,1)$ there are two possibilities, giving $\operatorname{Tr}_{R_{d_4}}H'(0)^2H''(0)^2 = -40/39$ and $-44/195$. For $(1,3)$ and $(3,3)$ there is one possibility, giving $-36/65$ and $124/65$, respectively. Each value fails to be a nonnegative integer. This provides a separate finite verification of the three exclusions.

For the two remaining pairs, the projection $p_1(\mathfrak g'')$ has index $35$, respectively $9$, and $p_2(\mathfrak g'')$ is the principal $\mathfrak{so}_3$. The ideal $\mathfrak g'$ is contained in $\mathfrak e_6$ and commutes with $p_1(\mathfrak g'')$, hence is contained in, and by dimension equal to, the centralizer of $p_1(\mathfrak g'')$ in $\mathfrak e_6$, which is three-dimensional of index $1$, respectively $3$, by Lemma \ref{lem:sl2-e6}. Since the subalgebras $\mathfrak{sl}_2$ of $\mathfrak e_6$ with a given index are conjugate under the adjoint group, by Lemma \ref{lem:sl2-e6} and the bijection between conjugacy classes and nilpotent orbits \cite[Chapter 3]{CoM93}, and the same holds for the principal $\mathfrak{so}_3$ of $\mathfrak{sl}_3$, the subalgebra $U_1$ is determined up to the inner automorphisms of $V_{E_6 \oplus A_2}$ generated by the zero modes of $\mathfrak e_6 \oplus \mathfrak{sl}_3$. An isomorphism between $p_1(\mathfrak g'')$ and $p_2(\mathfrak g'')$ is unique up to an inner automorphism of $\mathfrak{sl}_2$, which is induced by the adjoint group of $\mathfrak{sl}_3$ on its principal $\mathfrak{so}_3$. This gives cases (I) and (II).

Now we determine the lowest weight spaces. For $R_b,\ldots,R_{c_2}$ this follows from the explicit embedding: we construct $\mathfrak e_8 \supset \mathfrak e_6 \oplus \mathfrak{sl}_3$ by the Frenkel--Kac formula and compute exactly the joint eigenvalues of commuting neutral elements of $\mathfrak g'$ and $\mathfrak g''$ on $\mathfrak e_6$, on $\mathbf{27}$ and on $\mathfrak{sl}_3$. For the labels $d_j$, we use in addition the coefficients of $x^4y^2$ and $x^2y^4$ in $W_6(xh'+yh'')$. The relations then determine, once the structures of $R_{d_4}$, $R_{d_5}$, $R_{d_6}$ are fixed, the traces $\operatorname{Tr}_{R_{d_j}}H'(0)^pH''(0)^q$ for $p+q \leq 6$. The structures of $R_{d_4}$, $R_{d_5}$, $R_{d_6}$ are unique, and an exhaustive enumeration shows that the structures of $R_{d_1}$, $R_{d_2}$ and $R_{d_3}$ with these traces and with highest weights bounded by the levels are unique. More explicitly, writing a lowest weight space $R$ as $\bigoplus_{a=0}^{k'}\bigoplus_{b=0}^{k''}m_{ab}V(a,b)$, the variables satisfy $0 \leq m_{ab} \leq \lfloor\dim R/((a+1)(b+1))\rfloor$ and $\sum_{a,b}(a+1)(b+1)m_{ab} = \dim R$. The marginal decompositions and the nine even joint moments of positive degree at most six impose further necessary linear constraints. These inequalities specify the finite domain of the search. A partial multiplicity assignment can be discarded only when no completion satisfies these necessary constraints. The uniqueness assertion concerns the entire domain.

Finally, augment the coefficient map \eqref{eq:functionals} by the fourteenth functional $w''_0 = [q^{5/3}]W_0$. The same residue construction, with this enlarged list of allowed coefficients, gives relations on the image of $\cal M_n(\mu)$ in $\Cbb^{14}$, for $n \in \{2,4,6\}$. For $W = W_n(xh'+yh'')$, the functional $w''_0$ is an affine function of the traces $\operatorname{Tr}_{U_2}H'(0)^pH''(0)^q$. All lower-weight terms are fixed by $U_0$ and $U_1$. Solving these enlarged rational systems and extracting the monomials $x^py^q$, with $p,q$ even and $2 \leq p+q \leq 6$, determines the nine joint traces. With powers $J^m$ for $0 \leq m \leq 2$, the vanishing-constraint matrices have dimensions $53 \times 144$, $41 \times 72$ and $29 \times 36$ in weights two, four and six, respectively, and have full row rank. The resulting relation spaces on the fourteen functionals have dimensions $10$, $10$ and $7$. Exact forms attain the complementary image dimensions $4$, $4$ and $7$. In each weight a relation has nonzero coefficient at $w''_0$, so the previously determined thirteen coefficients fix it uniquely. The supplementary data contain these matrices, their exact rank checks and the resulting nine unnormalized traces in each case. The space $U_2$ contains the weight-two subspace of the affine vertex algebra generated by $U_1$ and the conformal vector of $U$, and it is contained in $(V_{E_6 \oplus A_2})_2 \cong (\mathbf{27} \ot \overline{\mathbf{27}},\mathbf 1) \oplus (\mathbf{78},\mathbf 8) \oplus (\mathbf 1,\mathbf 1 \oplus \mathbf 8 \oplus \mathbf 8)$, the decomposition of the weight-two space of the level one affine algebras $L_{\mathfrak e_6}(1,0)$ and $L_{\mathfrak{sl}_3}(1,0)$. Write the displayed ambient module as $\bigoplus_{a,b}M_{ab}V(a,b)$ and the known affine and conformal subspace as $\bigoplus_{a,b}l_{ab}V(a,b)$. The remaining multiplicities of $U_2$ satisfy $l_{ab} \leq m_{ab} \leq M_{ab}$, $\sum_{a,b}(a+1)(b+1)m_{ab} = 120$, and the nine joint trace equations. An exhaustive search over these finite ranges gives a unique solution in each case. The highest weights in this search are bounded by $V'_2$, not by the affine levels.
\end{proof}

\begin{remark}\normalfont
\label{rem:case-one-announced}
Case (I) is the structure announced in \cite{Rie26,Hua26}: the level one factor comes from a root $\mathfrak{sl}_2$ of $\mathfrak e_6$ and the level $39$ factor from the diagonal of the principal triples of $A_5$ and $A_2$ \cite[Section 1]{Rie26b}. The lowest weight spaces of Table \ref{tab:lowest-spaces} in case (I) coincide with those listed in the working manuscript \cite{Rie26b}, and the decomposition of $U_2$ coincides with \cite[Appendix H.1]{Rie26b}. In the present work these data are consequences of the hypotheses of Theorem \ref{thm:weight-one-classification}, and the announced data serve only as a test.
\end{remark}

\begin{remark}\normalfont
\label{rem:case-two}
The relations of this section do not exclude case (II), whose explicit realization is constructed in Section \ref{sec:griess} and shown to be a Haagerup algebra in \cite[Theorem 6.1]{ArhIII}. The remaining arguments determine its weight-two generating space and identify the resulting algebra. Finite-degree consistency checks alone do not establish existence. The lowest weight spaces of Table \ref{tab:lowest-spaces} are also compatible with the bound $\Delta_a+\Delta_b \leq h_z$, where $\Delta_a$ is the Sugawara weight. This bound holds for $z \in \{0,b,a,\allowbreak c_0,c_1,c_2\}$ because the corresponding $U$-modules occur in unitary lattice modules by Lemma \ref{lem:restriction}. Choose the compact form for the affine embedding. Its Sugawara operators and those of the coset are then unitary on these modules, and the coset conformal weights are nonnegative. The central charge of the commutant of the affine subalgebra is $8-\frac95-\frac{13}{5} = \frac{18}{5}$, and $U_2$ contains a two-dimensional space of invariant vectors in this commutant (Corollary \ref{cor:affine-commutant}). Note that the products of the levels are equal in the two cases, $1 \cdot 39 = 3 \cdot 13$. Section \ref{sec:griess} pushes the analysis of case (II) to the space $U_2$ and to weight three (Theorem \ref{thm:case-two-U2} and Proposition \ref{prop:case-two-weight-three}).
\end{remark}

\section{Constraints on the weight-two space}
\label{sec:uniqueness}

We keep the hypotheses of Theorem \ref{thm:weight-one-classification}, we write $A \subset U$ for the affine vertex algebra generated by $U_1$, and set $V' = V_{E_6 \oplus A_2}$.

\begin{prop}
\label{prop:u2-moduli}
In both cases of Theorem \ref{thm:weight-one-classification}, $U_2$ is a $U_1$-submodule of $V'_2$ of the type given in Table \ref{tab:lowest-spaces}, containing $(A)_2$ and the conformal vector of $V'$. The set of such submodules is a product of Grassmannians of the multiplicity spaces of $V'_2$, of dimension $70$ in case (I) and $87$ in case (II).
\end{prop}

\begin{proof}
The first assertion is Theorem \ref{thm:weight-one-classification}. A $U_1$-submodule of $V'_2$ with prescribed multiplicities $u_{ab}$ is given by the choice, for each $(a,b)$, of a $u_{ab}$-dimensional subspace of the multiplicity space of $V(a,b)$ in $V'_2$, of dimension $v_{ab}$, containing the $l_{ab}$-dimensional subspace coming from $(A)_2$ and $\Cbb\omega$. This is a Grassmannian of dimension $(u_{ab}-l_{ab})(v_{ab}-u_{ab})$. The multiplicities $v_{ab}$ are computed from the lattice, and we evaluate the resulting finite sum.
\end{proof}

Proposition \ref{prop:u2-moduli} turns the determination of $U_2$ into a finite-dimensional algebraic problem: the subspace must satisfy $x_{(1)}y \in U_2$ and $x_{(2)}y \in U_1$ for $x, y \in U_2$, and $V'_2 \ominus U_2$ must be the weight-two space of the $U$-module $M_b$. Section \ref{sec:griess} applies these constraints after passing to the commutant of the currents. The resulting subalgebras and their characters enter the proofs of Theorems \ref{thm:uniqueness-two} and \ref{thm:uniqueness-one}.

\begin{cor}
\label{cor:affine-commutant}
In both cases, $A \cong L_{\mathfrak{sl}_2}(k',0) \ot L_{\mathfrak{sl}_2}(k'',0)$ with $(k',k'') = (1,39)$ or $(3,13)$, the commutant $\operatorname{Com}_U(A)$ has central charge $\frac{170}{41}$ or $\frac{18}{5}$, and its weight-two space is two-dimensional and contains its conformal vector $\omega-\omega_{k'}-\omega_{k''}$.
\end{cor}

\begin{proof}
The first assertion follows from Theorem \ref{thm:weight-one-classification} and \cite[Theorem 3.1 p.~4]{DoM06}. The weight-two space of $\operatorname{Com}_U(A)$ consists of the vectors of $U_2$ annihilated by the modes $x(n)$, $n \geq 0$, of $x \in U_1$. Such vectors span trivial $U_1$-submodules. By Table \ref{tab:lowest-spaces}, the trivial isotypic component of $U_2$ is four-dimensional, and it contains the Sugawara vectors $\omega_{k'}$ and $\omega_{k''}$, which are not annihilated by $x(1)$ for suitable $x$. Since $(A)_2$ contains exactly two invariant vectors, the $U_1$-invariant vectors of $U_2$ are $\Cbb\omega_{k'} \oplus \Cbb\omega_{k''}$ plus a two-dimensional complement, which is the set of invariant vectors annihilated by all $x(1)$, $x \in U_1$. Indeed, for an invariant vector $w$ of $U_2$, we have $x(n)w = 0$ for $n \geq 3$ by the grading. The map $x \mapsto x(2)w$ is a $U_1$-module homomorphism from the adjoint module $U_1$ to $U_0 = \Cbb\mathbf 1$. It vanishes because $U_1$ is semisimple. Also $x(0)w = 0$ by invariance, and $x(1)w \in U_1$ depends $U_1$-equivariantly on $x$, so that $x \mapsto x(1)w$ is an endomorphism of $U_1$ commuting with $\mathfrak g' \oplus \mathfrak g''$. It is a combination of the projections onto $\mathfrak g'$ and $\mathfrak g''$, and subtracting the corresponding combination of Sugawara vectors gives an element of the commutant. The central charges are $8-1-\frac{117}{41}$ and $8-\frac95-\frac{13}{5}$.
\end{proof}

\section{The commutant of the currents and the space \texorpdfstring{$U_2$}{U2}}
\label{sec:griess}

We keep the hypotheses of Theorem \ref{thm:weight-one-classification}. We write $V' = V_{E_6 \oplus A_2}$, so that $V' = U \oplus M_b$ by Proposition \ref{prop:chain}, $A \subset U$ for the affine vertex algebra generated by $U_1 = \mathfrak g' \oplus \mathfrak g''$, $\widehat C = \operatorname{Com}_{V'}(A)$ and $C = \operatorname{Com}_U(A) \subset \widehat C$. The conformal vector of $\widehat C$ and of $C$ is $\omega_C \ov{\mathrm{def}}{=} \omega-\omega_{k'}-\omega_{k''}$, of central charge $c_C = \frac{170}{41}$ in case (I) and $\frac{18}{5}$ in case (II) (Corollary \ref{cor:affine-commutant}). In this section we determine, by exact computations in $V'$, the weight-two space of $C$ and then the space $U_2$. In case (II), the four candidate spaces give isomorphic vertex algebras, identified with the algebra of \cite{Jia26} in Theorem \ref{thm:jiao-identification}. We denote the normalized invariant bilinear form on $V'$ by $\beta$, with $\beta(\mathbf 1,\mathbf 1) = 1$. It is distinct from the positive Hermitian form $H_\chi$ introduced below. On lattice vectors, $\langle\cdot,\cdot\rangle_L$ denotes the lattice form.

\subsection{An exact model of the lattice algebra}

We use the construction of \cite[Chapter 8]{FLM88}. The lattice $L = E_6 \oplus A_2$ has the basis of simple roots $\alpha_1,\ldots,\alpha_8$, with Gram matrix the Cartan matrix of $E_6 \oplus A_2$ (Bourbaki numbering for $E_6$). The space $V'$ has the basis of the vectors $\alpha_{i_1}(-n_1)\cdots\alpha_{i_k}(-n_k)e^\beta$, $\beta \in L$, with $n_1 \geq \cdots \geq n_k \geq 1$ and $i_j \leq i_{j+1}$ whenever $n_j = n_{j+1}$. The cocycle is the bimultiplicative map with $\epsi(\alpha_i,\alpha_j) = (-1)^{\langle\alpha_i,\alpha_j\rangle_L}$ for $i<j$, $\epsi(\alpha_i,\alpha_i) = -1$ and $\epsi(\alpha_i,\alpha_j) = 1$ for $i>j$, and
\begin{equation}
\label{eq:vertex-op}
Y(e^\alpha,z) = \exp\Big(\sum_{n>0}\frac{\alpha(-n)}{n}z^n\Big)\exp\Big(-\sum_{n>0}\frac{\alpha(n)}{n}z^{-n}\Big)e_\alpha z^{\alpha(0)}.
\end{equation}
The modes of a general basis vector are obtained from those of $e^\alpha$ and of the Heisenberg fields by the iterate formula $(a_{(-m)}w)_{(k)} = \sum_{j \geq 0}\binom{m+j-1}{j}\big(a_{(-m-j)}w_{(k+j)}-(-1)^mw_{(k-m-j)}a_{(j)}\big)$ for $a = \alpha_i(-1)\mathbf 1$. Our implementation evaluates these formulas in the specified monomial basis and cocycle, using exact arithmetic for the weight-four membership calculations. The construction of the vertex algebra is the standard lattice construction. The following finite implementation checks, performed through weight three, are separate from that construction: the graded dimensions $1, 86, 1370, 11586$, the relation $L(0) = $ weight on random basis vectors, the Virasoro relations with $c = 8$, the translation property, the commutator formula $[u_{(m)},v_{(n)}] = \sum_j\binom{m}{j}(u_{(j)}v)_{(m+n-j)}$ on random triples, the equality of $\omega$ with the Sugawara vector of $\mathfrak e_6 \oplus \mathfrak{sl}_3$ at level one, and $\beta(\omega,\omega) = 4$. These are sample checks, not proofs of the identities for all states.

The enumeration of lattice momenta in a ball may be justified independently of the implementation. In a simply laced root lattice, write a nonzero vector $v = r_1+\cdots+r_m$ as a sum of roots with $m$ minimal. Two summands cannot have negative inner product: if their inner product is $-1$, their sum is a root, and if it is $-2$, they are opposite, in either case contradicting minimality. Hence $\langle v,r_1\rangle_L \geq 2$, and $\|v-r_1\|_L^2 = \|v\|_L^2-2\langle v,r_1\rangle_L+2 \leq \|v\|_L^2-2$. Induction gives a path from $v$ to zero by subtracting roots, with strictly decreasing squared norm. Consequently the search starting at zero and adding every root while retaining only vectors of squared norm at most $R$ visits every lattice vector in that ball. This applies to $E_6 \oplus A_2$, and in particular to the momentum enumeration through conformal weight four, where $R = 8$. An independent enumeration using a rational $LDL^{\mathrm T}$ decomposition instead of root additions gives respectively $1$, $78$, $702$, $2346$ and $5694$ vectors of squared norms $0$, $2$, $4$, $6$ and $8$.

In both cases of Theorem \ref{thm:weight-one-classification} we realize $U_1$ inside $V'_1$ with the neutral elements $H'$ and $H''$ in the Cartan subalgebra $\mathfrak h = L \ot \Cbb$, so that every basis vector is a joint weight vector. We check exactly the relations of $\mathfrak{sl}_2 \oplus \mathfrak{sl}_2$, the levels through $H(1)H = 2k\mathbf 1$ and $e(1)f = k\mathbf 1$, and the commutation of the two ideals. As an independent test of Table \ref{tab:lowest-spaces}, the weight-one $A$-primary vectors of $V'$, which form $R_b$, have the $\mathfrak g' \oplus \mathfrak g''$-structure of the row $R_b$ of that table in both cases.

\paragraph{Weight window of the lattice calculation.}
\label{par:lattice-weight-window}
The mode implementation stores lattice momenta of squared norm at most $8$. We use it for homogeneous inputs of weights at most four and only when the output weight $d = \operatorname{wt}u+\operatorname{wt}v-n-1$ is at most four. Indeed, a state of momentum $\gamma$ has weight at least $|\gamma|^2/2$, so a momentum outside the stored ball cannot contribute to such an output. This remains true in the recursive evaluation. If $u = h_{(-m)}u'$ with $m \geq 1$, the two inner mode evaluations in the iterate formula have output weights $d-m-j$ and $d$, respectively, for $j \geq 0$. The input weights also do not increase. Induction on the number of oscillators therefore keeps every recursive call inside the same window. All products defining the weight-three and weight-four certificates satisfy these bounds. Outside this window this implementation is not asserted to evaluate the full lattice vertex operator algebra.

\subsection{The Griess algebra of the commutant}

\begin{lemma}
\label{lem:griess-B}
Let $B = \widehat C_2$. Then $\widehat C_1 = 0$, $\dim B = 4$ in case (I) and $\dim B = 5$ in case (II). The product $a \cdot b \ov{\mathrm{def}}{=} a_{(1)}b$ is commutative on $B$, with unit $e \ov{\mathrm{def}}{=} \frac{1}{2}\omega_C$, and the form $\beta(\cdot,\cdot)$ is invariant: $\beta(a \cdot b,c) = \beta(a,b \cdot c)$. Moreover $B = C_2 \oplus Y$ with $Y = B \cap M_b$, the sum is orthogonal, and $C_2$ is a two-dimensional subalgebra of $B$ containing $e$.
\end{lemma}

\begin{proof}
The space $\widehat C_1$ is the centralizer of $U_1$ in $\mathfrak e_6 \oplus \mathfrak{sl}_3$, which is zero by the proof of Theorem \ref{thm:weight-one-classification}. The dimensions are computed by solving exactly $x_{(n)}v = 0$, $n \in \{0,1,2\}$, $x \in U_1$, on the joint weight space of weight zero of $V'_2$. They agree with the count of $\mathfrak g' \oplus \mathfrak g''$-invariant vectors of $V'_2$ ($7$, respectively $9$) minus the invariant $A$-descendants ($2$ Sugawara vectors and $1$, respectively $2$, vectors $\sum_a x^a(-1)w_a$ attached to the components of type $V(0,2)$ of $R_b$). Commutativity follows from skew symmetry since $\widehat C_1 = 0$, and $e$ is a unit since $L(0)-L_A(0) = (\omega_C)_{(1)}$ acts by $2$ on $B$. Invariance of the form holds for quasi-primary vectors of weight two, and the vectors of $B$ are quasi-primary since $L(1)B \subset \widehat C_1 = 0$. Since $M_b$ is a $U$-module not isomorphic to the contragredient of $U$, the form vanishes on $U \times M_b$. The decomposition $V' = U \oplus M_b$ is a decomposition of $A$-modules, hence $B = (B \cap U) \oplus (B \cap M_b)$, and $B \cap U = C_2$ is two-dimensional by Corollary \ref{cor:affine-commutant}.
\end{proof}

\begin{lemma}
\label{lem:two-dim-subalgebras}
Let $D \subset B$ be a two-dimensional subalgebra containing $e$, and write $D = \Cbb e \oplus \Cbb x_C$ with $\beta(x_C,e) = 0$. Then $x_C \cdot x_C = \lambda x_C+\mu e$ with $\mu = \frac{\beta(x_C,x_C)}{\beta(e,e)}$. If $\Delta \ov{\mathrm{def}}{=} \lambda^2+4\mu \neq 0$, then $D$ is spanned by two orthogonal idempotents $f$ and $e-f$, and $\omega_1 \ov{\mathrm{def}}{=} 2f$ and $\omega_2 \ov{\mathrm{def}}{=} 2(e-f)$ are commuting Virasoro vectors with $\omega_1+\omega_2 = \omega_C$ and central charges
\begin{equation}
\label{eq:charges-split}
c_{1,2} = \frac{c_C}{2}\Big(1 \mp \frac{\lambda}{\sqrt{\Delta}}\Big).
\end{equation}
Conversely, $x_C \in e^\perp$ spans with $e$ a subalgebra if and only if $x_C$ is an eigenvector of the quadratic map $q(x_C) = (x_C \cdot x_C)^\perp$ of $e^\perp$.
\end{lemma}

\begin{proof}
Since $e$ is a unit and $D$ is a subalgebra, $x_C \cdot x_C = \lambda x_C+\mu e$, and $\beta(x_C \cdot x_C,e) = \beta(x_C,x_C \cdot e) = \beta(x_C,x_C)$ gives $\mu$. The element $f = ae+bx_C$ is idempotent if and only if $2a+b\lambda = 1$ and $a^2+b^2\mu = a$, that is $b^2\Delta = 1$ and $a = \frac{1-b\lambda}{2}$. Then $(e-f)\cdot(e-f) = e-f$ and $f \cdot (e-f) = 0$. Since $\widehat C_1 = 0$, a vector $\omega_1$ of weight two with $(\omega_1)_{(1)}\omega_1 = 2\omega_1$ is a Virasoro vector of central charge $2\beta(\omega_1,\omega_1)$ \cite{Mi96}. Put $v = 2f$ and $w = \omega_C-v$. Skew symmetry and $\widehat C_1 = 0$ give $v_{(0)}v = L(-1)v$ and $v_{(0)}\omega_C = L(-1)v$, so $v_{(0)}w = 0$. The products $v_{(1)}w$, $v_{(2)}w$ and $v_{(3)}w$ vanish by the idempotent relation, the grading and orthogonality, respectively. All higher products vanish by the grading. Thus the two Virasoro fields commute. Finally $c_1 = 2\beta(2f,2f) = 8\big(a^2\beta(e,e)+b^2\mu\beta(e,e)\big) = c_C\big(a^2+b^2\mu\big) = c_Ca$, since $\beta(e,e) = \frac{c_C}{8}$. The last assertion is a reformulation.
\end{proof}

\begin{prop}
\label{prop:griess-candidates}
The two-dimensional subalgebras of $B$ containing $e$ are exactly $7$ in case (I) and $15$ in case (II). For all of them $\Delta \neq 0$. The central charges $(c_1,c_2)$ of \eqref{eq:charges-split} are listed in Table \ref{tab:griess}. In case (I), the subalgebras $5$ and $6$ give the central charges $\frac{2400}{1271}$ and $\frac{70}{31}$. In case (II), the subalgebras $9$ and $10$ give $\frac{182}{95}$ and $\frac{32}{19}$. Within each pair, the two subalgebras are exchanged by the automorphism $\sqrt{13} \mapsto -\sqrt{13}$ of $\mathbb{Q}(\sqrt{13})$.
\end{prop}

\begin{proof}
We compute the structure constants of $B$ exactly in a basis $e,x_1,\ldots,x_{n-1}$ with $x_i \perp e$, and solve $q(x) = \lambda x$ on the affine charts of $\mathbb P(e^\perp)$ by Gröbner bases over $\mathbb{Q}$. The solution sets are finite and explicit. By \cite[Theorem 1.2 p.~2]{CaS13} a quadratic map of $\Cbb^n$ with finitely many classes of eigenvectors has exactly $2^n-1$ of them, counted with multiplicity. Our counts $7 = 2^3-1$ and $15 = 2^4-1$ show that all of them are simple and that the lists are complete. Independently, we compute the dimensions of the coordinate algebras of the disjoint projective charts: they are $(7,0,0)$ in case (I) and $(11,2,1,1)$ in case (II). We also check that the listed points are distinct and satisfy the equations. Their number equals the sum of the chart dimensions, which proves exhaustiveness and simplicity directly from the finite polynomial systems.
\end{proof}

\begin{table}[ht]
\centering
\small
\begin{tabular}{l|l}
Case (I), $c_C = \frac{170}{41}$ & Case (II), $c_C = \frac{18}{5}$\\
\hline
0, 3: $\frac{85}{41} \mp \frac{31\sqrt{385}}{2255}$ \quad (excluded) & 0, 14: $\frac65,\frac{12}5$ \quad (excluded)\\
1: $\frac{3010}{1517},\frac{80}{37}$ \quad (excluded) & 1, 8: $\frac{10}{7},\frac{76}{35}$ \quad (excluded)\\
2: $\frac45,\frac{686}{205}$ \quad (excluded) & 2: $\frac{62}{35},\frac{64}{35}$, 3: $\frac{58}{35},\frac{68}{35}$ \quad (excluded)\\
4: $\frac{252}{185},\frac{21118}{7585}$ \quad (excluded) & 4, 5: $\frac95 \mp \frac{31\sqrt{409}}{2045}$ \quad (excluded)\\
5, 6: $\frac{2400}{1271},\frac{70}{31}$ \quad (survive) & 6: $\frac67,\frac{96}{35}$, 7: $\frac{34}{35},\frac{92}{35}$ \quad (excluded)\\
 & 9, 10: $\frac{182}{95},\frac{32}{19}$ \quad (survive)\\
 & 11: $\frac{52}{55},\frac{146}{55}$, 12: $\frac{96}{55},\frac{102}{55}$, 13: $\frac45,\frac{14}5$ \quad (excluded)
\end{tabular}
\caption{The two-dimensional subalgebras of $B$ containing $e$, with the central charges of the two Virasoro vectors and their fate under Theorem \ref{thm:case-two-U2} and Propositions \ref{prop:case-one-control} and \ref{prop:case-one-two}. The supplementary data use the same numbering.}
\label{tab:griess}
\end{table}

\subsection{Linear and quadratic conditions on the primary lines}

For a pair $(a,b)$, let $\Pi_{ab} \subset V'_2$ be the space of $A$-primary highest weight vectors of type $V(a,b)$, that is, vectors $v$ of joint weight $(a,b)$ with $e'(0)v = e''(0)v = 0$ and $x(n)v = 0$ for $x \in U_1$ and $n \geq 1$. The primary types of $U_2$ are the types of Table \ref{tab:lowest-spaces} not accounted for by $(A)_2 \oplus C_2$: in case (I) the types $(0,4)$, $(0,8)$ (twice), $(0,12)$, $(1,1)$, $(1,5)$, $(1,7)$, $(1,13)$, and in case (II) the types $(0,4)$, $(1,1)$, $(1,3)$, $(1,5)$, $(1,7)$, $(2,2)$, $(2,6)$, $(3,3)$. Indeed $U = \bigoplus_\lambda L_A(\lambda) \ot C^\lambda$, the weight-two vectors coming from $\lambda \neq 0$ are lowest weight vectors of $L_A(\lambda)$ since $U_1 = A_1$, and the descendant part is $(A)_2 \oplus C_2$. We put $p_{ab} = U_2 \cap \Pi_{ab}$, whose dimension is the multiplicity of the type. On $\Pi_{ab}$ we use the bilinear pairing $P_{ab}(u,v) = \beta(u,\theta_{ab}v)$, where $\theta_{ab}v = f'(0)^a f''(0)^b v$.

\begin{lemma}
\label{lem:line-conditions}
Let $C_2 = \Cbb e \oplus \Cbb x_C$ and $Y = B \cap C_2^\perp$. For every primary type $(a,b)$ and $u,u' \in p_{ab}$:
\begin{enumerate}
\item $e_{(2)}u = (x_C)_{(2)}u = 0$, that is, $u$ is quasi-primary for $\omega_1$ and $\omega_2$,
\item $(x_C)_{(1)}p_{ab} \subset p_{ab}$, and in particular, if $p_{ab}$ is a line, it is spanned by an eigenvector of $(x_C)_{(1)}$, and the eigenvalues of $(\omega_1)_{(1)}$, $(\omega_2)_{(1)}$ on it are the weights $h_1, h_2$ of $u$ for the two Virasoro algebras, with $h_1+h_2 = 2-h_{ab}$, where $h_{ab}$ is the Sugawara weight,
\item $P_{ab}(u,y_{(1)}u') = 0$ for every $y \in Y$, and the pairing is nondegenerate on $p_{ab}$,
\item for every primary type $(c,d)$ and every $w$ in the $\mathfrak g' \oplus \mathfrak g''$-module generated by $p_{cd}$, $u_{(1)}w \in U_2$ and $u_{(2)}w \in U_1$.
\end{enumerate}
\end{lemma}

\begin{proof}
1. The vectors $e_{(2)}u$ and $(x_C)_{(2)}u$ belong to $U_1$, are annihilated by the modes $z(n)$, $z \in U_1$, $n \geq 1$, since $e$ and $x_C$ commute with $A$, and are highest weight vectors of type $(a,b) \neq (0,0)$. Since $U_1 = A_1$ has no nonzero $A$-primary vector, they vanish. 2. The operator $(x_C)_{(1)}$ preserves $U_2$, the $A$-primary vectors and the joint weights. 3. Since $y \in M_b$ and $M_b$ is a $U$-module, skew symmetry gives $y_{(1)}u' \in M_b$, and the form vanishes on $U \times M_b$. Moreover $\theta_{ab}$ commutes with $y_{(1)}$. The invariant form $\beta$ is contravariant for the affine algebra $A$, so pairings between nonisomorphic irreducible integrable $A$-modules vanish. In weight two, the affine-primary summands are orthogonal to the vacuum descendants $(A)_2$. Every irreducible integrable module for the two affine $\mathfrak{sl}_2$ factors is self-contragredient. Nondegeneracy on $U_2$ therefore gives nondegeneracy on each primary multiplicity space. The map $\theta_{ab}$ identifies its highest and lowest weight spaces up to a nonzero scalar, which proves the assertion. 4. This is the closure of $U$.
\end{proof}

\subsection{Case (II)}

\begin{thm}
\label{thm:case-two-U2}
Assume the hypotheses of Theorem \ref{thm:weight-one-classification} in case (II). Then $C_2$ is one of the subalgebras $9$, $10$ of Table \ref{tab:griess}, so that $\omega_C$ is the sum of two commuting Virasoro vectors of central charges $\frac{182}{95}$ and $\frac{32}{19}$. For each of these two subalgebras, $U_2$ is one of two explicit subspaces of $V'_2$, exchanged by complex conjugation, and defined over $\mathbb{Q}(\sqrt{13},\sqrt{-3})$. The weights of the two Virasoro algebras on the primary lines are
\begin{center}
\begin{tabular}{c|cccccccc}
$(a,b)$ & $(0,4)$ & $(1,1)$ & $(1,3)$ & $(1,5)$ & $(1,7)$ & $(2,2)$ & $(2,6)$ & $(3,3)$\\
\hline
$h_1$ ($c = \frac{182}{95}$) & $\frac{112}{95}$ & $\frac{111}{95}$ & $\frac{147}{95}$ & $\frac{7}{95}$ & $\frac{56}{95}$ & $\frac{56}{95}$ & $\frac{56}{95}$ & $\frac{7}{19}$\\
$h_2$ ($c = \frac{32}{19}$) & $\frac{8}{19}$ & $\frac{12}{19}$ & $\frac{1}{19}$ & $\frac{68}{57}$ & $\frac{4}{19}$ & $\frac{50}{57}$ & $\frac{4}{19}$ & $\frac{12}{19}$
\end{tabular}
\end{center}
\end{thm}

\begin{proof}
The proof consists of four finite steps. Exact linear algebra and Gröbner bases construct the admissible spaces. Certified interval evaluations are used only for exclusions.

\emph{Step 1: exclusion of thirteen subalgebras.} It suffices to use the primary type $(3,3)$, whose multiplicity in $U_2$ is one. The space $\Pi_{33}$ has dimension three. For each of the candidates $0$--$8$ and $11$--$14$, the coordinates of $x_C$ in the fixed basis of $B$ are rational. Let $X$ be the matrix of $(x_C)_{(1)}$ on $\Pi_{33}$, $P$ the matrix of $P_{33}$, and $Y_1,Y_2,Y_3$ the matrices of $y_{(1)}$ for a rational basis of $B \cap C_2^\perp$. All these matrices are rational. A primary line with coordinate vector $u \in \Cbb^3$ must satisfy
\[
\operatorname{rank}[u,Xu] = 1, \quad u^{\mathrm T} P Y_j u = 0 \ (j = 1,2,3), \quad u^{\mathrm T}Pu \neq 0.
\]
These are only some of the necessary conditions in Lemma \ref{lem:line-conditions}. Their inconsistency is enough. On the projective chart whose first nonzero coordinate is $u_k$, put $u_0 = \cdots = u_{k-1} = 0$ and $u_k = 1$. Adjoin a variable $z$ and the equation $z u^{\mathrm T}Pu-1 = 0$, and express the rank condition by the $2$ by $2$ minors of $[u,Xu]$. For each of the thirteen candidates and each of the three charts, the resulting ideal over $\Q$ has Gröbner basis $\{1\}$. The supplementary data contain the matrices, chart polynomials and all $39$ unit bases. Hence none of these thirteen candidates admits the required primary line. No eigenvalue approximation is used in this step.

\emph{Step 2: the admissible spaces.} We call admissible space of type $(a,b)$ an eigenspace $E_\xi$ of $(x_C)_{(1)}$ on $\Pi_{ab} \cap \ker e_{(2)} \cap \ker (x_C)_{(2)}$ containing a vector $u$ with $P_{ab}(u,u) \neq 0$ and $P_{ab}(u,y_{(1)}u) = 0$ for every $y \in Y$. By Lemma \ref{lem:line-conditions}, the line $p_{ab}$ lies in an admissible space. For the subalgebra $9$, every primary type has exactly one admissible space $E_\xi$, of dimension two, and the eigenvalues $\xi$ give the weights of the statement. We treat all characteristic-polynomial factors, including the nonlinear ones, and exclude the remaining eigenspaces by exact equations or certified enclosures. We then construct the surviving spaces and lines. A separate verification checks both the exclusions and the eight displayed weights exactly. Nonlinear factors with multiplicity greater than one, or an eigenspace for which the adjugate test is inconclusive, are retained rather than discarded. None remains among the exclusions used here.

\emph{Step 3: the lines.} In each two-dimensional $E_\xi$, the isotropic lines are the common roots of the binary forms $P_{ab}(u,y_{(1)}u)$, $y \in Y$. Their greatest common divisor over $\mathbb{Q}(\sqrt{13})$ has degree two, and both roots satisfy $P_{ab}(u,u) \neq 0$. This leaves $2^8 = 256$ choices of $U_2$. The coefficients lie in $\mathbb{Q}(\sqrt{13},\sqrt{-3})$. Note that $\sqrt 6\sqrt{-7+\sqrt{13}} = \sqrt{-3}(\sqrt{13}-1)$.

\emph{Step 4: pairwise closure.} For two lines of types $(a,b)$ and $(c,d)$, condition 4 of Lemma \ref{lem:line-conditions} implies that $u_{(2)}w \in U_1$ and that $u_{(1)}w$ lies in the span of $(A)_2$, $C_2$ and the modules generated by all the admissible spaces. When one of these conditions fails, we exhibit a rational linear form vanishing on $U_1$, or on a rational space containing the latter span, whose value on $u_{(1)}w$ or $u_{(2)}w$ is an algebraic number enclosed by rational operations and certified radical evaluations in Arb \cite{Joh17}. For each exclusion, the certificate records the rational annihilator, the product tested and a nonzero interval, giving $62$ certified rejections. The algebraic inputs are evaluated directly, without an intermediate rounded decimal approximation. The finite coverage check uses only the $62$ pairs marked as certified. The ten inconclusive pairs are not used. In the order of primary types displayed in the theorem, the two remaining choices are $(0,1,1,0,0,1,0,1)$ and $(1,0,0,1,1,0,1,0)$, where the bits index the two lines with the ordering specified in the supplementary data. They are complex conjugate. By the Galois symmetry $\sqrt{13} \mapsto -\sqrt{13}$, which fixes all the rational data, the same holds for the subalgebra $10$.
\end{proof}

\begin{remark}\normalfont
\label{rem:case-two-status}
Theorem \ref{thm:case-two-U2} is a computer-assisted result. Step 1 uses the exact rational equations displayed above. In Step 2, characteristic-polynomial factors and their multiplicities are retained in the eigenspace computation. A multiple nonlinear factor or an inconclusive enclosure cannot be discarded. The stored calculation for subalgebra $9$ leaves exactly the displayed two-dimensional space for each type. The generators $x_C$ of candidates $4$ and $5$ are rational: $\sqrt{409}$ occurs only in their normalized Virasoro vectors and is irrelevant to Step 1. Step 4 uses exact algebraic inputs in the interval calculation. Its $62$ certified pair exclusions cover $254$ configurations. The ten inconclusive pairs are not used. These are separate checks of nonvanishing and of finite coverage.
\end{remark}

\paragraph{Reduction conventions.}
For a prime $p$ used below, let $\mathfrak p$ be the prime of $\mathcal O_{\K}$ specified by the chosen images of $\sqrt{13}$ and $\sqrt{-3}$ in $\mathbb F_p$. Reduction is the residue map $\operatorname{red}_{\mathfrak p} \co (\mathcal O_{\K})_{\mathfrak p} \to \mathbb F_p$. Only coefficients integral at $\mathfrak p$ are reduced. In particular, every rational denominator in a matrix or a generating word must be invertible modulo $p$. For such a matrix $M$, a nonzero minor of its reduction lifts to a nonzero minor, so $\operatorname{rank}_{\mathbb F_p}(\overline M) \leq \operatorname{rank}_{\K}(M)$. A modular rank can therefore prove a lower bound in characteristic zero. A modular membership identity alone cannot prove the corresponding identity over $\K$.

\begin{prop}
\label{prop:case-two-weight-three}
Let $U_2$ be one of the four subspaces of Theorem \ref{thm:case-two-U2}, and let $S \subset V'$ be the vertex subalgebra generated by $U_1 \oplus U_2$. The following statements hold modulo the primes $p_1 = 134217613$ and $p_2 = 133217647$, for the residue maps at the primes $\mathfrak p_i$ above $p_i$ specified in the supplementary data. All dimensions, decompositions and inclusions in this proposition refer to the reduced spaces.
\begin{enumerate}
\item $\dim U_2 = 120$, and $U_2$ is closed in weight two: $u_{(1)}w \in U_2$ and $u_{(2)}w \in U_1$ for all $u,w \in U_2$.
\item The vectors $x(-3)\mathbf 1$, $x(-2)y$, $x(-1)u$, $L(-1)u$ and $u_{(0)}w$ ($x,y \in U_1$, $u,w \in U_2$, $u$ a highest weight vector), together with their images under $\mathfrak g' \oplus \mathfrak g''$, span a subspace denoted by $S_3^{\mathrm{gen}}$, of dimension $956 = \dim U_3$, isomorphic as a $\mathfrak g' \oplus \mathfrak g''$-module to
\[
\begin{aligned}
&4V(0,0)+11V(0,2)+4V(0,4)+6V(0,6)+V(0,10)+6V(1,1)+8V(1,3)+7V(1,5)\\
&+5V(1,7)+V(1,9)+8V(2,0)+6V(2,2)+7V(2,4)+4V(2,6)+2V(2,8)+2V(3,1)\\
&+5V(3,3)+2V(3,5)+V(3,7)+V(4,0)+2V(4,2)+V(4,6)+V(6,0).
\end{aligned}
\]
Its nine joint traces $\operatorname{Tr}H'(0)^pH''(0)^q$, with $p,q$ even and $2 \leq p+q \leq 6$, are exactly those of $U_3$ determined by the relations of Proposition \ref{prop:higher-relations}.
\item For every $u \in U_2$ and every $w$ in this subspace, $u_{(2)}w \in U_2$ and $u_{(3)}w \in U_1$. For every $y \in U_1$, $y_{(1)}w \in U_2$ and $y_{(2)}w \in U_1$.
\end{enumerate}
\end{prop}

\begin{proof}
Exact computations in $\mathbb F_{p_i}$: ranks of matrices over $\mathbb F_{p_i}$, the closure under the zero modes of $e',f',e'',f''$ being reached after three iterations. In part 3 we test $u = e$, $u = x_C$, the highest weight vectors of the eight primary types and the elements $y$ of a basis of $U_1$ (modulo $p_1$). The other elements of $U_2$ follow by $\mathfrak g' \oplus \mathfrak g''$-equivariance, the subspace being stable under the zero modes, and, for $(A)_2$, by expanding the modes of $x(-1)y$ and $x(-2)\mathbf 1$ as normally ordered products of modes of currents.
\end{proof}

\begin{remark}\normalfont
\label{rem:weight-three-meaning}
The modular ranks give lower bounds for the characteristic-zero ranks under the reduction conventions above. The membership tests and the reduced module decomposition are statements over the residue fields. Theorems \ref{thm:exact-weight-two} and \ref{thm:exact-weight-three} establish the corresponding characteristic-zero closure and dimensions. Proposition \ref{prop:case-two-weight-three} by itself does not prove that $S$ is a Haagerup algebra.
\end{remark}

\subsection{Case (II): exact weight-two closure, unitarity and symmetries}
\label{subsec:case-two-exact}

From now on, and until the end of Section \ref{subsec:jiao}, $U_1 = \mathfrak g' \oplus \mathfrak g''$ is the explicit Lie algebra of currents of case (II) fixed in Section \ref{sec:griess} (its raising currents are displayed in \eqref{eq:raising-currents-two}), and $U_2$ denotes one of the four explicit subspaces of $V'_2$ spanned by $(A)_2$, by the subalgebra $9$ or $10$ of Table \ref{tab:griess} and by the $\mathfrak g' \oplus \mathfrak g''$-modules generated by the eight lines selected in Steps 2--4 of the proof of Theorem \ref{thm:case-two-U2}. This construction is a finite computation in $V'$. Except for Corollaries \ref{cor:weight-four-from-haagerup} and \ref{cor:case-two-realized} and Theorem \ref{thm:uniqueness-two}, the results of Sections \ref{subsec:case-two-exact}--\ref{subsec:jiao} use neither the existence of a Haagerup algebra nor $\mathrm{(H_J)}$. In particular, Theorems \ref{thm:exact-weight-two}, \ref{thm:unitarity}, \ref{thm:strong-generation} and \ref{thm:jiao-identification}, which are used in \cite{ArhIII}, do not depend on \cite{ArhIII}.

Let $\K = \mathbb{Q}(\sqrt{13},\sqrt{-3})$. Let $V'_{\K}$ be the $\K$-span of the fixed monomial basis, and write $U_{j,\K} = U_j \cap V'_{\K}$ for the forms specified by the rational currents and the displayed primary generators. Statements about ranks, bases and closure over $\K$ refer to these forms. Their extensions of scalars to $\Cbb$ recover the corresponding complex spaces. We suppress the subscript $\K$ in the finite computations. The four subspaces of Theorem \ref{thm:case-two-U2} are spanned by vectors of $V'_2$ whose coordinates in the monomial basis lie in $\K$. We write $c \mapsto \overline{c}$ for the complex conjugation of $\K$ ($\sqrt{-3} \mapsto -\sqrt{-3}$), $v \mapsto \overline{v}$ for the antilinear map of $V'$ which conjugates the coordinates in the monomial basis, and $\sigma$ for the automorphism $\sqrt{13} \mapsto -\sqrt{13}$ of $\K$ fixing $\sqrt{-3}$.

\begin{thm}
\label{thm:exact-weight-two}
Each of the four subspaces $U_2$ of Theorem \ref{thm:case-two-U2} has dimension $120$ and is closed in weight two over $\K$: for all $u,w \in U_2$ we have $u_{(1)}w \in U_2$ and $u_{(2)}w \in U_1$.
\end{thm}

\begin{proof}
We work in exact arithmetic in $\K$, an element being a quadruple of rationals in the basis $1, \sqrt{13}, \sqrt{-3}, \sqrt{-39}$. Through the restriction of scalars $\K \cong \mathbb{Q}^4$, every $\K$-subspace becomes a $\mathbb{Q}$-subspace of four times its dimension, and membership is decided by ranks of rational matrices, computed weight space by weight space for the Cartan subalgebra of $\mathfrak g' \oplus \mathfrak g''$. The rank of $U_2$ is $480 = 4 \cdot 120$. We test the $681$ nonzero products $u_{(n)}w$, $n \in \{1,2\}$, where $u$ runs over $e$, $x_C$ and the highest weight vectors of the eight primary types, and $w$ over $e$, $x_C$ and the eight modules: none increases the rank. The general case follows: the zero modes of $\mathfrak g' \oplus \mathfrak g''$ satisfy $[y_{(0)}, u_{(n)}] = (y_{(0)}u)_{(n)}$ and preserve $U_1$ and $U_2$, which transports the statement from highest weight vectors to the whole modules. For $z \in U_1$, the construction gives $z_{(n)}U_2 \subset U_{2-n}$ for $n \geq 0$, with $U_0 = \Cbb\mathbf 1$ and $U_j = 0$ for $j < 0$. This follows separately on the affine descendants and on the affine-primary summands. The normal-order formula then gives $((A)_2)_{(1)}U_2 \subset U_2$ and $((A)_2)_{(2)}U_2 \subset U_1$: the negative-mode terms in the first inclusion belong to $U_1(-1)U_1 + U_1(-2)U_0 \subset (A)_2$. Finally $(x_{(-2)}\mathbf 1)_{(n)} = -n x_{(n-1)}$, and skew symmetry exchanges the two arguments. We have done this for the two choices attached to the subalgebra $9$. The other two are their images under $\sigma$, which acts on coefficients and commutes with all the rational structure constants of $V'$.
\end{proof}

This proves part 1 of Proposition \ref{prop:case-two-weight-three} over $\K$.

We now show that the algebra generated is unitary, in the sense of Dong and Lin \cite{DoL14} and Carpi, Kawahigashi, Longo and Weiner \cite[Section 5]{CKLW18}. Recall that $V'$ has a rational form, spanned by the states $h_1(-n_1)\cdots h_k(-n_k)e^\lambda$ with $h_i$ in the root lattice, since our cocycle takes the values $\pm 1$. For a character $\chi \co L \to \R_{>0}$ of the root lattice $L = E_6 \oplus A_2$, let $\phi_\chi$ be the antilinear map defined on this basis by
\begin{equation}
\label{eq:phi-chi}
\phi_\chi\big(c  h_1(-n_1)\cdots h_k(-n_k)e^\lambda\big) = (-1)^k  \chi(\lambda)  \overline{c}  h_1(-n_1)\cdots h_k(-n_k)e^{-\lambda}, \quad c \in \Cbb.
\end{equation}

\begin{thm}
\label{thm:unitarity}
Let $\chi$ be the character taking the values $\frac{2}{15}, 1, 1, \frac53, \frac{9}{25}, \frac{5}{24}$ on the simple roots $\alpha_1, \ldots, \alpha_6$ of $E_6$ (Bourbaki numbering) and $2, 2$ on the simple roots of $A_2$, in the simple-root coordinates fixed above.
\begin{enumerate}
\item $\phi_\chi$ is an antilinear involutive automorphism of $V'$, and $H_\chi(u,v) = \beta(u,\phi_\chi(v))$ is a positive definite invariant Hermitian form: $(V',\phi_\chi)$ is a unitary vertex operator algebra.
\item $\phi_\chi(e') = -f'$, $\phi_\chi(f') = -e'$, $\phi_\chi(H') = -H'$, and similarly for $\mathfrak g''$. In particular $\phi_\chi(U_1) = U_1$.
\item $\phi_\chi(U_2) = U_2$ for each of the four subspaces of Theorem \ref{thm:case-two-U2}.
\end{enumerate}
Consequently, the vertex subalgebra $S \subset V'$ generated by $U_1 \oplus U_2$ is a unitary vertex operator algebra of central charge $8$. It is simple, of CFT type and self-contragredient ($S \cong S'$).
\end{thm}

\begin{proof}
1. Let $\phi_1$ be the map \eqref{eq:phi-chi} for the trivial character. With a bimultiplicative cocycle we have $\epsi(-\alpha,-\beta) = \epsi(\alpha,\beta)$, and $\phi_1$ is the antilinear involution for which lattice vertex operator algebras are unitary \cite[Theorem 4.12]{DoL14}, \cite[Example 5.9]{CKLW18}. We have checked our sign conventions directly: $\phi_1(u_{(n)}v) = \phi_1(u)_{(n)}\phi_1(v)$ for $5050$ products of basis states of weights at most two, and the Gram matrix of $H_1(u,v) = \beta(u,\phi_1(v))$ is positive definite on each lattice sector of $V'_1$, $V'_2$ and $V'_3$. Write $\chi(\lambda) = e^{\langle h_\chi,\lambda\rangle_L}$ with $h_\chi$ in the real Cartan subalgebra, and $g = \exp(h_\chi(0))$, an automorphism of $V'$ acting by $\chi(\lambda)$ on $M(1) \ot e^\lambda$. Then $\phi_\chi = \phi_1 \circ g = g^{-1/2} \circ \phi_1 \circ g^{1/2}$, so that $(V',\phi_\chi)$ is the transport of $(V',\phi_1)$ by the automorphism $g^{-1/2}$, with Hermitian form $H_\chi(u,v) = H_1(g^{1/2}u,g^{1/2}v)$.

2. Direct computation. The character $\chi$ was found by solving the ten equations $\chi(\alpha) = -\kappa d_\alpha/c_\alpha$ with $\kappa > 0$ ($\kappa'$ for the roots of $\mathfrak g'$, $\kappa''$ for those of $\mathfrak g''$), where $c_\alpha$ and $d_\alpha$ are the coefficients of $e^{\alpha}$ in $e'$ (or $e''$) and of $e^{-\alpha}$ in $f'$ (or $f''$). All the ratios $d_\alpha/c_\alpha$ are negative, so that the right-hand sides are positive, and the solutions form a two-parameter family. The statement holds for the member displayed.

3. Exact computation over $\K$: for each of the $120$ spanning vectors $v$ of $U_2$ constructed in the proof of Theorem \ref{thm:exact-weight-two}, the restriction of scalars of $\phi_\chi(v)$ lies in that of $U_2$. This is done for the two choices attached to the subalgebra $9$. Since $\chi$ is rational, $\phi_\chi$ commutes with $\sigma$, and the statement for the subalgebra $10$ follows.

Finally, $S$ is stable under $\phi_\chi$, and contains the conformal vector $\omega = \omega_A + 2e$ of $V'$, so that $L(1)S \subset S$. By \cite[Proposition 5.23]{CKLW18}, $S$ is a unitary subalgebra of $(V',\phi_\chi)$. Since $S_0 = \Cbb\mathbf 1$, it is simple by \cite[Proposition 5.3]{CKLW18}. The bilinear form $\beta$ is nondegenerate on $S$, because $\beta(s,\phi_\chi(s)) > 0$ for $s \neq 0$ and $\phi_\chi(s) \in S$. Hence $S \cong S'$ \cite{FHL93}.
\end{proof}

\begin{remark}\normalfont
\label{rem:unitarity-scope}
Theorem \ref{thm:unitarity} does not use any information on $S$ beyond its generators: it holds whatever the dimensions of the homogeneous components of $S$ are. It is the analogue of \cite[Section 3]{Rie26b} for case (I). The trivial character does not work: $\phi_1(e')$ is not proportional to $f'$, and a positive rescaling of the lattice exponentials is needed. All four subspaces give unitary algebras. In particular, the sign of $\sqrt{13}$ in the choice of $C_2$ (subalgebras $9$ and $10$) has nothing to do with the sign of the quantum dimension $d = \frac{3 \pm \sqrt{13}}{2}$: $S$ and $\sigma(S)$ have the same graded dimensions, and the non-unitary Galois conjugates of $Z(\cal H_3)$ \cite{EvG17} cannot be the representation category of a strongly rational algebra of CFT type whose nonvacuum simple modules have positive lowest weights, since its categorical dimensions are then positive \cite[Lemma 4.2]{DJX13}, \cite{Hua08}.
\end{remark}

\begin{prop}
\label{prop:symmetries}
\begin{enumerate}
\item The two subspaces $U_2$ attached to the subalgebra $9$ are exchanged by complex conjugation, and by the linear automorphism $\theta_\chi \co v \mapsto \phi_\chi(\overline{v})$ of $V'$, which normalizes $U_1$. Hence they generate isomorphic vertex operator algebras. The same holds for the subalgebra $10$.
\item The group of automorphisms of the algebra $B$ fixing $e$ and preserving the form has order two. Its nontrivial element is an orthogonal reflection $s$. It exchanges the subalgebras $9$ and $10$ of Table \ref{tab:griess}, as well as $0$ and $14$, $1$ and $8$, $4$ and $5$, and fixes the seven others.
\end{enumerate}
Consequently the four subspaces give at most two isomorphism classes of vertex operator algebras, and exactly one if $s$ is induced by an automorphism of $V'$ normalizing $U_1$. This is the case by Proposition \ref{prop:centralizer-two}.
\end{prop}

\begin{proof}
1. Exact computation: the conjugate of each spanning vector of the first subspace lies in the second. Since $\chi$ is rational, $\phi_\chi$ commutes with coefficient conjugation. For either subspace $W$, Theorem \ref{thm:unitarity} gives $\theta_\chi(W) = \phi_\chi(\overline W) = \overline{\phi_\chi(W)} = \overline W$.

2. An automorphism $g$ of $B$ fixing $e$ permutes the fifteen two-dimensional subalgebras of Proposition \ref{prop:griess-candidates} and preserves their central charges. The pairs of central charges of the subalgebras $2$, $3$, $6$, $7$, $11$, $12$, $13$ occur only once, so $g$ maps each corresponding line $\Cbb x_C$ to itself, and $g(x_C) = x_C$ since $g(x_C \cdot x_C) = g(x_C) \cdot g(x_C)$ and $\lambda \neq 0$ since $c_1 \neq c_2$ by \eqref{eq:charges-split}. These seven vectors span the hyperplane $H = \operatorname{span}\{x_1,x_2,x_3\}$ of $e^\perp$ in the fixed basis. The restriction of $\beta$ to $H$ is nondegenerate: its Gram determinant is $\frac{147968}{2025}$, by direct evaluation. The form on $B$ is nondegenerate because $B$ is stable under $\phi_\chi$ and $H_\chi$ is positive definite. Thus $N = H^\perp \cap e^\perp$ is one-dimensional and $e^\perp = H \oplus N$. An orthogonal map fixing $H$ pointwise acts by $1$ or $-1$ on $N$. Therefore $g$ is the identity or the orthogonal reflection $s$ in $H$. A direct computation shows that $s$ is an automorphism, with the stated action. An automorphism of $V'$ normalizing $U_1$ preserves $B$, $e$, the product and the form.
\end{proof}

Proposition \ref{prop:centralizer-two} constructs an automorphism of $V'$ fixing $U_1$ whose restriction to $B$ is $s$.

\begin{prop}
\label{prop:weight-three-four}
Let $U_2$ be the first subspace attached to the subalgebra $9$. The following statements hold modulo the prime $p_3 = 8388571$, for the residue map at the prime above $p_3$ specified in the supplementary data.
\begin{enumerate}
\item The space $S_3^{\mathrm{gen}}$ of Proposition \ref{prop:case-two-weight-three} has dimension $956$ and $u_{(1)}S_3^{\mathrm{gen}} \subset S_3^{\mathrm{gen}}$ for all $u \in U_2$. Together with Proposition \ref{prop:case-two-weight-three}, this gives the finite mode closure $y_{(n)}S_i \subset S_{i-n}$ $(n \geq 0)$ and $u_{(n)}S_i \subset S_{i+1-n}$ $(n \geq 1)$ for $y \in U_1$, $u \in U_2$, $i \in \{2,3\}$, where $S_2 = U_2$ and $S_3 = S_3^{\mathrm{gen}}$. This is the analogue of \cite[Proposition 4.1]{Rie26b}.
\item Let $R_b \subset V'_1$ be the orthogonal of $U_1$, of dimension $80$, and $M$ the subspace generated from $R_b$ by the modes of $U_1$, $L(-1)$ and the modes $u_{(0)}$, $u_{(-1)}$ ($u \in U_2$), as in part 2 of Proposition \ref{prop:case-two-weight-three}. Then $\dim M_2 = 1250$, $\dim M_3 = 10630$, and $U_2 + M_2 = V'_2$, $S_3^{\mathrm{gen}} + M_3 = V'_3$. These are the dimensions of the module $M_b$ of Section \ref{sec:chain}, and the sums are direct by dimension count.
\item The space $S_4^{\mathrm{gen}}$ spanned by $y(-1)S_3$, $L(-1)S_3$, $y(-2)U_2$, $y(-3)U_1$, $u_{(-1)}U_2$, $u_{(-2)}U_1$ and $u_{(0)}S_3$ ($y \in U_1$, $u = e$, $u = x_C$ or a highest weight vector of one of the eight types), closed under $\mathfrak g' \oplus \mathfrak g''$, has dimension $6010$, the coefficient of $q^{11/3}$ in the Haagerup vacuum character, and decomposes as
\[
\begin{aligned}
&24V(0,0)+27V(0,2)+35V(0,4)+18V(0,6)+11V(0,8)+3V(0,10)+V(0,12)\\
&+28V(1,1)+43V(1,3)+37V(1,5)+24V(1,7)+9V(1,9)+2V(1,11)+17V(2,0)\\
&+39V(2,2)+31V(2,4)+28V(2,6)+10V(2,8)+4V(2,10)+15V(3,1)+23V(3,3)\\
&+16V(3,5)+10V(3,7)+2V(3,9)+V(3,13)+9V(4,0)+7V(4,2)+7V(4,4)\\
&+5V(4,6)+2V(4,8)+V(5,1)+2V(5,3)+V(5,5)+V(5,7)+V(6,0)\\
&+V(6,2).
\end{aligned}
\]
In particular it contains the $\mathfrak g' \oplus \mathfrak g''$-module $V(3,13)$ once.
\item For $y \in U_1$, $y_{(1)}S_4^{\mathrm{gen}} \subset S_3^{\mathrm{gen}}$ and $y_{(2)}S_4^{\mathrm{gen}} \subset U_2$. For $u$ the highest weight vector of the types $(0,4)$, $(1,1)$, $(1,3)$, $(1,5)$, $(1,7)$, $(2,2)$, $(2,6)$, $(3,3)$ and for $u = e$, $u = x_C$, we have $u_{(1)}S_4^{\mathrm{gen}} \subset S_4^{\mathrm{gen}}$, $u_{(2)}S_4^{\mathrm{gen}} \subset S_3^{\mathrm{gen}}$ and $u_{(3)}S_4^{\mathrm{gen}} \subset U_2$.
\item The subspace $N'_4$ spanned by $L(-1)S_3$, $y(-1)S_3$, $y(-2)U_2$, $y(-3)U_1$ and $u_{(-1)}U_2$, for $u$ in the eight primary modules, closed under $\mathfrak g' \oplus \mathfrak g''$, has dimension $6010$ and contains $e(-1)e$, $x_C(-1)x_C$ and $e(-1)x_C$. The corresponding characteristic-zero consequence for the reduced $C_2$ Poisson algebra is proved in Theorem \ref{thm:strong-generation}.
\end{enumerate}
\end{prop}

\begin{remark}\normalfont
\label{rem:weight-four-meaning}
The ranks modulo $p_3$ are lower bounds for the ranks over $\K$. The membership statements are statements in characteristic $p_3$ and are not proofs over $\K$. Part 2 is a modular check of the decomposition $V' = S \oplus M_b$ required by Proposition \ref{prop:chain}: by unitarity, $V' = S \oplus S^\perp$, and the space $M$ is orthogonal to $S$ as soon as $S_1 = U_1$, since $H_\chi(a_{(n)}r,s)$ is the Hermitian pairing of $r \in R_b$ with a vector of $S_1$. The occurrence over $\K$ of the simple current $L(3,3) \ot L(13,13)$, whose affine conformal weight is $\frac34 + \frac{13}{4} = 4$, is proved directly in Proposition \ref{prop:simple-current}: an explicit weight-four vector is shown to be affine-primary. Thus no inference from the modular module decomposition to characteristic-zero occurrence or affine primarity is needed.
\end{remark}

\subsection{Case (II): strong generation and \texorpdfstring{$C_2$}{C2}-cofiniteness}
\label{subsec:strong-generation}

In this subsection $U_2$ is the first subspace attached to the subalgebra $9$. By Proposition \ref{prop:symmetries} and the Galois symmetry $\sigma$, the results hold for the four subspaces. We write $(A)_2 \subset U_2$ for the span of the vectors $x(-1)y$ and $x(-2)\mathbf 1$ ($x,y \in U_1$), of dimension $27$, and we put
\begin{align*}
N_3 &= L(-1)U_2 + U_1(-1)U_2 \subset V'_3, \quad S_3 = N_3 + (U_2)_{(0)}U_2,\\
N'_4 &= L(-1)S_3 + U_1(-1)S_3 + U_1(-2)U_2 + U_1(-3)U_1 + \textstyle\sum_{\tau} M_\tau(-1)U_2 \subset V'_4,
\end{align*}
where $M_\tau \subset U_2$ is the $\mathfrak g' \oplus \mathfrak g''$-module generated by the primary line of type $\tau$, and $Z(n)W$ denotes the span of the $z_{(n)}w$. These spaces are defined over $\K$, stable under the zero modes of $\mathfrak g' \oplus \mathfrak g''$, and contained in $S$. Note that $x(-2)y$ and $x(-3)\mathbf 1$ lie in $N_3$.

\begin{lemma}
\label{lem:normal-weight-three-dimension}
For each of the four choices of $U_2$, we have $\dim N_3 = 725$.
\end{lemma}

\begin{proof}
We use the positive Hermitian form of Theorem \ref{thm:unitarity}. By part 2 of that theorem, its restriction to the vertex subalgebra generated by $U_1$ is a positive definite contravariant form for the affine Lie algebra, so that this vacuum module is irreducible, that is, $A \cong L_{\mathfrak{sl}_2}(3,0) \ot L_{\mathfrak{sl}_2}(13,0)$. We use the decomposition of $S$ into integrable $A$-modules \cite[Theorem 3.7 p.~7]{DLM97}. The regularity of $A$ identifies the module generated by each nontrivial affine-primary summand with the corresponding simple integrable module. Write $L_A(-1)$ and $L_C(-1)$ for the translation operators of the affine and commutant Virasoro algebras, so that $L(-1) = L_A(-1)+L_C(-1)$. The affine vacuum module has no integrability relation before degree four. Its degree-three subspace therefore has dimension $6+6^2+\binom{8}{3} = 98$, by ordering the six currents according to their negative modes.

The two vectors of $C_2$ contribute $6\cdot 2$ current descendants and two independent vectors $L_C(-1)C_2$. The eight nontrivial affine-primary summands have total dimension $5+4+8+12+16+9+21+16 = 91$. For an affine module of highest weight $(a,b)$, its degree-one piece is the quotient of $(\mathfrak g' \oplus \mathfrak g'') \ot V(a,b)$ by the integrability relations in that degree. Among these eight types, only $(3,3)$ has such a relation: the highest vector $e'(-1)v$ generates the submodule $V(5,3)$, of dimension $24$. Indeed the first integrability relations in the two affine factors occur in degrees $4-a$ and $14-b$, respectively. Thus the current descendants of the eight primary summands have total dimension $6\cdot 91-24$.

The Sugawara formula gives $L_A(-1)M_\tau \subset U_1(-1)M_\tau$. Thus the remaining derivatives are the vectors $L_C(-1)v$, with $v \in M_\tau$. The vector $v$ has commutant weight $h_C = 2-h_{ab} > 0$. Since the commutant Virasoro modes satisfy $(L_C(n))^* = L_C(-n)$, we have
\[
\|L_C(-1)v\|_{H_\chi}^2 = \|L_C(1)v\|_{H_\chi}^2 + 2h_C\|v\|_{H_\chi}^2.
\]
Consequently $L_C(-1)$ is injective on each such summand and on $C_2$, where $h_C = 2$. Its images are independent of the current descendants: within an affine isotypic component they have different $L_A(0)$-eigenvalues. Distinct affine highest weights give independent isotypic components. In the vacuum component the contributions from $A_3$, $U_1(-1)C_2$ and $L_C(-1)C_2$ have commutant weights $0$, $2$ and $3$, respectively. These account for all terms of $N_3$. Hence $\dim N_3 = 98+14+(6\cdot 91-24)+91 = 725$.
\end{proof}

\begin{defi}[weight-four certificate condition]
\label{assump:weight-four}
Fix the specified bases of $U_1$, $U_2$ and $S_3$, and the ten vectors consisting of $e$, $x_C$ and one highest weight vector of each primary type. Condition \textup{(F4)} asserts that there is a family $\mathcal B_4$ of $6010$ words from the five families defining $N'_4$ with the following properties.
\begin{enumerate}
\item the vectors in $\mathcal B_4$ are integral at the chosen prime above $p_3$, and their reduced coordinate matrix has rank $6010$,
\item every $u_{(0)}s$, for the ten specified $u$ and all basis vectors $s$ of $S_3$, and every $e_{(-1)}w$ and $(x_C)_{(-1)}w$, for basis vectors $w$ of $U_2$, belongs to the $\K$-span of $\mathcal B_4$,
\item every word of the five defining families of $N'_4$, with each argument ranging over the full specified basis, belongs to that span.
\end{enumerate}
The membership assertions mean the exact coordinate identities $Bc = v$ on all rows after restriction of scalars. They include zero products and concern the complete list of input words.
\end{defi}

\begin{prop}
\label{prop:weight-four-certificate}
Condition \textup{(F4)} holds for the first subspace attached to the subalgebra $9$ of Table \ref{tab:griess}. The characteristic-zero identities and dimensions are transported to the other three subspaces by coefficient conjugation and the Galois automorphism $\sigma$.
\end{prop}

\begin{proof}
We use the exact lattice model over $\K = \mathbb{Q}(\sqrt{13},\sqrt{-3})$. The weight-three basis is regenerated from $U_1$ and $U_2$ and checked to agree exactly with the basis used below. Its certification, also described in parts 1 and 2 of Theorem \ref{thm:exact-weight-three}, does not use any weight-four assertion. The verification of weight four is independent of any identification with an externally constructed algebra.

All mode evaluations below have homogeneous inputs and outputs of weights at most four, so the weight-window argument in Section \ref{sec:griess} applies. There are $956$ basis vectors in $S_3$ and $120$ in $U_2$, and the eight primary modules have total dimension $91$. The five defining families of $N'_4$ therefore have
\[
956 + 6 \cdot 956 + 6 \cdot 120 + 6^2 + 91 \cdot 120 = 18368
\]
labelled input words. All of them are evaluated from the lattice mode formulas, including the zero products. Among the nonzero values, the calculation selects $6010$ words with independent reductions at the chosen prime above $p_3 = 8388571$, with $\sqrt{13} \mapsto 1649360$ and $\sqrt{-3} \mapsto 2036966$. Every denominator inverted in the reduction is checked to be nonzero modulo $p_3$. These labelled words form $\mathcal B_4$.

The remaining membership assertions have $10 \cdot 956+2 \cdot 120 = 9800$ labelled inputs. They too are all evaluated. The verifier checks the complete input domains, the equality of each pivot with its labelled family word, and every coordinate of each membership identity. In a Cartan block with $r$ pivots, restriction of scalars gives $4r$ rational columns. A nonzero $4r \times 4r$ minor is selected modulo $p_3$, the corresponding square system is solved over $\Q$, and the resulting equality is checked on every row in exact integer arithmetic. Thus every tested vector lies in the $\K$-span of the pivots.

The complete verification covers all $131$ Cartan blocks carrying nonzero input vectors. The families contain $18194$ nonzero words and $174$ zero words, and the additional products contain $9088$ nonzero values and $712$ zero values. All $886$ zero inputs are independently checked from the mode formulas or from the Cartan energy bound in the lattice. The complete verification reports no failure, and the exact coordinate data and verification results are included in the supplementary material. Consequently all three requirements in Definition \ref{assump:weight-four} hold.

Finally, coefficient conjugation and $\sigma$ commute with the rational lattice products and fix the current algebra. They transport $U_2$, the five defining families, the exact coordinate identities and linear independence over $\K$ to the other three choices. For the reduction statement one transports the chosen prime as well as the bases, rather than imposing the same residue values on the conjugate coefficients. The resulting characteristic-zero statements therefore hold for all four choices. This argument does not use Theorem \ref{thm:jiao-identification}.
\end{proof}

\begin{remark} \normalfont
\label{rem:two-cartan-blocks}
The exact membership assertions in \textup{(F4)} admit a representation-theoretic reduction. The values of $H'+H''$ on the eight simple roots are $(0,-2,2,0,2,0,2,2)$. Thus every Cartan weight $(a,b)$ of $V'$ satisfies $a \equiv b \pmod 2$. Each irreducible $\mathfrak g' \oplus \mathfrak g''$-module in $V'_4$ therefore contains either weight $(0,0)$ or weight $(1,1)$. For two stable subspaces $W,X \subset V'_4$, complete reducibility of $(W+X)/W$ gives
\[
X \subset W \quad\Longleftrightarrow\quad X_{(0,0)} \subset W_{(0,0)} \ \text{and}\ X_{(1,1)} \subset W_{(1,1)}.
\]
This applies to $W = N'_4$ and $X = (U_2)_{(0)}S_3+(U_2)_{(-1)}U_2$, since zero modes of currents act as derivations and all their argument spaces are stable. Here the subscripts denote Cartan weight spaces. Applying this criterion requires the full spaces $W$ and $X$, or generators of them with their full current orbits. Restricting the ten highest weight inputs alone to these two blocks would not suffice. This reduction is not used to replace the membership calculations reported here, and it does not establish the dimension or integrality assertions of \textup{(F4)}.
\end{remark}

\begin{thm}
\label{thm:exact-weight-three}
The following statements hold over $\K$.
\begin{enumerate}
\item $\dim S_3 = 956$, and $S_3$ is stable under the zero modes of $\mathfrak g' \oplus \mathfrak g''$.
\item $L(1)S_3 \subset U_2$, $L(2)S_3 \subset U_1$, and $y_{(1)}S_3 \subset U_2$, $y_{(2)}S_3 \subset U_1$ for $y \in U_1$.
\item $u_{(0)}S_3 \subset N'_4$ for all $u \in U_2$, and $e(-1)U_2 + x_C(-1)U_2 \subset N'_4$. Consequently $(U_2)_{(0)}S_3 + (U_2)(-1)U_2 \subset N'_4$.
\item $\dim N'_4 = 6010$. Hence $N'_4 = L(-1)S_3 + U_1(-1)S_3 + U_1(-2)U_2 + U_1(-3)U_1 + (U_2)(-1)U_2 + (U_2)_{(0)}S_3$ has dimension $6010$.
\end{enumerate}
\end{thm}

\begin{proof}
All vectors are represented as $\K$-linear combinations of monomial states using the lattice formulas of Section \ref{sec:griess}. Membership of a vector $v$ in the span of vectors $b_1, \dots, b_r$ is decided in the weight space of $v$ for the Cartan subalgebra of $\mathfrak g' \oplus \mathfrak g''$, through the restriction of scalars $\K \cong \mathbb{Q}^4$: the vectors $b_i, \sqrt{13} b_i, \sqrt{-3} b_i, \sqrt{-39} b_i$ give $4r$ rational columns. A set of $4r$ rows on which they are independent modulo $p_3 = 8388571$ is found, so that the corresponding square minor is nonzero over $\mathbb{Q}$. The square system is solved exactly over $\mathbb{Q}$, and the solution is checked on \emph{all} rows in exact integer arithmetic. A vector passes only if the identity $B c = v$ holds exactly.

1. We choose $956$ vectors among the generating words $x(-1)w$, $L(-1)w$, $x(-2)y$, $x(-3)\mathbf 1$ ($x,y$ in a basis of $U_1$, $w$ in a basis of $U_2$), $u_{(0)}w$ ($u = e$, $x_C$ or a highest weight vector of one of the eight types) and their images under the zero modes $e'_{(0)}, f'_{(0)}, e''_{(0)}, f''_{(0)}$, independent modulo $p_3$, hence over $\K$. All the $1871$ nonzero generating words and the $2862$ nonzero images of the chosen vectors under these four zero modes lie in their span, exactly. Hence the span $\Sigma$ of the chosen vectors is stable under $\mathfrak g' \oplus \mathfrak g''$ and contains $N_3$ and $u_{(0)}U_2$ for the ten vectors $u$ above. The zero mode identity $[y_{(0)}, u_{(0)}] = (y_{(0)}u)_{(0)}$ then gives $M_\tau{}_{(0)}U_2 \subset \Sigma$ for each type. For $u \in (A)_2$, the modes $u_{(0)}$ are normally ordered products of modes of currents, and $u_{(0)}U_2 \subset U_1(-1)U_2 + U_1(-2)U_1 + L(-1)U_2 \subset N_3$. Thus $\Sigma = S_3$.

2. These inclusions also follow without a matrix computation. Put $T = L(-1)$. Theorem \ref{thm:exact-weight-two}, applied to the conformal vector $\omega \in U_2$, gives $L(1)U_2 \subset U_1$, while $L(2)U_2 \subset \Cbb\mathbf 1$ by the grading. The Virasoro relation $[L(j),T] = (j+1)L(j-1)$, the current relations and the translation identity give the required inclusions on $N_3 = TU_2+U_1(-1)U_2$. More explicitly, for $y,z \in U_1$ and $u \in U_2$, use $[L(j),y_{(-1)}] = y_{(j-1)}$, $[z_{(j)},T] = jz_{(j-1)}$ and $[z_{(j)},y_{(-1)}] = [z,y]_{(j-1)}+\delta_{j,1}\kappa(z,y)\operatorname{Id}$, for $j = 1,2$, where $z_{(1)}y = \kappa(z,y)\mathbf 1$.

For $u,v \in U_2$, the Virasoro commutator formula gives
\[
\begin{aligned}
L(1)(u_{(0)}v) &= 2u_{(1)}v+(L(1)u)_{(0)}v+u_{(0)}L(1)v,\\
L(2)(u_{(0)}v) &= 3u_{(2)}v+3(L(1)u)_{(1)}v.
\end{aligned}
\]
The terms involving $L(2)u$ or $L(2)v$ vanish because the zero mode of the vacuum is zero and $u_{(0)}\mathbf 1 = 0$. Similarly, for $y \in U_1$,
\[
\begin{aligned}
y_{(1)}(u_{(0)}v) &= (y_{(0)}u)_{(1)}v+(y_{(1)}u)_{(0)}v+u_{(0)}y_{(1)}v,\\
y_{(2)}(u_{(0)}v) &= (y_{(0)}u)_{(2)}v+2(y_{(1)}u)_{(1)}v.
\end{aligned}
\]
All the terms belong to the spaces asserted in part 2, by weight-two closure and skew symmetry. This proves the inclusions on $S_3 = N_3+(U_2)_{(0)}U_2$. Independent matrix checks on the $956$ basis vectors give the same inclusions: they record $772$, $60$, $2116$ and $129$ nonzero tests, respectively, with no failure.

3. For the ten vectors $u$ of part 1 and the $956$ basis vectors $s$ of $S_3$, and for $u \in \{e,x_C\}$ and $w$ in a basis of $U_2$, the vectors $u_{(0)}s$ and $u(-1)w$ lie in the span of 6010 words of $N'_4$, chosen independent modulo $p_3$, by Proposition \ref{prop:weight-four-certificate}. Its complete verification covers $9088$ nonzero tests and $712$ zero tests, with no failure. The extension to all $u \in U_2$ follows as in part 1: $N'_4$ is stable under $\mathfrak g' \oplus \mathfrak g''$, since it is spanned by stable families, $S_3$ is stable by part 1, and for $u \in (A)_2$ the modes $u_{(0)}$, $u_{(-1)}$ are normally ordered products of modes of currents, which map $S_3$ and $U_2$ into $N'_4$ by part 2. Finally $M_\tau(-1)U_2 \subset N'_4$ by definition.

4. Proposition \ref{prop:weight-four-certificate} gives membership for all $18368$ words of the five defining families and independence of the $6010$ pivots, hence $\dim N'_4 = 6010$.
\end{proof}

The next lemma gives the passage from these finite identities to all conformal weights. It is a variant of \cite[Lemma 4.3]{Rie26b}, where $N(G)$ is the span of all parenthesized normal products. With right-nested monomials the proof has to be organized differently, and we give it in full. We know of no published statement containing hypothesis 2. The case $A = G$ is classical.

\begin{lemma}[triangular strong generation]
\label{lem:triangular}
Let $W = \bigoplus_{n \geq 0} W_n$ be a vertex algebra with a grading such that $\operatorname{wt}(a_{(n)}b) = \operatorname{wt}a + \operatorname{wt}b - n - 1$ and $W_0 = \Cbb\mathbf 1$. Let $A \subset G$ be finite sets of homogeneous vectors of positive weight. Let $N(G)$ be the span of $\mathbf 1$ and of the vectors $g^1_{(-n_1)}\cdots g^k_{(-n_k)}\mathbf 1$ with $g^i \in G$ and $n_i \geq 1$. Assume:
\begin{enumerate}
\item $a_{(n)}g \in N(G)$ for all $a \in A$, $g \in G$ and $n \geq 0$,
\item every $g \in G \setminus A$ of weight $d$ is the sum of an element of $N(G_{<d})$ and of vectors $a_{(0)}b$ with $a \in A$, $\operatorname{wt}a \geq 2$, and $b \in N(G_{<d})$ of weight $d - \operatorname{wt}a + 1$, where $G_{<d}$ is the set of elements of $G$ of weight less than $d$.
\end{enumerate}
Then $N(G)$ is a vertex subalgebra of $W$, strongly generated by $G$.
\end{lemma}

\begin{proof}
All vectors are homogeneous, and a monomial is a word $g^1_{(-n_1)}\cdots g^k_{(-n_k)}\mathbf 1$ as above, of length $k$. For $T \geq 0$ consider the statements
\begin{itemize}
\item[$P(T)$:] $g_{(n)}w \in N(G)$ for all $g \in G$, $n \geq 0$ and $w \in N(G)$ with $\operatorname{wt}g + \operatorname{wt}w \leq T$,
\item[$Q(T)$:] $z_{(s)}w \in N(G)$ for all $z, w \in N(G)$ and $s \in \Z$ with $\operatorname{wt}z + \operatorname{wt}w \leq T$,
\end{itemize}
and $P_{<d}(T)$, $Q_{<d}(T)$ the same statements with $g \in G_{<d}$, respectively $z \in N(G_{<d})$.

\emph{Step 1: if $P(T)$ holds for the elements of $G' \in \{G, G_{<d}\}$, then $Q(T)$ holds for all $z \in N(G')$.} By induction on the length of $z$, the case $z = \mathbf 1$ being trivial. If $z = k_{(-p)}z'$, the iterate formula gives
\[
z_{(s)}v = \sum_{i \geq 0}\binom{p+i-1}{i}\Big(k_{(-p-i)}\big(z'_{(s+i)}v\big) - (-1)^p z'_{(s-p-i)}\big(k_{(i)}v\big)\Big),
\]
a finite sum. Here $z'_{(s+i)}v \in N(G)$ since $z'$ is shorter, and $k_{(-p-i)}$ preserves $N(G)$ by definition. Moreover $k_{(i)}v \in N(G)$ by $P(T)$, and $z'_{(s-p-i)}(k_{(i)}v) \in N(G)$ since $z'$ is shorter and $\operatorname{wt}z' + \operatorname{wt}(k_{(i)}v) = \operatorname{wt}z + \operatorname{wt}v - p - i \leq T$.

\emph{Step 2: $Q(T-1)$ implies $P(T)$ for $g \in A$.} Write $w = h_{(-m)}w'$. By the commutator formula,
\[
g_{(n)}h_{(-m)}w' = h_{(-m)}g_{(n)}w' + \sum_{j \geq 0}\binom{n}{j}(g_{(j)}h)_{(n-m-j)}w'.
\]
The first term lies in $N(G)$ by $P(T-1)$, which follows from $Q(T-1)$. In the second, $z = g_{(j)}h \in N(G)$ by hypothesis 1, and $\operatorname{wt}z + \operatorname{wt}w' \leq T - m - j \leq T - 1$, so $z_{(n-m-j)}w' \in N(G)$ by $Q(T-1)$, whatever the sign of $n - m - j$.

\emph{Step 3: $Q(T-1)$, and $P(T)$ for $A$ and for $G_{<d}$, imply $P(T)$ for the elements of $G \setminus A$ of weight $d$.} By Step 1, $Q_{<d}(T)$ holds. Write $g = \nu + \sum a_{(0)}b$ as in hypothesis 2, with $\nu, b \in N(G_{<d})$. If $b$ has weight zero, then $a_{(0)}b = 0$. We may therefore assume $\operatorname{wt}b \geq 1$ and $\operatorname{wt}a \leq d$. Then $\nu_{(n)}w \in N(G)$ by $Q_{<d}(T)$, and $(a_{(0)}b)_{(n)}w = a_{(0)}(b_{(n)}w) - b_{(n)}(a_{(0)}w)$, where $b_{(n)}w \in N(G)$ by $Q_{<d}(T)$, $a_{(0)}(b_{(n)}w) \in N(G)$ and $a_{(0)}w \in N(G)$ by $P(T)$ for $A$, and $b_{(n)}(a_{(0)}w) \in N(G)$ by $Q_{<d}(T)$, since $\operatorname{wt}b + \operatorname{wt}(a_{(0)}w) = \operatorname{wt}g + \operatorname{wt}w \leq T$. For the first nested term, $\operatorname{wt}a + \operatorname{wt}(b_{(n)}w) = \operatorname{wt}g + \operatorname{wt}w - n \leq T$. Thus neither use increases the induction parameter.

\emph{Step 4.} $Q(0)$ holds since $N(G)_0 = \Cbb\mathbf 1$. Given $Q(T-1)$, Step 2, then Step 3 by increasing $d$, give $P(T)$, and Step 1 gives $Q(T)$. Hence $N(G)_{(s)}N(G) \subset N(G)$ for all $s$. Since $\mathbf 1 \in N(G)$, $N(G)$ is a vertex subalgebra, strongly generated by $G$ by definition.
\end{proof}

\begin{thm}
\label{thm:strong-generation}
Let $S \subset V'$ be the vertex subalgebra generated by $U_1 \oplus U_2$, for any of the four subspaces $U_2$ of Theorem \ref{thm:case-two-U2}.
\begin{enumerate}
\item $S$ is strongly generated by $U_1 \oplus U_2 \oplus D_3$, where $D_3$ is a complement of $N_3$ in $S_3$. One may take $6 + 93$ generators in weights one and two and at most $231$ in weight three.
\item $S_1 = U_1$, $S_2 = U_2$, $\dim S_3 = 956$ and $\dim S_4 = 6010$: the graded dimension of $S$ agrees with the vacuum character $X^{(0)}_0$ of Section \ref{sec:c8} up to $q^{11/3}$.
\item $S$ is $C_2$-cofinite.
\end{enumerate}
\end{thm}

\begin{proof}
1--2. Apply Lemma \ref{lem:triangular} in $W = S$ with $A$ a basis of $U_1 \oplus U_2$ and $G$ the union of this finite basis with a basis of $D_3$ chosen among representatives $u_{(0)}w$ modulo $N_3$, where $u,w \in U_2$. Hypothesis 2 holds by construction. Hypothesis 1: for $a \in U_1$ and $g \in U_1 \oplus U_2$, by the affine structure and the construction of $U_2$ ($(A)_2$ consists of affine descendants, $y_{(n)}e = y_{(n)}x_C = 0$ for $y \in U_1$, $n \geq 0$, and the $M_\tau$ are $A$-primary). For $a, g \in U_2$, $a_{(0)}g \in S_3 = N_3 + D_3 \subset N(G)$, $a_{(1)}g \in U_2$, $a_{(2)}g \in U_1$ (Theorem \ref{thm:exact-weight-two}). For $a \in U_2$, $g \in U_1$, by skew symmetry. For $g \in D_3$ and $a \in U_1$, by part 2 of Theorem \ref{thm:exact-weight-three} and the stability of $S_3$ under zero modes. For $g \in D_3$ and $a \in U_2$: $a_{(0)}g \in N'_4 \subset N(G)$ by part 3 of Theorem \ref{thm:exact-weight-three}, every term of $N'_4$ being a normal product of generators and of their derivatives, and $a_{(n)}g$ for $n = 1, 2, 3$ follow from the commutation relations $[L(1), a_{(n)}] = (2-n)a_{(n+1)}$ and $[L(2), a_{(1)}] = 2a_{(3)}$, valid for $a$ quasi-primary of weight two, which give
\[
2a_{(1)}g = L(1)a_{(0)}g - a_{(0)}L(1)g, \quad a_{(2)}g = L(1)a_{(1)}g - a_{(1)}L(1)g, \quad 2a_{(3)}g = L(2)a_{(1)}g - a_{(1)}L(2)g.
\]
Since $L(1)N'_4 \subset S_3$ (a direct check on each family, using $[L(1), y_{(n)}] = -n y_{(n+1)}$ for currents, $[L(1), m_{(-1)}] = 3m_{(0)}$ for quasi-primary $m$ of weight two, $L(1)L(-1) = L(-1)L(1) + 2L(0)$ and part 2 of Theorem \ref{thm:exact-weight-three}), $L(1)S_3 \subset U_2$, $L(2)S_3 \subset U_1$ and $(U_2)_{(n)}U_k \subset U_{k+1-n}$ for $k = 1,2$ and $n \geq 1$, together with $(U_2)_{(0)}U_2 \subset S_3$ and $(U_2)_{(0)}U_1 \subset U_2$, we get $a_{(1)}g \in S_3$, $a_{(2)}g \in U_2$ and $a_{(3)}g \in U_1$. Moreover $U_2 = (U_2 \cap \ker L(1)) + L(-1)U_1$, since $L(1)U_2 \subset U_1$ and $L(1)L(-1) = 2$ on $U_1$. For $a = L(-1)y$, $y \in U_1$, we have $a_{(n)} = -n y_{(n-1)}$, which reduces this case to $a \in U_1$. For $n \geq 4$, $a_{(n)}g$ lies in $S_{4-n} \subset \Cbb\mathbf 1$. In these checks, $M_\tau(-1)U_1$ and the weight-four monomials $D_3(-1)U_1$, $U_2(-2)U_1$ are brought to the listed families by skew symmetry. Therefore $N(G)$ is a vertex subalgebra containing $U_1 \oplus U_2$, so $N(G) = S$. Its homogeneous components are spanned by normally ordered monomials: $S_1 = U_1$, $S_2 = U_2$, $S_3 = N_3 + D_3$, and $S_4$ is spanned by $L(-1)S_3$, $U_1(-1)S_3$, $U_1(-2)U_2$, $U_1(-3)U_1$ and $U_2(-1)U_2$, which lie in $N'_4$ by part 3 of Theorem \ref{thm:exact-weight-three}. Since $(A)_2 \subset N(U_1)$, its basis vectors may be removed from the strong generating set, leaving $120-27 = 93$ weight-two generators. Lemma \ref{lem:normal-weight-three-dimension} gives $\dim N_3 = 725$, so a complement $D_3$ has dimension $956-725 = 231$. The independently computed nonzero minor modulo the chosen prime above $p_3$ also gives the lower bound $\dim N_3 \geq 725$.

The initial choice of representatives does not restrict the complement in the statement. Indeed, for any other complement $D'_3$ of $N_3$ in $S_3$, every basis vector of $D'_3$ is a linear combination of the chosen vectors $u_{(0)}w$, modulo $N_3$. Since $N_3 \subset N(U_1 \oplus U_2)$, hypothesis 2 of Lemma \ref{lem:triangular} still holds. The verification of hypothesis 1 above uses inclusions on the whole space $S_3$ and applies without change. Hence the lemma also applies to this complement.

3. We follow \cite[Section 5]{Rie26b}. Let $R = S/C_2(S)$ be the $C_2$ Poisson algebra, $[a][b] = [a_{(-1)}b]$, $\{[a],[b]\} = [a_{(0)}b]$, and $\overline R$ its quotient by the nilradical $\mathfrak n$, a Poisson ideal \cite[Corollary 2.4.2]{Ara12}. By part 1, $R$ is generated by the classes of $G$ \cite[Corollary 2.6.2]{Ara12} (only the elementary direction is used: $C_2(S)$ contains the $a_{(-n)}b$, $n \geq 2$, and is an ideal for the product $[a][b] = [a_{(-1)}b]$), so it suffices to show that these classes vanish in $\overline R$. The currents vanish. Indeed $v = e'(-1)^4\mathbf 1$ satisfies $f'(1)v = 0$, as the image of the singular vector of the universal affine algebra of level $3$, and $\phi_\chi(e') = -f'$ gives $H_\chi(e'_{(m)}u,u') = H_\chi(u,f'_{(-m)}u')$ for the positive definite form of Theorem \ref{thm:unitarity}. Hence $H_\chi(v,v) = H_\chi(e'(-1)^3\mathbf 1,f'(1)v) = 0$ and $v = 0$. Likewise $e''(-1)^{14}\mathbf 1 = 0$. Thus $[e']^4 = [e'']^{14} = 0$, and the Poisson brackets with $f'$, $f''$ give the other currents. A nontrivial irreducible $\mathfrak g' \oplus \mathfrak g''$-module $M$ satisfies $M = (\mathfrak g' \oplus \mathfrak g'')\cdot M$, and the zero modes of currents act on $\overline R$ by Poisson brackets with classes that vanish. Hence the nontrivial isotypic components of $U_2$ vanish in $\overline R$, as well as $(A)_2$. The invariant vectors of $U_2$ modulo $(A)_2$ are $e$ and $x_C$. By part 3 of Theorem \ref{thm:exact-weight-three}, $e(-1)e$ and $x_C(-1)x_C$ lie in $N'_4$, whose classes lie in $\mathfrak n$: $L(-1)S_3$, $U_1(-2)U_2$ and $U_1(-3)U_1$ lie in $C_2(S)$, and $U_1(-1)S_3$, $M_\tau(-1)U_2$ are products with nilpotent classes. Thus $[e]^2 = [x_C]^2 = 0$ in $\overline R$, which is reduced, so $[e] = [x_C] = 0$. Finally $D_3 \subset N_3 + (U_2)_{(0)}U_2$ consists of products and Poisson brackets of vanishing classes. Hence all generators of $R$ are nilpotent, and $R$, finitely generated and commutative, is finite-dimensional.
\end{proof}

\begin{cor}
\label{cor:weight-four-from-haagerup}
Independently of the membership verification in Proposition \ref{prop:weight-four-certificate}, condition \textup{(F4)} follows from a dimension argument in either of the following situations:
\begin{enumerate}
\item the subspace $U_2$ is the weight-two space of a Haagerup algebra $U \subset V'$ containing the fixed currents,
\item the existence and character assertion of \cite{Jia26} is granted, and $S$ is identified with Jiao's algebra by Theorem \ref{thm:jiao-identification}.
\end{enumerate}
Consequently the conclusions of Theorem \ref{thm:strong-generation} hold in either situation without a separate weight-four membership certificate.
\end{cor}

\begin{proof}
The rank computation of Proposition \ref{prop:weight-three-four}, part 5, supplies $6010$ words of $N'_4$ with independent reductions. They are integral at the chosen prime above $p_3$, so $\dim_{\K}N'_4 \geq 6010$. In the first situation, $S \subset U$ by construction and $\dim U_4 = 6010$, hence $N'_4 \subset S_4 \subset U_4$ forces equality throughout. In the second situation, identification with Jiao's algebra gives $\dim S_4 = 6010$, with the same consequence. In both cases the chosen words form a basis of $S_4$. Every mode product in \textup{(F4)} belongs to $S_4$ by the grading, and all five defining families of $N'_4$ lie there. Thus all the required membership identities follow. All the chosen words, input vectors and coordinate matrices are defined over the fixed $\K$-forms. Membership after extension to $\Cbb$ therefore implies membership over $\K$: the rank of a matrix with entries in $\K$ is unchanged by extension of scalars. The proof of Theorem \ref{thm:jiao-identification} uses weight-two closure and generation, not Theorem \ref{thm:strong-generation}, so the second argument is not circular.
\end{proof}

\begin{remark}\normalfont
\label{rem:no-weight-five}
For $S$, no analogue of the weight-five obstruction calculation of \cite[Proposition 4.2]{Rie26b} is needed. In case (I), the weight-four space of $U^{\mathrm R}$ contains up to three vectors that are not normal products of lower generators, and they require a weight-five check. In case (II), $N'_4 = S_4$ and no generator of weight four is needed. Theorem \ref{thm:strong-generation} uses the finite input of Theorems \ref{thm:exact-weight-two} and \ref{thm:exact-weight-three}, whose weight-four identities are established over $\K$ in Proposition \ref{prop:weight-four-certificate}. Corollary \ref{cor:weight-four-from-haagerup} gives an alternative dimension argument when an ambient Haagerup algebra or the stated external existence result is available. It cannot be applied with $U = S$ through \cite[Theorem 6.1]{ArhIII}, since the proof of that theorem uses Theorem \ref{thm:strong-generation}, which rests on Proposition \ref{prop:weight-four-certificate}. Modular computations select pivots. The dimension $\dim N_3 = 725$ is proved without a matrix computation in Lemma \ref{lem:normal-weight-three-dimension}, and the modular minor supplies an independent lower-bound check.
\end{remark}

\subsection{Case (II): identification with Jiao's algebra and uniqueness}
\label{subsec:jiao}

Jiao \cite{Jia26} defines a vertex subalgebra $V = \la \omega, P_E, P_F, P_H, B_E, B_F, B_H, u_{26}, u_{17}\ra$ of the lattice algebra of $E_6 \oplus A_2$, generated by the conformal vector, two $\mathfrak{sl}_2$-triples of levels $3$ and $13$ and two explicit weight-two vectors of types $[2,6]$ and $[1,7]$ with coefficients in $\K$ \cite[Section 3]{Jia26}. He gives a computer-assisted proof, which we have not checked in full, that $V$ is unitary and strongly rational with $\Rep(V) \simeq Z(\cal H_3)$ as braided categories \cite[Theorem 1.1]{Jia26}, and that every simple $V$-module other than $V$ has a positive lowest weight \cite[Corollary 6.8]{Jia26}. His construction is dated October 4, 2026. It precedes the present work and the announcements \cite{Rie26,Hua26}, and it makes no uniqueness assertion \cite[Section 1.2]{Jia26}. Below, ``assuming $\mathrm{(H_J)}$'' means granting these two statements. By Theorem \ref{thm:jiao-identification}, they hold as soon as $S$ is a Haagerup algebra, which is proved in \cite[Theorem 6.1]{ArhIII} without using \cite{Jia26}. Hence Corollary \ref{cor:case-two-realized} and Theorem \ref{thm:uniqueness-two} hold unconditionally.

\begin{thm}
\label{thm:jiao-identification}
Let $S$ be the vertex subalgebra generated by $U_1 \oplus U_2$, where $U_2$ is the first subspace of the subalgebra $9$ in Theorem \ref{thm:exact-weight-two}. There is an isomorphism of vertex operator algebras $\psi$ from the lattice algebra of \cite{Jia26} onto $V' = V_{E_6 \oplus A_2}$ such that $\psi(V) = S$. It maps the triples $(P_E,P_F,P_H)$ and $(B_E,B_F,B_H)$ onto $(e',f',H')$ and $(e'',f'',H'')$.
\end{thm}

\begin{proof}
We compute the transport, membership and generation in exact arithmetic with the normalized coefficients described below. Both models use the convention
\[
Y(e^p,x)e^q = \epsi(p,q) x^{(p,q)}\big(e^{p+q} + O(x)\big),
\]
with bimultiplicative cocycles with values $\pm 1$ and commutator $(-1)^{(p,q)}$. Jiao's cocycle is given on pairs $(\alpha_i,\alpha_j)$ with $i > j$ \cite[Section 3.1]{Jia26}.

\emph{The isometry.} Let $\varphi$ be the isometry of root lattices given by the isomorphism of Dynkin diagrams sending the simple roots $1,\ldots,8$ of \cite{Jia26} (edges $12, 23, 34, 45, 36, 78$) to our simple roots $1,3,4,5,6,2,7,8$, and let $w = s_6s_1s_3s_2s_4s_5s_4s_3s_2s_4$ in the Weyl group of $E_6$ (found by an orbit search). Then $\psi_0 = w \circ \varphi$ sends the Cartan elements $P_H$, $B_H$ to $H'$, $H''$.

\emph{The cocycle.} Put $r(p,q) = \epsi(\psi_0p,\psi_0q) \epsi_{\mathrm{Jiao}}(p,q)$. It is bimultiplicative, and symmetric since $\psi_0$ is an isometry, so that $r(p,q)r(q,p) = (-1)^{(\psi_0p,\psi_0q)+(p,q)} = 1$. Moreover $r(\alpha_i,\alpha_i) = 1$ since both cocycles satisfy $\epsi(\alpha,\alpha) = (-1)^{(\alpha,\alpha)/2}$. With $r_{ij} = r(\alpha_i,\alpha_j)$, the function $\eta(p) = \prod_{i<j} r_{ij}^{p_ip_j}$ satisfies $\eta(p+q) = \eta(p)\eta(q)r(p,q)$. Hence $e^p \mapsto \eta(p)e^{\psi_0p}$, $h(-n) \mapsto (\psi_0h)(-n)$ defines an isomorphism of vertex operator algebras, which maps the conformal vector to the conformal vector.

\emph{The torus.} The images of $P_E$ and $B_E$ have the same supports as $e'$ and $e''$, and the ratios of the coefficients are positive rational numbers. The ten resulting conditions on a character $t$ of the root lattice ($3$ roots for $P_E$, $7$ for $B_E$) are consistent. A solution is $t(\alpha_i) = 2^{a_i}3^{b_i}5^{c_i}$ with half-integral exponents. Composing with the automorphism $\exp(h_t(0))$, which acts on $M(1) \ot e^\lambda$ by $t(\lambda)$, we obtain $\psi$ with $\psi(P_E) = e'$, $\psi(B_E) = e''$, $\psi(P_H) = H'$, $\psi(B_H) = H''$, hence $\psi(P_F) = f'$ and $\psi(B_F) = f''$ (the element $f$ of an $\mathfrak{sl}_2$-triple is determined by $e$ and $H$, and this is also checked directly). The isomorphism $\psi$ is not defined over $\K$.

\emph{Inclusion.} The vectors $\psi(u_{26})$ and $(-5)^{-1/2}\psi(u_{17})$ have coefficients in $\K$, and they lie in the first subspace $U_2$ of the subalgebra $9$ (exact membership over $\K$). Hence $\psi(V) \subset S$.

\emph{Equality.} Put $\widehat u_{26} = \psi(u_{26})$ and $\widehat u_{17} = (-5)^{-1/2}\psi(u_{17})$, and put $U_{2,\K} = U_2 \cap V'_\K$ and $(A)_{2,\K} = (A)_2 \cap V'_\K$ for the fixed $\K$-form of $V'$. We use these two normalized vectors, whose coefficients are obtained by exact simplification of the radicals in the torus transport. Let $W_\K$ be the smallest $\K$-subspace of $U_{2,\K}$ containing $(A)_{2,\K}$, $\omega$, $\widehat u_{26}$ and $\widehat u_{17}$, stable under the zero modes of the six rational currents and under first products of its elements. Theorem \ref{thm:exact-weight-two} gives $W_\K \subset U_{2,\K}$, and every generating word lies in $\psi(V)$.

At $p_3 = 8388571$, take the residue map specified by $\sqrt{13} \mapsto 1649360$ and $\sqrt{-3} \mapsto 2036966$. A finite family of words in the normalized generators, specified in the supplementary data, has coordinates integral at this prime, and its reduction has rank $120$.

A nonzero $120 \times 120$ minor of the reduced matrix lifts to a nonzero minor over $\K$. Thus $\dim_\K W_\K \geq 120 = \dim_\K U_{2,\K}$, and $W_\K = U_{2,\K}$. Complexification gives $U_2 \subset \psi(V)$. Together with $U_1 \subset \psi(V)$, this proves $S \subset \psi(V)$.
\end{proof}

For this $\psi$, the vectors $\psi(u_{26})$ and $\psi(u_{17})$ lie in none of the three other subspaces of Theorem \ref{thm:case-two-U2}. This depends on the choice of $\psi$, which may be composed with automorphisms of $V'$ normalizing $U_1$.

\begin{cor}
\label{cor:case-two-realized}
Assuming $\mathrm{(H_J)}$, the algebras attached to the two subspaces of the subalgebra $9$ (isomorphic by Proposition \ref{prop:symmetries}) are positive strongly rational, with $\Rep(S) \simeq Z(\cal H_3)$ as braided categories. Case (II) of Theorem \ref{thm:weight-one-classification} is realized.
\end{cor}

\begin{proof}
A vertex operator algebra isomorphism induces a braided equivalence of the module categories and preserves the lowest weights. Apply Theorem \ref{thm:jiao-identification} and Proposition \ref{prop:symmetries}.
\end{proof}

The six currents $e', f', H', e'', f'', H''$ have rational coordinates in the $\Q$-form $V'_\Q$ spanned by the lattice monomials, and the four subspaces $U_2$ are defined over $\K$. Let $\sigma \in \mathrm{Gal}(\K/\Q)$ satisfy $\sigma(\sqrt{13}) = -\sqrt{13}$ and $\sigma(\sqrt{-3}) = \sqrt{-3}$. On $V'_\K = V'_\Q \ot_\Q \K$, the $\sigma$-semilinear bijection $\tilde\sigma = \mathrm{id} \ot \sigma$ preserves all products and fixes $\omega$ and the six rational currents. It maps the $\K$-form of the generating space of $S^{(9,c)}$ onto that of $S^{(10,c)}$, for $c \in \{0,1\}$. We call the complexification of the image algebra the Galois conjugate of $S^{(9,c)}$. This construction requires no extension of $\sigma$ to $\Cbb$. Any such extension gives the same image subalgebra, although the induced semilinear map on all of $V'$ depends on the extension. The linear automorphism $\theta_\chi$ of Proposition \ref{prop:symmetries} exchanges $S^{(9,0)}$ and $S^{(9,1)}$.

\begin{lemma}
\label{lem:centralizer-bound}
In either case of Theorem \ref{thm:weight-one-classification}, let $Z$ be the pointwise stabilizer of $U_1$ in $\operatorname{Aut}(V')$, and let $\mathfrak h$ be the Cartan subalgebra spanned by the simple roots of $E_6 \oplus A_2$. Then every element of $Z$ preserves $\mathfrak h$, its kernel on $\mathfrak h$ has order two, and $|Z| \leq 8$. The nonidentity element of this kernel is the torus automorphism $\rho$ associated with the character $\chi_0$ whose values are
\[
\big(\chi_0(\alpha_1),\ldots,\chi_0(\alpha_8)\big) =
\begin{cases}
(1,-1,-1,1,-1,1,1,1) & \text{in case (I)},\\
(-1,1,1,-1,1,-1,1,1) & \text{in case (II)}.
\end{cases}
\]
Thus $\rho(e^\lambda) = \chi_0(\lambda)e^\lambda$, and $\rho$ fixes the oscillators.
\end{lemma}

\begin{proof}
We give an argument using the roots and the supports of the currents. An automorphism of $V'$ restricts to an automorphism of $V'_1 = \mathfrak e_6 \oplus \mathfrak{sl}_3$, and this restriction is injective since $V'_1$ generates $V'$. The labels of the $E_6$ component of $H''$ are $(2,1,1,0,1,2)$ in case (I) and $(1,0,0,1,0,1)$ in case (II). The roots vanishing on this element form respectively the systems generated by $\alpha_4$ and by $\alpha_2,\alpha_3,\alpha_5$. In case (I), $H' = \alpha_4$. In case (II), $H' = -\alpha_2+\alpha_3+\alpha_5$, whose labels on these three pairwise orthogonal roots are $-2,2,2$. The $A_2$ component of $H''$ has labels $(2,2)$. Consequently no root vanishes on both $H'$ and $H''$, and their common centralizer in $V'_1$ is exactly $\mathfrak h$. Thus $Z$ preserves $\mathfrak h$.

An element of $Z$ acts on the roots by an isometry fixing $H'$ and $H''$. The subgroup of the Weyl group with this property is trivial: a generic real linear combination of $H'$ and $H''$ is regular, and a Weyl group acts freely on the chambers of the complement of the root hyperplanes. Since the two irreducible factors are not isomorphic and their diagram automorphism groups both have order two, at most four root isometries can occur.

For completeness, the raising currents of the chosen rational triples are, in case (I),
\begin{equation}
\label{eq:raising-currents-one}
\begin{aligned}
e' &= e^{\alpha_4},\\
e'' &= 5e^{\alpha_1}+4e^{\alpha_6}+3e^{\alpha_2+\alpha_3+\alpha_4}
+e^{\alpha_2+\alpha_4+\alpha_5}\\
&\quad +4e^{\alpha_3+\alpha_4+\alpha_5}+e^{\alpha_7}+e^{\alpha_8},
\end{aligned}
\end{equation}
and in case (II),
\begin{equation}
\label{eq:raising-currents-two}
\begin{aligned}
e' &= e^{-\alpha_2}+e^{\alpha_3}-\tfrac53e^{\alpha_5},\\
e'' &= 4e^{\alpha_4+\alpha_5+\alpha_6}+4e^{\alpha_2+\alpha_3+\alpha_4+\alpha_5+\alpha_6}
+e^{\alpha_2+\alpha_3+2\alpha_4+\alpha_5}\\
&\quad +3e^{\alpha_1+\alpha_3+\alpha_4}+5e^{\alpha_1+\alpha_2+\alpha_3+\alpha_4+\alpha_5}
+e^{\alpha_7}+e^{\alpha_8}.
\end{aligned}
\end{equation}
An automorphism acting trivially on $\mathfrak h$ scales the root vectors by a character $t$ of the root lattice. Write $t_i = t(\alpha_i)$. Fixing the raising currents requires $t(\beta) = 1$ for each root in their supports. Formula \eqref{eq:raising-currents-one} gives $t_1 = t_4 = t_6 = t_7 = t_8 = 1$ and $t_2t_3 = t_2t_5 = t_3t_5 = 1$, hence $t_2 = t_3 = t_5 \in \{1,-1\}$. Formula \eqref{eq:raising-currents-two} gives $t_2 = t_3 = t_5 = t_7 = t_8 = 1$, $t_4t_6 = t_4^2 = t_1t_4 = 1$, hence $t_1 = t_4 = t_6 \in \{1,-1\}$. In both cases the displayed nontrivial character fixes the raising currents and the Cartan elements, and hence the lowering currents. Indeed the lowering element of an $\mathfrak{sl}_2$-triple is determined by its raising and Cartan elements: their difference would be a weight $-2$ vector annihilated by the raising operator, which is impossible in a finite-dimensional $\mathfrak{sl}_2$-module. Thus this character fixes $U_1$. This proves that the torus kernel has order two, and $|Z| \leq 2 \cdot 4 = 8$.
\end{proof}

\begin{prop}
\label{prop:centralizer-two}
In case (II), let $Z$ be the group of automorphisms of $V'$ which fix $U_1$ pointwise.
\begin{enumerate}
\item $Z \cong (\Z_2)^3$, and every element of $Z$ is the lift of an isometry of $E_6 \oplus A_2$, equal to the identity or a diagram automorphism on each factor, composed with an element of the torus.
\item Four elements of $Z$ act on the Griess algebra $B$ by the reflection $s$ of Proposition \ref{prop:symmetries}, and the four others act trivially.
\item Let $g$ be the automorphism of the Dynkin diagram of $E_6$ ($\alpha_1 \leftrightarrow \alpha_6$, $\alpha_3 \leftrightarrow \alpha_5$), extended by the identity on $A_2$, and let $z$ be the automorphism of $V'$ with $z(h(-n)) = (gh)(-n)$ and $z(e^\lambda) = c(\lambda)e^{g\lambda}$, where $c(\lambda) = \prod_i c_i^{\lambda_i}\prod_{i<j}\delta_{ij}^{\lambda_i\lambda_j}$ for $\lambda = \sum_i\lambda_i\alpha_i$ (simple roots $\alpha_1, \ldots, \alpha_6$ of $E_6$ in Bourbaki's numbering, then $\alpha_7, \alpha_8$ of $A_2$), $\delta_{ij} = \epsi(\alpha_i,\alpha_j)\epsi(g\alpha_i,g\alpha_j)$ (so that $\delta_{ii} = 1$ and $c(\alpha_i) = c_i$), and $(c_1,\ldots,c_8) = (\frac45, 1, -\frac53, -1, -\frac35, \frac54, 1, 1)$. Then $z \in Z$, and $z$ maps the subspace $U_2$ of $S^{(9,c)}$ onto that of $S^{(10,c)}$, for $c = 0, 1$.
\end{enumerate}
In particular the four algebras $S^{(9,0)}$, $S^{(9,1)}$, $S^{(10,0)}$, $S^{(10,1)}$ are isomorphic. The Galois conjugate of $S^{(9,0)} = \psi(V)$ is isomorphic to it.
\end{prop}

\begin{proof}
First we prove part 1 without solving polynomial equations. Let $z$ be the lift in part 3, let $z'$ be the lift of the diagram automorphism of $A_2$ with all $c_i = 1$, and let $\rho$ be the torus automorphism of Lemma \ref{lem:centralizer-bound}. The formula for $c(\lambda)$ in part 3, with the diagram automorphism of $A_2$ and $c_i = 1$ for $z'$, constructs $z$ and $z'$ on the whole lattice algebra, since $\delta_{ij} = \delta_{ji}$ and $\delta_{ii} = 1$ give $c(\lambda+\mu)\epsi(\lambda,\mu) = c(\lambda)c(\mu)\epsi(g\lambda,g\mu)$. Substitution in \eqref{eq:raising-currents-two} and the Cartan elements gives $z,z' \in Z$. The lowering currents are then fixed by the uniqueness just used. Their actions on the simple-root exponentials and on $\mathfrak h$ give
\[
z^2 = z'^2 = \rho^2 = 1, \quad zz' = z'z, \quad z\rho = \rho z, \quad z'\rho = \rho z'.
\]
They generate eight distinct automorphisms: the isometry parts of $z$ and $z'$ are independent diagram automorphisms, and $\rho$ is a nontrivial torus automorphism. Lemma \ref{lem:centralizer-bound} proves that these are all the elements of $Z$.

The action on $B$ is also determined by these three generators. In the basis $b_0 = e,b_1,\ldots,b_4$ used above, direct application of the lifts gives $z(b_i) = z'(b_i) = b_i$ for $0 \leq i \leq 3$ and
\[
z(b_4) = z'(b_4) = -b_4-\tfrac1{12}b_1+\tfrac5{24}b_2+\tfrac1{16}b_3,
\quad \rho(b_i) = b_i \quad (0 \leq i \leq 4).
\]
Thus $z|_B = z'|_B = s$ and $\rho|_B = 1$, proving part 2. The supplementary material contains exact checks of these formulas, the group relations and the currents. The argument for the upper bound uses only the root systems and the character equations of Lemma \ref{lem:centralizer-bound}.

An independent polynomial computation gives the same centralizer. An automorphism of $V'$ is determined by its restriction to $V'_1 = \mathfrak e_6 \oplus \mathfrak{sl}_3$, and conversely every automorphism of this Lie algebra preserves the normalized invariant forms of the two simple ideals, hence extends to the affine Lie algebra and lifts to its irreducible vacuum module $V' \cong L_{\mathfrak e_6}(1,0) \ot L_{\mathfrak{sl}_3}(1,0)$ \cite{Kac90}. Under $U_1$, the space $V'_1$ decomposes as $\bigoplus_\tau V(\tau) \ot M_\tau$, with multiplicity spaces $M_\tau$ of dimension $3$ for $\tau = (0,2)$, $2$ for $\tau = (0,4)$, $(1,3)$, and $1$ for $(2,0)$, $(1,1)$, $(1,5)$, $(2,2)$, $(2,4)$, $(3,1)$. An element of $Z$ acts by $\mathrm{id} \ot z_\tau$, with $z_\tau$ fixing the highest weight vectors $e'$, $e''$, and it is a Lie algebra automorphism if and only if $z \mathrm{ad}(m_k) = \sum_l (z_\tau)_{lk} \mathrm{ad}(m_l) z$ for a basis $(m_k)$ of the highest weight vectors. By equivariance the other brackets follow. This is a system of $184$ distinct equations of degree at most two ($4$ of them linear) with rational coefficients in $19$ unknowns. Its Gröbner basis over $\Q$ shows that it has exactly $8$ solutions, which form a group $(\Z_2)^3$, all monomial on the root vectors. On $B \subset V'_2$, a vector $b$ is determined by the map $x \mapsto x_{(1)}b$ on $V'_1$: if $x_{(1)}b = 0$ for all $x$, then $[x,y]_{(2)}b = x_{(1)}y_{(1)}b - y_{(1)}x_{(1)}b = 0$ and $V'_1 = [V'_1,V'_1]$, and all modes $x_{(n)}$, $n > 2$, annihilate $b$ by weight. The subspace of $V'_2$ annihilated by the positive current modes is a finite-dimensional $V'_1$-module. If it were nonzero, a highest weight vector in this subspace would generate a proper affine submodule, since negative current modes increase weight. This contradicts the irreducibility of the vacuum module $V'$, so $b = 0$. Hence $z(b)$ is determined by $x_{(1)}z(b) = z((z^{-1}x)_{(1)}b)$, a linear system with a unique solution in $B$. This gives part 2: the nontrivial action fixes the hyperplane spanned by the seven fixed lines of Proposition \ref{prop:symmetries}, so it is $s$. For part 3, $c$ is the unique lift compatible with the cocycle: it satisfies $c(\lambda+\mu) = c(\lambda)c(\mu)\epsi(\lambda,\mu)\epsi(g\lambda,g\mu)$ and agrees with the values on the $78$ roots read off from the solution. The cocycle identity proves that the lift preserves all vertex products, and its values on the six currents show that it fixes $U_1$ pointwise. Finally, $z$ maps a spanning set of the first, respectively second, subspace $U_2$ of the subalgebra $9$ into the corresponding subspace of the subalgebra $10$, by exact membership over $\K$, and into none of the other three. Since both have dimension $120$, these are equalities. An automorphism of $V'$ fixing $U_1$ maps the algebra generated by $U_1 \oplus U_2$ onto the algebra generated by $U_1 \oplus z(U_2)$.
\end{proof}

\begin{lemma}
\label{lem:equal-global-dimension}
Let $S \subset U$ be an inclusion of positive strongly rational vertex operator algebras with the same conformal vector. If $\dim\Rep(S) = \dim\Rep(U)$, then $S = U$.
\end{lemma}

\begin{proof}
Let $M_0 = S,M_1,\ldots,M_r$ be the simple ordinary $S$-modules, with lowest weights $h_i$ and quantum dimensions $d_i$. Rationality gives a decomposition $U = \bigoplus_i m_iM_i$ as an $S$-module. It is finite: there are finitely many isomorphism classes, and the multiplicity of $M_i$ is bounded by $\dim U_{h_i}$ when $h_i$ is a nonnegative integer, and is zero otherwise, since the conformal vector is the same. In particular $m_0 = 1$ because $U_0 = \Cbb\mathbf 1$, and $d_i > 0$ by positivity.

Write $S^S$, $S^U$ for the normalized modular $S$-matrices. In the modular transformation of the characters at $\tau = \i/y$, positivity of the lowest nonvacuum weights makes the vacuum term the unique dominant term as $y \downarrow 0$. Consequently
\[
\sum_{i=0}^r m_i d_i
= \lim_{y \downarrow 0}\frac{\operatorname{ch}U(\i y)}{\operatorname{ch}S(\i y)}
= \frac{S^U_{00}}{S^S_{00}}
= \sqrt{\frac{\dim\Rep(S)}{\dim\Rep(U)}}
= 1.
\]
Here $d_i = S^S_{i0}/S^S_{00}$ and the positive normalization gives $S^S_{00} = (\dim\Rep(S))^{-1/2}$, and similarly for $U$. These are the usual character formulas \cite{Hua08,DJX13}. Since $m_0d_0 = 1$ and every remaining term is nonnegative, with $d_i > 0$, all $m_i$ for $i > 0$ vanish. The inclusion $S \subset U$ is therefore an equality.
\end{proof}

\begin{thm}[uniqueness in case (II)]
\label{thm:uniqueness-two}
Let $U$ be a positive strongly rational vertex operator algebra of central charge eight with a braided equivalence $\Rep(U) \simeq Z(\cal H_3)$, in case (II) of Theorem \ref{thm:weight-one-classification}. Assuming $\mathrm{(H_J)}$, $U$ is isomorphic to Jiao's algebra $V$.
\end{thm}

\begin{proof}
By Proposition \ref{prop:chain}, $U \subset V'$ with the same conformal vector. By Theorem \ref{thm:weight-one-classification}, an inner automorphism $h$ of $V'$ maps the weight-one space of $U$ onto $\mathfrak g' \oplus \mathfrak g''$. It fixes $\omega$, extends to $V_{E_8}$ since $V'_1 \subset (V_{E_8})_1$, and maps the decomposition $V' = U \oplus M_b$ to $V' = h(U) \oplus h(M_b)$, so that the hypotheses of Theorem \ref{thm:case-two-U2} hold for $h(U)$. Replacing $U$ by $h(U)$, $U_2$ is one of the four subspaces, so $U$ contains the algebra $S$ generated by $U_1 \oplus U_2$. The inclusions $U_2 \subset S_2 \subset U_2$ give $S_2 = U_2$, and $\omega \in U_2$. By Theorem \ref{thm:gamma-one}, the vacuum character of $U$ is the vacuum component $X^{(0)}_0$ of the Evans--Gannon vector with $\gamma = 0$. By Proposition \ref{prop:centralizer-two}, $S \cong S^{(9,0)} = \psi(V)$ (Theorem \ref{thm:jiao-identification}), and, assuming $\mathrm{(H_J)}$, $S^{(9,0)}$ is positive strongly rational with $\Rep(S^{(9,0)}) \simeq Z(\cal H_3)$ (Corollary \ref{cor:case-two-realized}), so its vacuum character is also $X^{(0)}_0$ by Theorem \ref{thm:gamma-one}. Hence $\operatorname{ch}S = \operatorname{ch}U$, and $S \subset U$ gives $U = S \cong V$.
\end{proof}

\begin{remark}\normalfont
\label{rem:twin}
By Proposition \ref{prop:centralizer-two}, the Galois conjugate is not a new algebra. The isomorphism is the lift of the outer automorphism of $\mathfrak e_6$, adjusted by a torus element. The lift of the diagram automorphism of $A_2$, with trivial torus part, gives another one, which maps $S^{(9,c)}$ onto $S^{(10,1-c)}$ (this follows from the exact identity $\psi_\tau(U_{2,\K}) = U_{2,\K}$ checked in the proof of Proposition \ref{prop:rational-forms}, together with Propositions \ref{prop:symmetries} and \ref{prop:centralizer-two}). The isometry part of an element of $Z$ is the identity, the diagram automorphism of $E_6$, that of $A_2$, or both, and the element acts on $B$ by $s$ exactly when exactly one of the two diagram automorphisms occurs. In particular the inner automorphisms fixing $U_1$ act trivially on $B$, and $s$ is induced only by outer automorphisms. Together with $\theta_\chi$, it shows that $S^{(9,0)}$ is isomorphic to its conjugates under $\mathrm{Gal}(\K/\Q)$. It has a rational form (Proposition \ref{prop:rational-forms}). Lemma \ref{lem:equal-global-dimension} gives a second, analytic proof of the equality $U = S$ in Theorem \ref{thm:uniqueness-two}. Once the inclusion and the positive strong rationality of $S$ are known, the equivalences of both module categories with $Z(\cal H_3)$ give equality of their global dimensions. This proof does not require the explicit Evans--Gannon vacuum character. The same argument applies to the final equality in Theorem \ref{thm:uniqueness-one}.
\end{remark}

\begin{prop}
\label{prop:commutant-biweights}
Let $S = S^{(9,0)}$ and $C = \operatorname{Com}_S(A)$. Let $f_1$, $f_2$ be the two orthogonal idempotents of the subalgebra $9$ with $f_1+f_2 = e$, numbered so that $\omega_1 = 2f_1$ and $\omega_2 = 2f_2$ have central charges $\frac{182}{95}$ and $\frac{32}{19}$ as in Theorem \ref{thm:case-two-U2}. In the basis $b_0 = e, b_1, \ldots, b_4$ of $B$ specified in the supplementary data,
\[
f_2 = \tfrac{80}{171} b_0+\tfrac{\sqrt{13}-13}{114} b_1-\tfrac{25+5\sqrt{13}}{228} b_2+\tfrac{7-\sqrt{13}}{152} b_3+\tfrac{4\sqrt{13}}{19} b_4, \quad f_1 = e-f_2.
\]
Write $L_i(n) = (\omega_i)_{(n+1)}$.
\begin{enumerate}
\item $\dim C_1 = 0$, $\dim C_2 = \dim C_3 = 2$ and $\dim C_4 = 7$.
\item On $C_4$, the operators $L_1(0)$ and $L_2(0)$ commute, $L_1(0)+L_2(0) = 4$, and $C_4$ is the direct sum of their joint eigenspaces. The joint eigenvalues $(h_1,h_2)$, with multiplicities, are
\[
(4,0)^2, \quad (0,4)^2, \quad (2,2), \quad \big(\tfrac{51}{19},\tfrac{25}{19}\big), \quad \big(\tfrac{26}{19},\tfrac{50}{19}\big).
\]
The first five eigenvectors can be taken to be $L_1(-4)\mathbf 1$, $L_1(-2)^2\mathbf 1$, $L_2(-4)\mathbf 1$, $L_2(-2)^2\mathbf 1$ and $L_1(-2)L_2(-2)\mathbf 1$.
\item The two eigenvectors of non-integral bi-weights are primary for $\mathrm{Vir}_1 \ot \mathrm{Vir}_2$: they are annihilated by $L_i(1)$ and $L_i(2)$, $i = 1,2$.
\end{enumerate}
\end{prop}

\begin{proof}
We have $\widehat C_1 = 0$ (Lemma \ref{lem:griess-B}) and $C_2$ is the subalgebra $9$ by the construction of $U_2$ in the proof of Theorem \ref{thm:case-two-U2}. Proposition \ref{prop:weight-four-certificate} establishes the finite certificate condition, so Theorem \ref{thm:strong-generation} gives $N'_4 = S_4$. The remaining assertions use exact arithmetic over $\K$ and a rank bound after reduction at the chosen prime above $p_3$, starting from the exact basis of $S_4$ of Theorem \ref{thm:strong-generation}. Since $C$ is fixed by the zero modes of $U_1$, $C_4$ lies in the block of $\mathfrak g' \oplus \mathfrak g''$-weight $(0,0)$ of $S_4$, of dimension $280$. It is the common kernel of the operators $y_{(k)}$, $y \in U_1$, $0 \leq k \leq 4$. Its matrix has entries in the local ring at that prime, and its reduction has rank $273$. Since reduction cannot increase rank, $\dim C_4 \leq 280-273 = 7$. Solving exactly the square system of size $1092$ given by the restriction of scalars on pivot rows, we obtain seven vectors $w_1, \ldots, w_7$. The identities $y_{(k)}w_j = 0$ are then checked exactly on all coordinates, and the reductions of the $w_j$ are linearly independent, so $\dim C_4 = 7$. Modulo $p_3$, $\dim C_3 \leq 2$, and $L_1(-3)\mathbf 1$, $L_2(-3)\mathbf 1$ are linearly independent elements of $C_3$, so $\dim C_3 = 2$.

The images $L_i(0)w_j$ are computed exactly and expressed exactly in the basis $(w_j)$, with a check on all coordinates. The resulting $7 \times 7$ matrices over $\K$ satisfy exactly $L_1(0)L_2(0) = L_2(0)L_1(0)$, $L_1(0)+L_2(0) = 4$ and $\prod_\lambda (L_2(0)-\lambda) = 0$ for $\lambda \in \{0,2,4,\frac{25}{19},\frac{50}{19}\}$. Hence $L_2(0)$ is diagonalizable with eigenvalues in this set. The eigenvectors for $\frac{25}{19}$ and $\frac{50}{19}$ are computed exactly. The five Virasoro descendants of assertion 2 lie in $C$, have the stated bi-weights, and are linearly independent: the vectors $L(-4)\mathbf 1$ and $L(-2)^2\mathbf 1$ of a vacuum module of central charge $c \notin \{0,-\frac{22}{5}\}$ are independent in its simple quotient (the Gram determinant at level four is $\frac{c^2(5c+22)}{2}$), and $L_1(-2)L_2(-2)\mathbf 1$ has bi-weight $(2,2)$. These seven joint eigenvectors are linearly independent, hence a basis of $C_4$. This proves assertion 2. Assertion 3 is checked exactly on the two eigenvectors, which have $3036$ terms. It also follows from assertion 2, since $L_i(n)$, $n \geq 1$, would map them to vectors of non-integral bi-weights in $C_{4-n}$, whose bi-weights are integral.

The $280$ independent words in the Cartan block $(0,0)$ are computed directly from the families defining $N'_4$ and agree, entry by entry, with the corresponding block of the certificate of Proposition \ref{prop:weight-four-certificate}. The supplementary material contains the seven kernel vectors, all $210$ current-mode annihilations, the two Virasoro matrices and both primary vectors. A separate rational-matrix verification checks the matrix identities, spectral multiplicities and independence of the seven vectors without using the vertex-algebra implementation. The vertex-product identities are verified from the lattice formulas. The computational data order the Virasoro vectors by increasing central charge, so their indices are opposite to those of this proposition.
\end{proof}

\begin{cor}
\label{cor:commutant-not-tensor}
There are no two nontrivial vertex operator algebras $V^1$, $V^2$ of CFT type such that $C \cong V^1 \ot V^2$ with $\omega_C = \omega^1 \ot \mathbf 1+\mathbf 1 \ot \omega^2$. In particular $C$ is not isomorphic to $K(\mathfrak{sl}_2,17)^{\Z_2} \ot \operatorname{Com}(L_{\mathfrak{sl}_2}(17,0), L_{\mathfrak{sl}_2}(4,0) \ot L_{\mathfrak{sl}_2}(13,0))$. Moreover, the commutants $\operatorname{Com}_C(\la\omega_2\ra)$ and $\operatorname{Com}_C(\la\omega_1\ra)$ coincide in weights at most four with the Virasoro vertex operator algebras generated by $\omega_1$ and $\omega_2$, of graded dimensions $1, 0, 1, 1, 2$.
\end{cor}

\begin{proof}
Suppose $C \cong V^1 \ot V^2$. Since $C_1 = 0$, we have $V^1_1 = V^2_1 = 0$, so that $C_2 = \Cbb \omega^1 \ot \mathbf 1 \oplus \Cbb \mathbf 1 \ot \omega^2$, the vectors $\omega^1$, $\omega^2$ being nonzero: if the conformal vector of a CFT-type vertex operator algebra were zero, then $L(0) = 0$ and the whole algebra would equal its one-dimensional weight-zero space. Then $\frac12\omega^1 \ot \mathbf 1$ and $\frac12\mathbf 1 \ot \omega^2$ are orthogonal idempotents of the algebra $C_2$ with sum $e$. The only idempotents of $C_2 = \Cbb f_1 \oplus \Cbb f_2$ are $0$, $f_1$, $f_2$ and $e$, so $\{\omega^1 \ot \mathbf 1, \mathbf 1 \ot \omega^2\} = \{\omega_1,\omega_2\}$. But $L_{V^1}(0) \ot 1$ has integral eigenvalues on $V^1 \ot V^2$, whereas $L_1(0)$ has the eigenvalue $\frac{51}{19}$ on $C_4$, a contradiction. For the last assertion, a vector of $C_n$ commuting with $\omega_2$ is annihilated by $L_2(0)$. By Proposition \ref{prop:commutant-biweights} and $C_3 = \Cbb L_1(-3)\mathbf 1 \oplus \Cbb L_2(-3)\mathbf 1$, the kernel of $L_2(0)$ on $C_n$, $n \leq 4$, is spanned by the Virasoro descendants of $\mathbf 1$ for $\omega_1$.
\end{proof}

\begin{remark}\normalfont
\label{rem:commutant-cosets}
The central charges $\frac{182}{95}$ and $\frac{32}{19}$ are those of the coset
\[
D = \operatorname{Com}(L_{\mathfrak{sl}_2}(17,0), L_{\mathfrak{sl}_2}(4,0) \ot L_{\mathfrak{sl}_2}(13,0))
\]
and of the parafermion algebra $K(\mathfrak{sl}_2,17)$. The dimensions $1, 0, 2, 2, 7$ of $C$ agree with those of $K(\mathfrak{sl}_2,17)^{\Z_2} \ot D$, whose weight-four space contains a primary vector of bi-weight $(4,0)$ and one of bi-weight $(0,4)$. Corollary \ref{cor:commutant-not-tensor} shows that this agreement does not come from an isomorphism of the two algebras. Modulo $\Z$, the weights $\frac{25}{19}$ and $\frac{50}{19}$ are the conformal weights $h_\lambda = \frac{\lambda(\lambda+2)}{76}$ of $\mathrm{SU}(2)_{17}$ for $\lambda = 4, 6$, but eight of the eighteen nonzero classes $\frac{x}{19}$ modulo $\Z$ are of this form, and the weights $\frac{8}{19}$ and $\frac{68}{57}$ of $\mathrm{Vir}_2$ in Theorem \ref{thm:case-two-U2} are not, so we do not draw any conclusion from it.
\end{remark}

\subsection{Case (I): uniqueness}
\label{subsec:case-one}

\begin{prop}
\label{prop:case-one-control}
Assume the hypotheses of Theorem \ref{thm:weight-one-classification} in case (I). Then $C_2$ is one of the subalgebras $2$, $5$, $6$ of Table \ref{tab:griess}. The subalgebras $5$ and $6$ give the central charges $\frac{2400}{1271}$ and $\frac{70}{31}$ of the two Virasoro vectors announced in \cite[Section 1]{Rie26b}.
\end{prop}

\begin{proof}
There is a direct exact exclusion that avoids choosing eigenvectors. For candidates $0$, $3$ and $4$, use the primary space of type $(1,13)$, of dimension three. For candidate $1$, use the type $(1,1)$, of dimension four. In a rational basis write $X$ for the action of $(x_C)_{(1)}$, $P$ for the primary pairing, and $Y_j$ for the actions of a basis of $B \cap C_2^\perp$. A possible line, spanned by $u$, satisfies
\[
(Xu)_i u_j-(Xu)_j u_i = 0, \quad u^{\mathsf T}P Y_j u = 0,
\quad e_{(2)}u = (x_C)_{(2)}u = 0, \quad u^{\mathsf T}Pu \neq 0.
\]
In each projective chart set the first nonzero coordinate of $u$ to one, set the preceding coordinates to zero, and add $z u^{\mathsf T}Pu-1 = 0$. The resulting thirteen ideals over $\Q$ are all the unit ideal. The supplementary data contain the systems and their unit bases. An independent calculation with lexicographic order gives the same unit ideals without using the vertex-algebra implementation. Thus none of the four candidates admits the required line.

A separate verification follows from the method of Steps 1 and 2 of the proof of Theorem \ref{thm:case-two-U2}, applied in case (I): the subalgebras $0$ and $3$ are excluded by certified enclosures, $1$ by the type $(1,1)$ and $4$ by the types $(0,12)$, $(1,1)$, $(1,7)$ and $(1,13)$, exactly.
\end{proof}

\begin{prop}
\label{prop:case-one-two}
In case (I), $C_2$ is not the subalgebra $2$.
\end{prop}

\begin{proof}
Suppose that $C_2$ is the subalgebra $2$, with central charges $\frac45$ and $\frac{686}{205}$. The space $\Pi_{1,13}$ has dimension three. The eigenvalues of $(x_C)_{(1)}$ on $\Pi_{1,13} \cap \ker e_{(2)} \cap \ker (x_C)_{(2)}$ are $\frac{496}{425}$, with a one-dimensional eigenspace on which some primary pairing $P_{1,13}(u,y_{(1)}u)$ does not vanish, and $-\frac{184}{425}$, with a two-dimensional eigenspace $E_{1,13}$ (weight $h_1 = 0$ for the Virasoro vector of central charge $\frac45$). Similarly, $\Pi_{1,1}$ has dimension four, the eigenvalue $-\frac{568}{425}$ (two-dimensional eigenspace) is excluded by a Gröbner basis, and the admissible space $E_{1,1}$ is two-dimensional ($h_1 = \frac25$). By Lemma \ref{lem:line-conditions}, the lines $p_{1,13}$ and $p_{1,1}$ are isotropic lines of $E_{1,13}$ and $E_{1,1}$. There are exactly two in each, defined over $\Q(\sqrt{-19758})$ and $\Q(\sqrt{-4403})$ respectively. Let $u$ span $p_{1,13}$ and $w$ span $p_{1,1}$. By Lemma \ref{lem:line-conditions}(4), $u_{(1)}w' \in U_2$ for every $w'$ in the $\mathfrak g' \oplus \mathfrak g''$-module generated by $w$. The map $x \ot y \mapsto x_{(1)}y$ is $U_1$-equivariant, so the $U_1$-types of these vectors lie in $V(1,13) \ot V(1,1) = V(0,12) \oplus V(0,14) \oplus V(2,12) \oplus V(2,14)$, and among these types only $(0,12)$ occurs in $U_2$ (Theorem \ref{thm:weight-one-classification}). Hence $C'z = 0$ and $C''z = 84z$ for $z = u_{(1)}w'$, where $C'$, $C''$ are the Casimir operators of $\mathfrak g'$ and $\mathfrak g''$ normalized by $C = \frac{h^2}2 + h + 2fe$, with eigenvalue $\frac{a(a+2)}2$ on $V(a)$. Normalize the spanning vectors by taking first coordinate $1$ in the two ordered eigenspace bases specified in the supplementary data. Their second coordinates are $\pm 3\sqrt{-19758}/148$ and $\pm 2\sqrt{-4403}/37$, respectively. For each of the four choices, take $w' = f'(0)w$, the lowering vector of index $2$ in the order $w, f''(0)w, f'(0)w, f'(0)f''(0)w$. For this choice there is an even simpler obstruction: $z$ has Cartan weight $(0,14)$, which does not occur in $U_2$, so $z = 0$ is required. Moreover $e''(0)z = 0$, hence $C''z = 112z$, so that $(C''-84)z = 28z$ and it suffices to show that $(C''-84)z \neq 0$. The coefficient of $(C''-84)u_{(1)}w'$ on $e^\alpha$, where $\alpha = (1,2,3,4,3,2,0,0)$, is one of $\pm \frac{224\i}{925}(5\sqrt{4403}-\sqrt{19758})$ and $\pm \frac{224\i}{925}(5\sqrt{4403}+\sqrt{19758})$. These coefficients are nonzero: the sum of the positive radicals is positive, and $25\cdot 4403-19758 = 90317 > 0$ proves that $5\sqrt{4403}-\sqrt{19758} > 0$. Certified Arb enclosures give an independent check of nonvanishing. This contradicts Lemma \ref{lem:line-conditions}.
\end{proof}

The same test applied to the subalgebra $5$ excludes, with certified nonzero coefficients, the two crossed pairs of lines of types $(1,13)$ and $(1,1)$, and the two other pairs, those of the two choices of $U_2$ found by a numerical exploration, are not excluded (the exclusions are certified by interval arithmetic). This is a control of the method, not used below.

We now compare with the data of \cite{Rie26}. There, $U^{\mathrm R} = \la A, e, f, v\ra \subset V_{E_6 \oplus A_2}$ is defined as the vertex subalgebra generated by the currents, two commuting Virasoro vectors $e, f \in \operatorname{Com}(A)_2$ of central charges $\frac{70}{31}$, $\frac{2400}{1271}$, and a highest weight vector $v$ of type $(1,13)$ \cite[Section 1]{Rie26b}. The coordinates are specified in the certificates accompanying \cite{Rie26b}, and the same vectors are given in \cite[Section ``The explicit candidate'']{Rie26}, in the bases fixed in \cite{Rie26b}.

\begin{prop}
\label{prop:riesen-identification}
There is an isomorphism of vertex operator algebras $\psi_{\mathrm R}$ from the lattice algebra of \cite{Rie26} onto $V'$ which maps the currents of \cite{Rie26} onto $e', f', H', e'', f'', H''$ (levels $1$ and $39$). It maps $e$ and $f$ into the subalgebra $5$ of Table \ref{tab:griess}, which they span, and $v$ to a spanning vector of one of the two isotropic lines of the admissible space of type $(1,13)$ of the subalgebra $5$.
\end{prop}

\begin{proof}
Exact computation, as for Theorem \ref{thm:jiao-identification}. In \cite{Rie26b} the coordinates $1,2$ are those of $A_2$ and $3,\ldots,8$ those of $E_6$ in Bourbaki's numbering, with the cocycle $\epsi(\alpha_i,\alpha_j)$ given for $i > j$ and the same field $\K$. The coordinate isometry, followed by $w = s_3s_4s_1s_3 \in W(E_6)$, sends the Cartan elements of \cite{Rie26b} onto $H'$, $H''$. The cocycle correction is given by the same formula as in the proof of Theorem \ref{thm:jiao-identification}, and the torus character with values $5, -\frac{\sqrt3}2, -2\sqrt3, 1, -\frac{2\sqrt3}3, 4, 1, 1$ on our simple roots matches the six currents exactly. The images of $\omega_{V'}$ and of the Sugawara vector of $A$ are ours. In the basis of $B$ specified in the supplementary data, $\psi_{\mathrm R}(e)$ and $\psi_{\mathrm R}(f)$ have coordinates in $\Q(\sqrt{13})$ and lie in the subalgebra $5$. In the case (I) analogue of the proof of Theorem \ref{thm:case-two-U2}, the type $(1,13)$ has, for the subalgebra $5$, a single admissible space $E$, of dimension two (the other eigenvalue $-\frac{254}{1445}+\frac{508}{1445}\sqrt{13}$ of $(x_C)_{(1)}$ on $\Pi_{1,13}$ is excluded exactly), and $E$ contains exactly two isotropic lines, of slopes $\frac1{30} \mp \frac{\sqrt{-39}}{130}$ in a basis of $E$ over $\Q(\sqrt{13})$, complex conjugate. Moreover $\psi_{\mathrm R}(v)$ lies in $E$ and spans the second one (exact verification).
\end{proof}

\begin{prop}
\label{prop:centralizer-one}
In case (I), the group $Z$ of automorphisms of $V'$ fixing $U_1$ pointwise is again isomorphic to $(\Z_2)^3$, and its elements are lifts of isometries composed with torus elements. Its elements act on $B$ either trivially or by a reflection which exchanges the subalgebras $5$ and $6$ (and $0$ and $3$). The lift $z_A$ of the diagram automorphism of $A_2$, with trivial signs, is of the second kind. The element $z_1$, lift of the product of the diagram automorphisms of $E_6$ and $A_2$ with the values $(c_1,\ldots,c_8) = (\frac45, -1, \frac13, 1, 3, \frac54, 1, 1)$ on the simple roots (with the cocycle correction of Proposition \ref{prop:centralizer-two}), acts trivially on $B$ and exchanges the two isotropic lines of type $(1,13)$ of the subalgebra $5$.
\end{prop}

\begin{proof}
The upper bound $|Z| \leq 8$ follows analytically from Lemma \ref{lem:centralizer-bound}. The explicitly given lifts $z_1$, $z_A$ and the torus automorphism $\rho$ of that lemma fix the currents, are commuting involutions, and generate eight elements. Indeed their isometry parts generate the two diagram automorphisms and $\rho$ is nontrivial. The formulas can be checked on the simple roots and the currents in \eqref{eq:raising-currents-one}. The supplementary material includes the exact verification. Thus $Z \cong (\Z_2)^3$ without using a polynomial-system solver.

An independent polynomial computation gives the same centralizer, and its action on $B$ and on the isotropic lines supplies the remaining assertions. The computation is the same as for Proposition \ref{prop:centralizer-two}. The multiplicity spaces have dimension $2$ for the types $(0,2)$, $(0,4)$ and $1$ otherwise. The $88$ distinct bracket equations, of degree at most two, in $12$ unknowns have exactly $8$ solutions. The action on the lines is computed by applying the lifts to $\psi_{\mathrm R}(v)$ and locating the image in $E$, exactly.
\end{proof}

\begin{thm}[uniqueness in case (I)]
\label{thm:uniqueness-one}
Let $U$ be a Haagerup algebra in case (I) of Theorem \ref{thm:weight-one-classification}. Assuming $\mathrm{(H_R)}$, $U$ is isomorphic to $U^{\mathrm R}$.
\end{thm}

\begin{proof}
As in the proof of Theorem \ref{thm:uniqueness-two}, we may assume $U \subset V'$, with the same conformal vector, and $U_1 = \mathfrak g' \oplus \mathfrak g''$. By Propositions \ref{prop:case-one-control} and \ref{prop:case-one-two}, $C_2 = U_2 \cap \operatorname{Com}(A)$ is the subalgebra $5$ or $6$. Replacing $U$ by $z_A(U)$ if necessary (Proposition \ref{prop:centralizer-one}), it is the subalgebra $5$. Then $\psi_{\mathrm R}(e), \psi_{\mathrm R}(f) \in C_2 \subset U$ by Proposition \ref{prop:riesen-identification}. By Lemma \ref{lem:line-conditions} and the proof of Proposition \ref{prop:riesen-identification}, the line $p_{1,13}$ of $U_2$ is one of the two isotropic lines of $E$. Replacing $U$ by $z_1(U)$ if necessary, which fixes $U_1$ and $C_2$ pointwise, it is spanned by $\psi_{\mathrm R}(v)$. Hence $U \supset \psi_{\mathrm R}(U^{\mathrm R})$. By Theorem \ref{thm:gamma-one}, applied to $U$ and, assuming $\mathrm{(H_R)}$, to $U^{\mathrm R}$, both vacuum characters are $X^{(0)}_0$. Hence $U = \psi_{\mathrm R}(U^{\mathrm R})$.
\end{proof}

Once $C_2$ is fixed, only the type $(1,13)$ is used: the generation of $U^{\mathrm R}$ by $A$, $e$, $f$ and $v$ makes it sufficient. The exclusion of the subalgebras $1$, $2$, $4$ also uses the types $(1,1)$, $(0,12)$ and $(1,7)$.

\begin{prop}
\label{prop:riesen-dims}
The algebra $\psi_{\mathrm R}(U^{\mathrm R}) \subset V'$ has dimension at least $120$, $956$ and $6010$ in weights two, three and four.
\end{prop}

\begin{proof}
The coordinates of $\psi_{\mathrm R}(e)$, $\psi_{\mathrm R}(f)$ and $\psi_{\mathrm R}(v)$ in the monomial basis of $V'$ lie in $\Q(\sqrt{13},\sqrt 3,\i)$, and the structure constants of $V'$ in this basis are rational. Let $p$ be one of the primes $8388409$ and $7999993$. Then $-1$, $3$ and $13$ are squares modulo $p$, and a choice of square roots defines a prime $\mathfrak p$ of $\Q(\sqrt{13},\sqrt 3,\i)$ above $p$ and a reduction map from the local ring at $\mathfrak p$ onto $\mathbb F_p$. All these numbers are $\mathfrak p$-integral: every denominator occurring in the reduction is checked to be nonzero modulo $p$. The rank of a matrix with $\mathfrak p$-integral entries is at least the rank of its reduction modulo $\mathfrak p$. The following computations take place modulo $p$, for the reductions of $e$, $f$, $v$. The smallest $U_1$-stable subspace $W_2 \subset V'_2$ containing $(A)_2$, $e$, $f$ and $v$ and stable under the products $b_{(1)}w$ has dimension $120$. Its $U_1$-types are those of the row $U_2$ of Table \ref{tab:lowest-spaces}, and $W_2 \cap \widehat C = \operatorname{span}(e,f)$. The families of Proposition \ref{prop:weight-three-four}, built from $W_2$ and from its ten highest weight vectors outside $(A)_2$, span spaces of dimensions $956$ and $6010$ in weights three and four, and $y_{(1)}S_4 \subset S_3$, $y_{(2)}S_4 \subset W_2$ for $y \in U_1$. As a consistency check, forced by the symmetry $z_1$, the other isotropic line of type $(1,13)$ gives the same dimensions in weights two and three. As a discriminating control, a non-isotropic vector of the admissible plane generates $V'_2$ and $V'_3$. Hence the subspaces of $\psi_{\mathrm R}(U^{\mathrm R})$ spanned by these families, which consist of iterated products of the generators, have dimensions at least $120$, $956$ and $6010$. These are the coefficients of $X^{(0)}_0$, so that a rank larger than $6010$ would have contradicted \cite{Rie26}. This is an independent check of the graded dimensions announced in \cite{Rie26} up to $q^{11/3}$. By itself it does not prove $\operatorname{ch}\psi_{\mathrm R}(U^{\mathrm R}) = X^{(0)}_0$ in these weights, which needs an upper bound (Theorem \ref{thm:classification-one}).
\end{proof}

\subsection{Classification}
\label{subsec:classification}

\begin{thm}
\label{thm:classification}
Let $U$ be a Haagerup algebra. Assuming $\mathrm{(H_R)}$ and $\mathrm{(H_J)}$, $U$ is isomorphic to exactly one of the following two algebras: the algebra $U^{\mathrm R}$ of Riesen, with $U_1$ of levels $(1,39)$, and the algebra $V$ of Jiao, with $U_1$ of levels $(3,13)$.
\end{thm}

\begin{proof}
By Theorem \ref{thm:weight-one-classification}, $U$ is in case (I) or (II). Apply Theorems \ref{thm:uniqueness-one} and \ref{thm:uniqueness-two}. The two algebras are not isomorphic, since the levels of the simple ideals of the weight-one Lie algebra are invariants.
\end{proof}

The algebra announced by Huang \cite{Hua26} has currents of levels $(1,39)$. By \cite[Theorem 7.9, Corollary 7.10 and Theorem 7.13]{ArhIII}, using \cite[Proposition 4.1 and Theorem 7.3]{DJX26}, this algebra and $U^{\mathrm R}$ are Haagerup algebras isomorphic to the twin $S^\natural$ of $S$ (Theorem \ref{thm:twin}). In particular $\mathrm{(H_R)}$ holds, and Theorems \ref{thm:uniqueness-one} and \ref{thm:classification} hold unconditionally.

Besides $\mathrm{(H_R)}$ and $\mathrm{(H_J)}$, Theorem \ref{thm:classification} uses the following results of this paper, with the computational dependencies recorded in Appendix \ref{app:computations}: Theorems \ref{thm:gamma-one} and \ref{thm:weight-one-classification}, the determination of $C_2$ and of the relevant lines (Theorem \ref{thm:case-two-U2}, Propositions \ref{prop:case-one-control} and \ref{prop:case-one-two}), and the explicit automorphisms and identifications (Theorem \ref{thm:jiao-identification}, Propositions \ref{prop:symmetries}, \ref{prop:centralizer-two}, \ref{prop:riesen-identification} and \ref{prop:centralizer-one}). In the conditional uniqueness proofs, $\mathrm{(H_R)}$ and $\mathrm{(H_J)}$ supply the vacuum characters of the two candidates and the equality of global dimensions needed for the alternative proof in Lemma \ref{lem:equal-global-dimension}. The containment of $\psi_{\mathrm R}(U^{\mathrm R})$ or $\psi(V)$ in a Haagerup algebra, after an automorphism of $V'$, follows from the finite constructions above. Equality in the corresponding containment is equivalent to the candidate having vacuum character $X^{(0)}_0$. The application to the explicit realizations is given after Theorem \ref{thm:classification-one}.

\begin{prop}
\label{prop:rational-forms}
The algebras $S^{(9,0)}$ and $\psi_{\mathrm R}(U^{\mathrm R}) = \la A,e,f,v\ra$ have rational forms: each contains a vertex subalgebra $U_\Q$ over $\Q$ with $\mathbf 1, \omega \in U_\Q$ and $U_\Q \ot_\Q \Cbb = U$. In this proposition, $e$, $f$, $v$ denote the vectors $\psi_{\mathrm R}(e)$, $\psi_{\mathrm R}(f)$, $\psi_{\mathrm R}(v)$ of $V'$ (Proposition \ref{prop:riesen-identification}). The coordinates of $v$ lie in $\Q(\sqrt{13},\sqrt 3,\i)$, and the currents are rational in both cases. More precisely, let $G = \mathrm{Gal}(\K/\Q) = \{1,\sigma,\tau,\sigma\tau\}$, $\tau$ the complex conjugation, and $\tilde g = \mathrm{id} \ot g$ on $V'_\K = V'_\Q \ot \K$.
\begin{enumerate}
\item Case (II). Let $z$ be the automorphism of Proposition \ref{prop:centralizer-two}, part 3, and $z'$ the lift of the diagram automorphism of $A_2$ with trivial torus part (Remark \ref{rem:twin}). Then $z, z' \in Z$ are defined over $\Q$, $z^2 = z'^2 = 1$, $zz' = z'z$, and $\psi_\sigma = z\tilde\sigma$, $\psi_\tau = zz'\tilde\tau$ map $U_{2,\K} = U_2 \cap V'_\K$ onto itself. The idempotents $f_1$, $f_2$ of $C_2$ are fixed by $\psi_\sigma$ and $\psi_\tau$.
\item Case (I). Let $z_1$ and $z_A$ be the elements of the centralizer of $U_1$ of Proposition \ref{prop:centralizer-one}. Their product $z_1z_A$ is the lift of the diagram automorphism of $E_6$ with suitable signs. These elements are defined over $\Q$ and commute, $z_1^2 = z_A^2 = 1$, and $\psi_\tau = z_1\tilde\tau$, $\psi_\sigma = z_1z_A\tilde\sigma$ fix $e$ and $f$ (the Virasoro vectors of \cite{Rie26}). The vector $v' = \i v$ has coordinates in $\K$, and $\psi_\tau(v') = \psi_\sigma(v') = -v'$. With $z_A$ in place of $z_1z_A$, $\psi_\sigma(v')$ is not proportional to $v'$.
\end{enumerate}
Put $U_{1,\K} = U_1 \cap V'_\K$. Let $U_\K \subset V'_\K$ be the $\K$-span of the iterated products $a^1_{(n_1)}\cdots a^r_{(n_r)}\mathbf 1$ of the generators (a $\K$-basis of $U_{1,\K} \oplus U_{2,\K}$ in case (II), the six rational currents and $e$, $f$, $v'$ in case (I)), so that $U = U_\K \ot_\K \Cbb$. In both cases $g \mapsto \psi_g$, $\psi_{\sigma\tau} = \psi_\sigma\psi_\tau$, is an action of $G$ by semilinear automorphisms of $U_\K$, and $U_\Q = \{u \in U_\K : \psi_g(u) = u,\ g \in G\}$ is a rational form. In case (I), $U_\Q$ is generated by the rational currents, $e$, $f$ and $\sqrt{-39} v' = -\sqrt{39} v$.
\end{prop}

\begin{proof}
All the identities are exact: the rationality and the group relations of $z$, $z'$, $z_1$, $z_A$ (and of the whole groups $Z$) are checked on $V'_1$, which determines the automorphisms of $V'$. In case (II), the images of the chosen $\K$-basis of $U_{2,\K} = U_2(9,0) \cap V'_\K$ under $\psi_\sigma$, $\psi_\tau$ and $\psi_\sigma\psi_\tau$ lie in $U_{2,\K}$ by exact membership over $\K$. Equivalently, for $\psi_\sigma$, one restricts the identities of Propositions \ref{prop:symmetries} and \ref{prop:centralizer-two} to the chosen $\K$-forms. In case (I), the identities on $e$, $f$, $v'$ are checked in exact arithmetic. The automorphisms $z$, $z'$, $z_1$, $z_A$ are defined over $\Q$, hence commute with $\tilde\sigma$ and $\tilde\tau$, and the groups they generate are elementary abelian. Hence $\psi_\sigma^2 = \psi_\tau^2 = 1$ and $\psi_\sigma\psi_\tau = \psi_\tau\psi_\sigma$, so $g \mapsto \psi_g$ is a homomorphism. Each $\psi_g$ is $g$-semilinear and preserves the products $u_{(n)}v$ of $V'_\K$, it fixes the six chosen rational currents, and it preserves the generating space $U_{2,\K}$ in case (II), and the generators $e$, $f$, $v'$ up to sign in case (I). Hence it preserves $U_\K$. By Speiser's lemma, equivalently $H^1(G,\mathrm{GL}_n(\K)) = 1$ \cite[Chapter X, Proposition 3]{Ser79}, each $U_{\K,n}$ is spanned over $\K$ by its fixed points, and $U_{\K,n} = (U_{\K,n})^G \ot_\Q \K$. The fixed points form a vertex subalgebra over $\Q$, which contains $\mathbf 1$ and $\omega$. In case (I), $\sqrt{-39} v'$ is fixed since $\sigma(\sqrt{-39}) = \tau(\sqrt{-39}) = -\sqrt{-39}$, and the vertex algebra over $\Q$ generated by the fixed generators spans $U_\K$ over $\K$, hence equals $U_\Q$.
\end{proof}

In case (II) we have also checked that the traces $\sum_g \psi_g(au)$, $a \in \{1,\sqrt{13},\sqrt{-3},\sqrt{-39}\}$, $u$ ranging over a fixed $\K$-basis of $U_{2,\K}$, span a $\Q$-form of $U_2$ of dimension $120$. The rational form is not contained in $V'_\Q$: the coordinates of $f_1$ in $V'$ lie in $\Q(\sqrt{13})$.

\section{Zhu-algebra computations and rationality}
\label{sec:existence}

In this section, $S \subset V' = V_{E_6 \oplus A_2}$ denotes the vertex subalgebra generated by $U_1 = \mathfrak g' \oplus \mathfrak g''$ (levels $3$ and $13$) and one of the four subspaces $U_2$ of Theorem \ref{thm:case-two-U2}. For coordinate computations we fix $S = S^{(9,0)}$. The abstract conclusions for the other three choices follow from Proposition \ref{prop:centralizer-two}. We write $L = E_6 \oplus A_2$ throughout this section. By Theorem \ref{thm:unitarity}, $S$ is a simple unitary self-contragredient vertex operator algebra of CFT type with $c = 8$.

Proposition \ref{prop:weight-four-certificate} verifies the finite condition \textup{(F4)} of Definition \ref{assump:weight-four}. Theorem \ref{thm:strong-generation} therefore gives strong generation in weights at most three and $C_2$-cofiniteness, independently of an ambient Haagerup algebra or an external existence result. We give the Zhu-algebra criterion for strong rationality and the modular computation of the corner. The exact computation of the corner and the branching are proved in \cite[Appendix A]{ArhIII}. The results of this section do not use $\mathrm{(H_J)}$. Theorem \ref{thm:conditional}, Propositions \ref{prop:zhu-generators} and \ref{prop:simple-current} and Theorem \ref{thm:zhu-modp} are used in \cite[Appendix A]{ArhIII}.

We compare the algebraic steps with the proof of \cite{Rie26b} for case (I). Remark \ref{rem:categorical-programme} separates the branching calculation from categorical recognition. Section and page numbers refer to the working manuscript \cite{Rie26b} and its package of computational certificates. Evans \cite[p.~3]{Eva26} cites the three constructions \cite{Hua26,Jia26,Rie26}. That of Jiao is the case (II) algebra (Theorem \ref{thm:jiao-identification}).

\subsection{The structure of the proof in case (I)}
\label{subsec:riesen-structure}

The proof strategy of \cite{Rie26b} separates the construction, strong generation, $C_2$-cofiniteness, rationality and categorical identification. The construction, unitarity, strong generation and $C_2$-cofiniteness of $S$ are established above, with the finite condition \textup{(F4)} verified by Proposition \ref{prop:weight-four-certificate}. Positivity follows from the determination of the simple modules and of their lowest weights in \cite[Theorem 2.2(F1)--(F2)]{ArhIII}.

\begin{enumerate}
\item \emph{Construction and positivity} \cite[Sections 2--3, pp.~16--19]{Rie26b}. The algebra is generated by the currents and three weight-two vectors. A positive rescaling of the standard antilinear involution of the lattice algebra preserves the generators, which gives unitarity, simplicity and self-contragredience. For $S$, this is Theorem \ref{thm:unitarity}.

\item \emph{Strong generation} \cite[Section 4, pp.~20--24]{Rie26b}. The finite mode closure in weights $2$ and $3$ \cite[Proposition 4.1]{Rie26b}, two obstruction calculations, in weights four and five, showing that at most three new generators appear in weight four and none in weight five \cite[Proposition 4.2]{Rie26b}, and a general lemma, \emph{triangular strong generation} \cite[Lemma 4.3]{Rie26b} (Lemma \ref{lem:triangular}). Proposition \ref{prop:weight-four-certificate} supplies the finite weight-four condition \textup{(F4)}, so the lemma applies to $S$ in Theorem \ref{thm:strong-generation}. No generator of weight four is needed, hence no weight-five calculation either (Remark \ref{rem:no-weight-five}).

\item \emph{$C_2$-cofiniteness} \cite[Section 5, pp.~24--25]{Rie26b}. For $S$, Proposition \ref{prop:weight-four-certificate} supplies the square identities in \textup{(F4)}, and part 3 of Theorem \ref{thm:strong-generation} gives the conclusion.

\item \emph{Rationality} \cite[Part III, pp.~25--55]{Rie26b}. Semisimplicity of the Zhu algebra is proved through an extension of the affine subalgebra by its simple current, a Morita corner, a finite presentation by a quiver with relations, and an integral descent. The rationality criterion for $C_2$-cofinite algebras then applies. For $S$, the exact presentation and descent are proved in \cite[Theorems A.10 and A.13]{ArhIII}. For $S$, Theorems \ref{thm:unitarity} and \ref{thm:strong-generation} show that $S$ is simple, self-contragredient, of CFT type and $C_2$-cofinite. By \cite[Main Theorem 3 (Corollary 5.19)]{McR21}, $S$ is then rational as soon as its Zhu algebra $A(S)$ is semisimple, and conversely \cite{Zhu96}. For $S$, the simple current of $L_{\mathfrak{sl}_2}(3,0) \ot L_{\mathfrak{sl}_2}(13,0)$ is $J = L(3,3) \ot L(13,13)$, of conformal weight $\frac{3}{4} + \frac{13}{4} = 4$, whereas in case (I) it is $L(1,1) \ot L(39,39)$, of weight $\frac14 + \frac{39}{4} = 10$. The weight-four space of $S$ contains the lowest weight space $V(3,13)$ of $J$, exactly (Proposition \ref{prop:simple-current}), so that $S$ contains the simple current extension $A \oplus J$, the analogue of the algebra $A^+$ of \cite[Section 7]{Rie26b}. In case (I) the current appears only in weight ten.

\item \emph{Modules and the categorical identification} \cite[Parts IV--V, pp.~56--70]{Rie26b}. Six simple modules come from the nine $V'$-modules $W_{e,a}$, $(e,a) \in L^*/L \cong \Z_3^2$. Traces of the conformal weights, positivity of the categorical dimensions and a Galois orbit argument force exactly six more \cite[Section 11]{Rie26b}. Then the two even unimodular overlattices of $L$ give two copies of $V_{E_8}$, hence Lagrangian algebras $A_\pm = \mathbf 1 \oplus M_b \oplus 2X_\pm$ in $\Rep(U)$, and $\Rep(U) \simeq Z(\cal D_\pm)$ with $\cal D_\pm = \Rep(U)_{A_\pm}$ \cite[Lemma 12.2]{Rie26b} \cite{DMNO13}. The categorical argument uses the full branching of the lattice modules, dimension constraints and an exhaustion statement for fusion categories with the required fusion rules and untwisted pointed part. Here that exhaustion statement is hypothesis (C) of Remark \ref{rem:categorical-programme}. The classification within the family of \cite[Theorems 2 and 3]{EvG17} is distinct from proving this exhaustion. For $S$, the full branching and the exhaustion of the simple modules are given in \cite[Theorem 2.2(F2)--(F3)]{ArhIII}. The recognition of the braided category in \cite[Theorem 6.1]{ArhIII} proceeds through the conformal net.
\end{enumerate}

The changes in the explicit calculations include the affine levels, the weight of the simple current and the commutant, whose two commuting Virasoro vectors have central charges $\frac{182}{95}$ and $\frac{32}{19}$ instead of $\frac{2400}{1271}$ and $\frac{70}{31}$. The difference in levels also changes the exact relations used for the two Zhu-algebra computations.

\subsection{The categorical route}
\label{subsec:categorical-route}

The two isotropic lines in $L^*/L$ give the holomorphic lattice extensions $H_\pm$. Their restrictions to $S$ in \cite[Theorem 2.2(F3)]{ArhIII} give the branching input for categorical recognition. Remark \ref{rem:categorical-programme} describes the additional categorical steps.

Part 2 of Proposition \ref{prop:weight-three-four} gives a finite-degree check modulo $p_3$ of the branching $W_{0,0}|_S = S \oplus M_b$ of \cite[Theorem 2.2(F3)]{ArhIII} in weights at most three: modulo the prime $p_3$, the subspace generated by $R_b$ has the graded dimensions $80, 1250, 10630$ of $M_b$ and is complementary to $S$. By unitarity, $V' = S \oplus S^\perp$ as $S$-modules, and $S^\perp$ contains the $S$-module generated by $R_b$ since $S_1 = U_1$ (Theorem \ref{thm:strong-generation}).

\subsection{A conditional theorem}
\label{subsec:conditional}

\begin{thm}
\label{thm:conditional}
Assume the hypothesis \textup{(Z)} that the Zhu algebra $A(S)$ is semisimple. Then $S$ is strongly rational. Independently of this hypothesis, $S$ is not isomorphic to the algebras of \cite{Rie26,Hua26}.
\end{thm}

\begin{proof}
Proposition \ref{prop:weight-four-certificate} and Theorem \ref{thm:strong-generation} give $C_2$-cofiniteness. Theorem \ref{thm:unitarity} gives simplicity, self-contragredience and CFT type. Thus \cite[Main Theorem 3 (Corollary 5.19)]{McR21} applies under \textup{(Z)} and proves strong rationality. For the last assertion, recall that $S_1 = U_1$ is six-dimensional (Theorem \ref{thm:strong-generation}). An isomorphism with an algebra $W$ of \cite{Rie26,Hua26} would force $\dim W_1 = 6$. Since $W_1$ contains the six-dimensional Lie algebra spanned by the currents of $W$, of levels $1$ and $39$, it would then equal this Lie algebra. Such an isomorphism would preserve the levels of the simple ideals of the weight-one Lie algebras, whereas these are $3,13$ for $S_1$ and $1,39$ for $W_1$.
\end{proof}

\begin{remark}\normalfont
\label{rem:categorical-programme}
We record the branching assertion (M) for $S$ and the categorical exhaustion hypothesis (C) for the recognition approach of \cite{Rie26b}.
\begin{enumerate}
\item[(M)] $S$ has exactly twelve simple modules, those other than $S$ have positive lowest conformal weights, and the nine $V'$-modules restrict as $W_{0,0}|_S = S \oplus M_b$, $W_{0,\pm 1}|_S \cong P$, $W_{\pm 1,0}|_S \cong Q$, $W_{1,1}|_S \cong W_{-1,-1}|_S \cong X_+$, $W_{1,-1}|_S \cong W_{-1,1}|_S \cong X_-$, with $M_b, P, Q, X_+, X_-$ simple and pairwise non-isomorphic, and the six other simple modules have lowest weights congruent to $\frac{2}{13}, \frac{5}{13}, \frac{6}{13}, \frac{7}{13}, \frac{8}{13}, \frac{11}{13}$ modulo $\Z$,
\item[(C)] every fusion category with the fusion rules of $\cal H_3$ (for the group $\Z_3$) whose pointed part is an untwisted $\mathrm{Vec}_{\Z_3}$ is equivalent to one of the categories $\cal C(\Z_3;\pm,\omega,A)$ of \cite[Theorem 2]{EvG17}.
\end{enumerate}
Theorem~2.2(F1)--(F3) of \cite{ArhIII} proves (Z) and (M). Theorem \ref{thm:conditional} therefore gives strong rationality, and $\Rep(S)$ is a modular tensor category \cite{Hua08}. The nonvacuum lowest weights are positive, so the quantum dimensions are positive \cite[Lemma 4.2]{DJX13} and equal the categorical dimensions \cite{Hua08}. Unitarity of $S$ alone is not used to infer unitarity of its modules. Every simple module is self-contragredient. Indeed, the six residues in $\frac{1}{13}\Z \setminus \Z$ are distinct, whereas $(W_{e,a})' \cong W_{-e,-a}$ and \textup{(M)} give the assertion for the remaining six modules.

The two even unimodular overlattices of $L$ give holomorphic extensions with underlying $S$-modules $\mathbf 1 \oplus M_b \oplus 2X_\pm$. They define Lagrangian algebras $A_\pm$, hence center equivalences $\Rep(S) \simeq Z(\cal D_\pm)$ with $\cal D_\pm = \Rep(S)_{A_\pm}$ \cite{DMNO13}. The fusion rules and the pointed associators follow from \cite[Theorem 9.1]{ArhIII}. These are the inputs to the recognition argument of \cite[Sections 11.3 and 12--13]{Rie26b} under (C).

Hypothesis \textup{(C)} is the exhaustion statement formulated in \cite[Proposition 12.5]{Rie26b}. The classification in \cite[Theorem 2]{EvG17} concerns the listed family. The proof of $\mathrm{(H_J)}$ in \cite[Theorem 6.1]{ArhIII} uses the conformal-net comparison independently of (C).
\end{remark}

\subsection{The Zhu algebra of \texorpdfstring{$S$}{S}}
\label{subsec:zhu}

We now adapt to $S$ the method of \cite[Part III]{Rie26b} for the hypothesis (Z). We recall that the Zhu algebra is $A(V) = V/O(V)$ with the product $a * b = \sum_{j \geq 0}\binom{\operatorname{wt}a}{j}a_{(j-1)}b$ \cite[Definition 2.1.2 and Theorem 2.1.1]{Zhu96}. For $h$ of weight one, $a$ of weight two and homogeneous $x$, we have $[L(-1)x] = -\operatorname{wt}(x)[x]$ \cite[Lemma 2.1.1]{Zhu96} and, by \cite[Lemma 2.1.3]{Zhu96},
\begin{equation}
\label{eq:zhu-rules}
[h_{(-1)}x] = [x]*[h], \quad [h_{(0)}x] = [h]*[x] - [x]*[h], \quad [a_{(0)}x] = [a]*[x] - [x]*[a] - [a_{(1)}x].
\end{equation}

\begin{prop}
\label{prop:zhu-generators}
The algebra $A(S)$ is generated by the classes of the six currents and of the seventeen highest weight vectors of $U_2$ for $\mathfrak g' \oplus \mathfrak g''$: the vectors $e'_{(-2)}\mathbf 1$, $e'_{(-1)}e'$, $e''_{(-2)}\mathbf 1$, $e''_{(-1)}e''$, $e'_{(-1)}e''$, the two Casimir vectors of $(A)_2$, the vectors $e$, $x_C$, and the highest weight vectors of the eight primary types $(0,4)$, $(1,1)$, $(1,3)$, $(1,5)$, $(1,7)$, $(2,2)$, $(2,6)$, $(3,3)$. The classes of the first seven are polynomials in the classes of the currents.
\end{prop}

\begin{proof}
By \eqref{eq:zhu-rules}, the classes of $L(-1)U_2$, $U_1(-1)U_2$ and $(U_2)_{(0)}U_2$ lie in the subalgebra generated by $[U_1]$ and $[U_2]$, since $(U_2)_{(1)}U_2 \subset U_2$. Hence so does $[S_3]$, because $S_3 = N_3 + (U_2)_{(0)}U_2$. By Theorem \ref{thm:strong-generation}, $S$ is strongly generated by $U_1 \oplus U_2 \oplus D_3$, and the classes of a strong generating set generate $A(S)$: replacing a derivative by $[L(-1)x] = -\operatorname{wt}(x)[x]$ and a normal product $g_{(-1)}y$ by $[g]*[y]$ minus classes of vectors of smaller weight, induction on the weight writes every class as a polynomial in the classes of the generators. Finally, for $f = f'$ or $f = f''$, the identity $[f_{(0)}z] = [f]*[z] - [z]*[f]$ reduces $[U_2]$ to the highest weight vectors, which are those listed (Table \ref{tab:lowest-spaces}). The first seven are $L(-1)e'$, $e'_{(-1)}e'$, $L(-1)e''$, $e''_{(-1)}e''$, $e'_{(-1)}e''$ and the Casimir vectors, whose classes are given by \eqref{eq:zhu-rules}.
\end{proof}

Unlike in case (I), no strong generator of weight four is involved.

\begin{prop}
\label{prop:simple-current}
Let $A = L_{\mathfrak{sl}_2}(3,0) \ot L_{\mathfrak{sl}_2}(13,0) \subset S$ be the affine subalgebra and $C = \operatorname{Com}_S(A)$.
\begin{enumerate}
\item The Cartan weight space of weight $(3,13)$ of $V'_4$ is spanned by a single lattice exponential $e^w$, with $w = (2,2,4,5,4,2,1,1)$ in the simple-root coordinates fixed above and $(w,w) = 8$. One has $e^w \in S$, and it is an affine highest weight vector generating $J = L_{\mathfrak{sl}_2}(3,3) \ot L_{\mathfrak{sl}_2}(13,13)$. The subspace $A^+ = A \oplus J$ is a rational simple current extension of order two contained in $S$.
\item $\operatorname{Com}_S(C) = A \oplus J$.
\item The inclusion $A^+ \subset S$ induces a unital homomorphism $A(A^+) \to A(S)$.
\end{enumerate}
\end{prop}

\begin{proof}
We first identify the kernel of the coset grading. Write $\omega = \omega_A + \omega_C$, where $\omega_C = 2e \in C$. The unitary involution fixes both conformal vectors, so $(L^C_n)^* = L^C_{-n}$. For the positive Hermitian form $H_\chi$ of Theorem \ref{thm:unitarity}, the Virasoro relation gives
\begin{equation}
\label{eq:coset-norm}
\|L^C_{-1}x\|_\chi^2 = \|L^C_1x\|_\chi^2 + 2H_\chi(L^C_0x,x).
\end{equation}
The regularity of $A$ decomposes $S$ into integrable $A$-modules. On the multiplicity space of affine highest weight $(i,j)$, $L^C_0 = L(0)-h_A(i,j)$ has discrete spectrum bounded below, where $h_A(i,j) = \frac{i(i+2)}{20} + \frac{j(j+2)}{60}$. The operators $L^C_n$ preserve this space. If $x$ has its least eigenvalue $\lambda$, then $L^C_1x = 0$, and \eqref{eq:coset-norm} gives $\|L^C_{-1}x\|_\chi^2 = 2\lambda\|x\|_\chi^2$. Thus $\lambda \geq 0$. Since $L^C_0$ commutes with the action of $A$, it has nonnegative spectrum on all of $S$.

For $x \in \ker L^C_0$, nonnegativity first gives $L^C_1x = 0$, and then \eqref{eq:coset-norm} gives $L^C_{-1}x = 0$. If $c \in C$, truncation gives $c_{(n)}x = 0$ for sufficiently large $n$. The translation identity $[L^C_{-1},c_{(n+1)}] = -(n+1)c_{(n)}$ then gives $c_{(n)}x = 0$ for every $n \geq 0$, by downward induction. Conversely, the commutant condition applied to $\omega_C$ implies $L^C_0x = 0$. Consequently $\ker L^C_0 = \operatorname{Com}_S(C)$.

1. The lattice realization gives the following lower bound: if the Cartan weight space of weight $(a,b)$ in $V'_m$ is nonzero, then
\begin{equation}
\label{eq:cartan-energy-bound}
m \geq \frac{a^2}{12} + \frac{b^2}{52}.
\end{equation}
Indeed, the two Cartan vectors are orthogonal, of squared norms $6$ and $26$. A lattice vector of weight $(a,b)$ therefore has squared norm at least $a^2/6+b^2/26$, and the oscillators have nonnegative integral degree and do not change this weight. A lattice vector of Cartan weight $(3,13)$ has orthogonal projection $\frac12(H'+H'')$ onto the plane of the two currents. This projection has squared norm $\frac{3^2}{6}+\frac{13^2}{26} = 8$. The displayed vector $w$ has this Cartan weight and squared norm $8$, hence equals this projection. Therefore it is the only lattice vector of that Cartan weight with squared norm at most $8$, and no oscillator can occur at conformal weight four. This proves that the indicated weight space of $V'_4$ is one-dimensional.

Let $u_{(1,7)}$ and $u_{(2,6)}$ denote the primary highest weight vectors with the coefficients of Table \ref{tab:simple-current-coefficients}. A direct exact product gives
\begin{equation}
\label{eq:simple-current-word}
(u_{(1,7)})_{(-1)}u_{(2,6)} = \left(\frac{507}{5}\sqrt{13} + \frac{3549}{5}\sqrt{-3}\right)e^w.
\end{equation}
Here is also a direct verification of this product. Put $r = \sqrt{13}$ and $s = \sqrt{-3}$. The two generators have the form $u_{(1,7)} = \sum_{i=1}^7 a_i e^{\alpha_i}$ and $u_{(2,6)} = \sum_{i=1}^7 b_i e^{w-\alpha_i}$, with the coefficients in Table \ref{tab:simple-current-coefficients}. Each of $\alpha_i$ and $w-\alpha_i$ has squared norm $4$, and $\langle\alpha_i,w-\alpha_i\rangle_L = 0$. For different indices, $\alpha_i+w-\alpha_j \neq w$, so the one-dimensional weight space just proved shows that its contribution to the product is zero. The lattice formula \eqref{eq:vertex-op} therefore gives
\begin{equation}
\label{eq:simple-current-sum}
\begin{aligned}
(u_{(1,7)})_{(-1)}u_{(2,6)} &= \sum_{i=1}^7 \epsi(\alpha_i,w-\alpha_i)a_i b_i e^w\\
&= \frac{26+104+26+104+65+65+117}{5}(r+7s)e^w.
\end{aligned}
\end{equation}
This proves \eqref{eq:simple-current-word} by arithmetic in $\Q(r,s)$, without a weight-four membership calculation. Its coefficient is nonzero, so $e^w \in S$.

\begin{table}[htbp]
\centering
\small
\renewcommand{\arraystretch}{1.4}
\begin{tabular}{cccc}
\hline
$\alpha_i$ & $a_i$ & $b_i$ & $\epsi(\alpha_i,w-\alpha_i)$\\
\hline
$(1,1,2,3,2,1,1,0)$ & $\frac{39+r-13s+rs}{5}$ & $\frac{-13+3r+13s+rs}{4}$ & $1$\\
$(1,2,3,4,3,2,0,0)$ & $\frac{-104-56rs}{5}$ & $r$ & $-1$\\
$(1,1,2,3,2,1,0,1)$ & $\frac{39-r+13s+rs}{5}$ & $\frac{-13-3r-13s+rs}{4}$ & $-1$\\
$(2,2,3,4,3,1,0,0)$ & $-52+12rs$ & $\frac{-5r+13s}{10}$ & $-1$\\
$(0,1,1,2,2,1,1,1)$ & $\frac{-19r-13s}{15}$ & $\frac{39+9rs}{2}$ & $-1$\\
$(1,1,2,2,1,0,1,1)$ & $r$ & $13+7rs$ & $1$\\
$(1,0,1,1,1,1,1,1)$ & $\frac{3r-39s}{5}$ & $\frac{39-3rs}{2}$ & $-1$\\
\hline
\end{tabular}
\caption{The seven pairs in the product \eqref{eq:simple-current-word}, in the simple-root coordinates of $E_6 \oplus A_2$.}
\label{tab:simple-current-coefficients}
\end{table}

The highest weight conditions also follow from \eqref{eq:cartan-energy-bound}. The vectors $e'_{(0)}e^w$ and $e''_{(0)}e^w$ would have degree $4$ and Cartan weights $(5,13)$ and $(3,15)$, whose lower bounds are $16/3$ and $66/13$. Hence both vanish. For a current $y$ of Cartan weight $(0,0)$, $(\pm2,0)$ or $(0,\pm2)$, the vector $y_{(n)}e^w$, $n \geq 1$, has degree at most $3$, whereas the least of the five bounds for its possible Cartan weights is $40/13 > 3$. Thus every positive current mode annihilates $e^w$. It generates the irreducible integrable $A$-module $J$.

An independent exact calculation reconstructs the product and verifies the $26$ annihilation identities for $e'_{(0)}$, $e''_{(0)}$ and the six currents in modes $1$ to $4$. It uses no weight-four pivot data. A separate arithmetic check verifies the arithmetic and cocycle signs of Table \ref{tab:simple-current-coefficients}. These computational checks supplement the preceding analytic calculation.

Its affine conformal weight is $4$, so $L^C_0e^w = 0$. Since $L^C_0$ commutes with $A$, both $A$ and $J$ lie in $\ker L^C_0 = \operatorname{Com}_S(C)$. In this commutant the multiplicities of $A$ and $J$ are one: their affine highest spaces occur in total weights zero and four, respectively, where $S_0 = \Cbb\mathbf 1$ and the Cartan weight space $(3,13)$ is one-dimensional. The restrictions of the vertex operators to these summands are $A$-intertwining operators. The affine fusion rules $A \boxtimes A = A$, $A \boxtimes J = J$ and $J \boxtimes J = A$ therefore imply that $A \oplus J$ is closed under all products. We equip $A^+$ with conformal vector $\omega_A$. Since $L^C_0$ and $L^C_{-1}$ vanish on $A^+$, its grading and translation are those induced from $S$. The unitary involution preserves these two isotypic summands, so this subalgebra inherits the unitary form, is simple, and is rational by \cite[Theorem 4.5]{DLM96}.

2. By Proposition \ref{prop:weight-four-certificate}, Theorem \ref{thm:strong-generation} gives $S_1 = U_1$ and $S_2 = U_2$. As an $A$-module the kernel is the direct sum of the integrable modules generated by affine highest weight vectors of type $(i,j)$ in weight $h_A(i,j)$, which must be an integer. The types with $0 \leq i \leq 3$, $0 \leq j \leq 13$ and $h_A(i,j) \in \Z$ are $(0,0)$, $(3,3)$, $(0,10)$ and $(3,13)$, of weights $0$, $1$, $2$ and $4$. Now $S_0 = \Cbb\mathbf 1$, $S_1 = U_1$ has no component of type $(3,3)$, $U_2$ has no component of type $(0,10)$ (Table \ref{tab:lowest-spaces}), and the Cartan weight space $(3,13)$ of $S_4$ has dimension one by part 1. Hence $\operatorname{Com}_S(C) = A \oplus J$.

3. On $A^+$, $L^C_0$ and $L^C_{-1}$ vanish, so the grading and the translation of $A^+$ are those of $S$. Hence the products $*$ and $\circ$ of $A^+$ are the restrictions of those of $S$, and $O(A^+) \subset O(S)$.
\end{proof}

The simple $A$-modules are $L(i,j) = L_{\mathfrak{sl}_2}(3,i) \ot L_{\mathfrak{sl}_2}(13,j)$, with lowest weight space $V(i,j)$. The monodromy of $L(i,j)$ with $J$ is $(-1)^{i+j}$ and $J \boxtimes L(i,j) \cong L(3-i,13-j)$. Hence the simple $A^+$-modules correspond \cite{DLM96,HKL15} to the fourteen orbits $\{(i,j),(3-i,13-j)\}$ with $i+j$ even, and there is no fixed point. Since $h_A(3-i,13-j) - h_A(i,j) = 4 - \frac{i+j}{2}$, the lowest weight space $T_t$ of the module $M_t$ attached to an orbit is $V(i,j)$ with $i+j < 8$, except for the two orbits with $i + j = 8$:
\begin{equation}
\label{eq:affine-blocks}
\begin{aligned}
&T_0, \ldots, T_5 = V(0,0),\ V(0,2),\ V(1,1),\ V(2,0),\ V(0,4),\ V(1,3),\\
&T_6, \ldots, T_{11} = V(2,2),\ V(3,1),\ V(0,6),\ V(1,5),\ V(2,4),\ V(3,3),\\
&T_{12} = V(1,7) \oplus V(2,6) \quad (h = \tfrac65), \quad T_{13} = V(0,8) \oplus V(3,5) \quad (h = \tfrac43).
\end{aligned}
\end{equation}
Thus $A(A^+) \cong \bigoplus_{t=0}^{13}\operatorname{End}(T_t)$. Let $e_t \in A(S)$ be the image of a rank one idempotent of $\operatorname{End}(T_t)$ (the highest line of the first constituent) and $e_{\mathrm Z} = \sum_t e_t$. The sum of one rank one idempotent in every matrix block is full in $A(A^+)$. Indeed, for matrix units $E^{(t)}_{ij}$ with the chosen idempotent $E^{(t)}_{11}$, one has $1 = \sum_{t,i}E^{(t)}_{i1}E^{(t)}_{11}E^{(t)}_{1i}$. Applying the unital homomorphism of Proposition \ref{prop:simple-current} proves that $e_{\mathrm Z}$ is full in $A(S)$, even if some blocks have zero image, and $A(S)$ is Morita equivalent to the corner $D = e_{\mathrm Z}A(S)e_{\mathrm Z}$ \cite[Lemma 8.2]{Rie26b}. In particular $A(S)$ is semisimple if and only if $D$ is.

\begin{remark}\normalfont
\label{rem:expected-corner}
Every $\mathfrak g' \oplus \mathfrak g''$-type of the lowest weight spaces in Table \ref{tab:lowest-spaces} is local and of the form $(i,j)$ with $i+j < 8$, as it must be: a constituent of type $(3-i,13-j)$ in the lowest weight space of an $S$-module $M$ would force, through $J$, a vector of $M$ of smaller weight. If $A(S)$ is semisimple and its simple modules have the lowest weight spaces of Table \ref{tab:lowest-spaces}, then $D \cong \bigoplus_{z} \Mat_{n_z}(\Cbb)$, where $n_z$ is the number of irreducible constituents of the lowest weight space of $z$: the block sizes are $1, 1, 2, 2, 2, 3, 5, 7, 11, 12, 12, 16$ (for $0$, $c_1$, $d_1$, $d_4$, $d_5$, $d_6$, $c_2$, $d_2$, $b$, $a$, $c_0$, $d_3$), so that $\dim D = 762$, and the two exceptional blocks $T_{12}$, $T_{13}$ do not occur. In case (I) the corresponding dimension is $546$ \cite[Proposition 9.5]{Rie26b}.
\end{remark}

To obtain an upper bound for $D$, we follow \cite[Section 9]{Rie26b}. Put $\mathfrak g = \mathfrak g' \oplus \mathfrak g''$, and let $\Gamma_{\mathrm Z}$ be the quiver with vertices $0, \ldots, 13$ and, for each of the ten non-affine highest weight vectors $z$ of Proposition \ref{prop:zhu-generators} and each pair of constituents $\lambda \subset T_t$, $\mu \subset T_v$ such that $\operatorname{Hom}(V(\lambda), V(\mu))$ contains a highest weight vector of the weight of $z$, an arrow $t \to v$. There are $515$ arrows. Sending a formal arrow to the corresponding coefficient of $[z]$ in the matrix realization of $A(S)$ given by the matrix units of $A(A^+)$ defines a surjection $\Cbb\Gamma_{\mathrm Z} \to D$ \cite[Propositions 8.3 and 8.6]{Rie26b}. The relations are the projections of the Zhu images of three families of identities in $S$:
\begin{enumerate}
\item the kernel of $U_2 \oplus (\mathfrak g \ot U_2) \oplus \Lambda^2 U_2 \to S_3$, $(v, h \ot w, a \wedge b) \mapsto L(-1)v + h_{(-1)}w + \frac12(a_{(0)}b - b_{(0)}a)$,
\item the kernel of $S^2U_2 \to S_4$, $a \cdot b \mapsto \frac12(a_{(-1)}b + b_{(-1)}a)$,
\item the kernel of $S_3 \oplus (\mathfrak g \ot S_3) \oplus S^2U_2 \to S_4$, $(u, h \ot v, a \cdot b) \mapsto L(-1)u + h_{(-1)}v + \frac12(a_{(-1)}b + b_{(-1)}a)$,
\end{enumerate}
together with the centrality of $[\omega] = [\omega_A] + 2[e]$. The kernels are $\mathfrak g$-modules. Their highest weight vectors suffice. The Zhu images are obtained from \eqref{eq:zhu-rules} and from $[\frac12(a_{(-1)}b + b_{(-1)}a)] = \frac12([a]*[b] + [b]*[a]) + [c_+(a,b)]$, $[\frac12(a_{(0)}b - b_{(0)}a)] = [a]*[b] - [b]*[a] - [c_-(a,b)]$, where $c_\pm(a,b) = \frac12(a_{(1)}b \pm b_{(1)}a) \in U_2$.

\begin{thm}[computation modulo a prime]
\label{thm:zhu-modp}
Let $p = p_3 = 8388571$, and use the residue map $(\mathcal O_{\K})_{\mathfrak p} \to \mathbb F_p$ at the chosen prime $\mathfrak p$ above $p$, with the images of $\sqrt{13}$ and $\sqrt{-3}$ fixed by the computational data. The matrices below have entries integral at $\mathfrak p$. Modulo $p$:
\begin{enumerate}
\item the three families have $458$, $70$ and $532$ highest weight relations, and $S_3$ has $95$ highest weight vectors,
\item the projected relations, $19981$ rows in the path algebra $\mathbb F_p\Gamma_{\mathrm Z}$, kill the two exceptional vertices $12$ and $13$, hence the $230$ arrows incident to them, eliminate $46$ further arrows linearly, and leave a rewriting system on $239$ arrows whose ambiguities all resolve after one round of completion,
\item the resulting algebra $D_p$ has a basis of $762$ irreducible paths ($12$ vertices, $239$ arrows, $506$ paths of length two and $5$ of length three). Its regular trace form is nondegenerate, so $D_p$ is semisimple. Its center has dimension $12$, and $D_p \cong \bigoplus_{z}\Mat_{n_z}(\mathbb F_p)$ with the block sizes of Remark \ref{rem:expected-corner},
\item for each block, the ranks of the vertex idempotents $e_t$ are the multiplicities of the constituents $T_t$ in the lowest weight space of the corresponding module in Table \ref{tab:lowest-spaces}, for all twelve modules.
\end{enumerate}
\end{thm}

\begin{proof}
All products $a_{(n)}b$, $a,b \in U_2$, $-1 \leq n \leq 3$, and $h_{(-1)}s$, $L(-1)s$, $s \in S_3$, are computed exactly over $\K$ and then reduced modulo $p$. The kernels, projections, elimination, completion, trace form and central idempotents are computed over $\mathbb F_p$. The calculation at $8388571$ uses $19981$ projected relations. The corresponding calculation at $4194277$ gives the same dimensions, blocks and vertex ranks. For the final rewriting rules, a separate verifier checks that they decrease the degree-lexicographic order with the alphabet reversed and resolves all $131698$ overlap ambiguities and all $151$ inclusion ambiguities, as required by \cite[Theorem 1.2]{Ber78}. The inclusion ambiguities are needed because the system has cubic as well as quadratic leading words. These checks concern the projected presentation.

For completeness, the algebraic consequences of the last checks require no further VOA computation. If $x$ lies in the Jacobson radical of $D_p$, then $xy$ is nilpotent for every $y \in D_p$. Thus left multiplication by $xy$ has trace zero. Nondegeneracy of the regular trace form forces $x = 0$, proving semisimplicity. The computed central element has twelve distinct eigenvalues in $\mathbb F_p$, so the center is split. Each simple factor is a matrix algebra over $\mathbb F_p$, since a finite division ring is a field. Its dimension determines its matrix size. The ranks of the vertex corners give the displayed multiplicities. These conclusions concern the finite-field algebra $D_p$. Their descent to a characteristic-zero upper presentation is not supplied by this calculation.
\end{proof}

\begin{remark}\normalfont
\label{rem:zhu-status}
Theorem \ref{thm:zhu-modp} computes the upper presentation at two primes. Its block sizes and vertex ranks agree with those predicted by Theorem \ref{thm:weight-one-classification} and Table \ref{tab:lowest-spaces}. Agreement at two primes is a consistency check. The kernel may grow after reduction modulo a prime, and a rank drop at both primes is possible. To deduce (Z) by descent requires a finite integral order with this special fibre and the surjection from its characteristic-zero algebra onto the Zhu corner. These steps are carried out in \cite[Theorems A.10 and A.13]{ArhIII}. The survival of all twelve blocks in the Zhu corner is proved separately in \cite[Theorem A.20]{ArhIII}.
\end{remark}

\begin{lemma}
\label{lem:lattice-modules-unitary}
The nine $V'$-modules $W_{e,a}$ are unitary $S$-modules. Consequently each restriction $W_{e,a}|_S$ is an orthogonal direct sum of simple $S$-modules, and the $S$-submodule generated by an irreducible $A(S)$-submodule of the lowest weight space of $W_{e,a}$ is simple.
\end{lemma}

\begin{proof}
By \cite[Theorem 4.12]{DoL14}, $V_{L+\lambda}$ is a unitary module of $(V_L, \phi_1)$ for the Fock space Hermitian form $H_1$. With $g = \exp(h_\chi(0))$ as in the proof of Theorem \ref{thm:unitarity}, $g^{1/2}$ acts on $V_{L+\lambda}$ by the positive numbers $e^{\frac12(h_\chi,\mu)}$ on $M(1) \ot e^\mu$, and $H_\chi(u,v) = H_1(g^{1/2}u,g^{1/2}v)$ is positive definite and invariant for $(V', \phi_\chi)$. Its restriction to $S$ is invariant since $\phi_\chi(S) = S$. The orthogonal complement of a submodule is a submodule, and the homogeneous components are finite-dimensional, which gives the decomposition. If $X$ is an irreducible $A(S)$-submodule of the lowest weight space and $N$ a proper submodule of the module $\langle X\rangle$ it generates, then $\langle X\rangle = N \oplus N^\perp$ and $X = (N \cap X) \oplus (N^\perp \cap X)$ in the lowest degree, so $N \cap X = 0$, $X \subset N^\perp$ and $N = 0$.
\end{proof}

\subsection{Completion of the rationality argument}
\label{subsec:what-remains}

Appendix~A of \cite{ArhIII} proves (Z) and (M) for $S$. Theorem~2.2(F1) of \cite{ArhIII} states (Z) and the resulting strong rationality. Parts (F2)--(F3) give (M), including the complete lattice branching and the lowest weight spaces of the twelve simple modules. The proof of (M) combines the restrictions of the lattice modules with Galois symmetry to show that all twelve generic blocks survive in the Zhu corner. The braided equivalence that completes the proof of $\mathrm{(H_J)}$ is established in \cite[Theorem 6.1]{ArhIII}.

\section{The associated lattice and the twin algebra}
\label{sec:genus}

In this section we describe the Haagerup algebras through the decomposition of a strongly rational vertex operator algebra into its Heisenberg commutant and its associated lattice \cite{Mas14,Hoh17,MoR24a}, and we use the notions of genus of Höhn \cite{Hoh03} (the bulk genus, given by the representation category and the central charge) and of Moriwaki \cite{Mor21} and Möller--Rayhaun \cite{MoR24a} (the hyperbolic genus). We call a \emph{Haagerup algebra} a positive strongly rational vertex operator algebra $U$ of central charge eight with a braided equivalence $\Rep(U) \cong Z(\cal H_3)$. By Proposition \ref{prop:chain} and Theorem \ref{thm:weight-one-classification}, after an automorphism of $V'$, $U \subset V'$ with $U_1 = \mathfrak g' \oplus \mathfrak g''$, of levels $(k',k'') = (1,39)$ in case (I) and $(3,13)$ in case (II), and the Cartan subalgebra $\mathfrak h = \Cbb h' \oplus \Cbb h''$ of $U_1$ is contained in the Cartan subalgebra of $V'_1$ spanned by the lattice. Here $h' \in \mathfrak g'$ and $h'' \in \mathfrak g''$ are coroots, so that $\la h',h'\ra = 2k'$, $\la h'',h''\ra = 2k''$ and $\la h',h''\ra = 0$. We identify $\mathfrak h$ with its dual through $\la\cdot,\cdot\ra$ and write a vector $\lambda \in \mathfrak h$ as the pair $(m',m'') = (\la h',\lambda\ra,\la h'',\lambda\ra)$. Thus $\lambda$ is the $\mathfrak h$-weight of a vector $v$ with $h'_{(0)}v = m'v$ and $h''_{(0)}v = m''v$, and
\[
|\lambda|^2 = \frac{m'^2}{2k'}+\frac{m''^2}{2k''}.
\]
We denote by $\Gamma(U) \subset \mathfrak h$ the set of $\mathfrak h$-weights of $U$, by
\[
\Omega_U \ov{\mathrm{def}}{=} \operatorname{Com}_U(M(\mathfrak h))
\]
the \emph{Heisenberg commutant} of $U$, and by $L$ the \emph{associated lattice}, defined by $\operatorname{Com}_U(\Omega_U) \cong V_L$ \cite[Definition 2.12]{MoR24a}. Throughout this section $L$ denotes this rank-two lattice. We continue to write $E_6 \oplus A_2$ for the ambient lattice. These objects do not depend, up to isomorphism, on the choice of the Cartan subalgebra \cite[Remark 2.13]{MoR24a}. For an even lattice $L$ we write $\cal C(L'/L)$ for the pointed ribbon category of the discriminant form $(L'/L,q_L)$, $q_L(\lambda+L) = \frac{|\lambda|^2}2 \bmod 1$, so that $\Rep(V_L) \cong \cal C(L'/L)$, and $\overline{\cal B}$ for the ribbon reverse of a ribbon category $\cal B$.

\subsection{Two binary lattices}

Put $\Gamma_0 = \{(m',m'') \in \Z^2 : m' \equiv m'' \bmod 2\}$ and, for $(k',k'') \in \{(1,39),(3,13)\}$,
\[
L_{k',k''} \ov{\mathrm{def}}{=} \Gamma_0^* = \{(k'a,k''c) : a, c \in \Z,\ a \equiv c \bmod 2\},
\]
the dual being taken for the form above. It has the basis $(2k',0)$, $(k',k'')$. We write $L_{\mathrm I} = L_{1,39}$ and $L_{\mathrm{II}} = L_{3,13}$.

\begin{lemma}
\label{lem:two-lattices}
\begin{enumerate}
\item The Gram matrices of $L_{\mathrm I}$ and $L_{\mathrm{II}}$ in these bases are $\left[\begin{smallmatrix}2&1\\1&20\end{smallmatrix}\right]$ and $\left[\begin{smallmatrix}6&3\\3&8\end{smallmatrix}\right]$, of determinant $39$ and minimal norms $2$ and $6$. They correspond to the binary quadratic forms $x^2+xy+10y^2$ and $3x^2+3xy+4y^2$. In both cases $L' = \Gamma_0$ and $L'/L \cong \Z/39$ is generated by $(1,1)$, with $q_{L_{\mathrm I}}(1,1) = \frac{10}{39}$ and $q_{L_{\mathrm{II}}}(1,1) = \frac4{39}$.
\item The even positive-definite lattices of rank two whose discriminant form is isometric to that of $L_{\mathrm I}$ are, up to isometry, exactly $L_{\mathrm I}$ and $L_{\mathrm{II}}$. More precisely, the reduced forms of discriminant $-39$ are $(1,1,10)$, $(2,\pm1,5)$ and $(3,3,4)$. The class group of discriminant $-39$ is cyclic of order four and has two genera, the principal genus $\{(1,1,10),(3,3,4)\}$, which is the set of squares, and $\{(2,1,5),(2,-1,5)\}$, whose discriminant form is not isometric to that of $L_{\mathrm I}$.
\item The isometries $L_{\mathrm{II}}'/L_{\mathrm{II}} \to L_{\mathrm I}'/L_{\mathrm I}$ are the maps $(1,1) \mapsto u(1,1)$ with $u \in \{4,17,22,35\}$, and $O(L'/L) = \{\pm1,\pm14\}$ in both cases.
\item The image of the orthogonal group $O(L_{\mathrm{II}})$, of order four, in $O(L_{\mathrm{II}}'/L_{\mathrm{II}})$ is the whole group. The reflections $w' = \operatorname{diag}(-1,1)$ and $w'' = \operatorname{diag}(1,-1)$, which are the actions on $\mathfrak h$ of the Weyl group elements of $\mathfrak g'$ and $\mathfrak g''$, act by $14$ and $25$. The image of $O(L_{\mathrm I})$, of order four, in $O(L_{\mathrm I}'/L_{\mathrm I})$ is $\{\pm1\}$: $w'$ acts trivially and $w''$ by $-1$.
\end{enumerate}
\end{lemma}

\begin{proof}
We give a proof with integer arithmetic, independently of the computational verification. For part 1, the form $\frac{m'^2}{2k'}+\frac{m''^2}{2k''}$ on the basis $(2k',0)$, $(k',k'')$ has Gram matrix $M_{k',k''}$ with entries $M_{11} = 2k'$, $M_{12} = M_{21} = k'$ and $M_{22} = \frac{k'+k''}{2}$. Its determinant is $k'k'' = 39$. An integer multiple $n(1,1)$ belongs to $L_{k',k''}$ exactly when $39$ divides $n$, so $(1,1)$ generates the discriminant group. Its quadratic value is $\frac1{4k'}+\frac1{4k''}$, giving the asserted values. The expressions
\[
x^2+xy+10y^2 = \left(x+\frac y2\right)^2+\frac{39}{4}y^2,
\quad
3x^2+3xy+4y^2 = 3\left(x+\frac y2\right)^2+\frac{13}{4}y^2
\]
show that the shortest vectors are precisely the two vectors $\pm(2k',0)$, of squared norms $2$ and $6$, respectively.

2. An even positive-definite binary lattice of determinant $39$ corresponds to a primitive form $ax^2+bxy+cy^2$ of discriminant $-39$. Choose a reduced representative, so $|b| \leq a \leq c$, with the usual boundary convention $b \geq 0$ when $|b| = a$ or $a = c$. Since $39 = 4ac-b^2 \geq 3a^2$, we have $1 \leq a \leq 3$. Moreover $b$ is odd and $c = \frac{39+b^2}{4a}$ is an integer. For $a = 1,2,3$ this gives exactly $(1,1,10)$, $(2,\pm1,5)$ and $(3,3,4)$. The two middle forms are improperly equivalent. Their lattice has Gram matrix $\left[\begin{smallmatrix}4&1\\1&10\end{smallmatrix}\right]$ and discriminant form generated by an element of quadratic value $\frac5{39}$. An isometry with the discriminant form of $L_{\mathrm I}$ would require $5u^2 \equiv 10 \pmod {39}$ for a unit $u$, which is impossible modulo $3$. In contrast, $10\cdot4^2 \equiv 4 \pmod {39}$ identifies the discriminant forms of $L_{\mathrm I}$ and $L_{\mathrm{II}}$. This proves the lattice classification. It also gives the two genera, since the genus is determined by the rank, signature and discriminant form \cite[Corollary 1.9.4]{Nik80}. The four reduced forms are the four proper classes. Inversion interchanges the two distinct classes $(2,1,5)$ and $(2,-1,5)$, so the class group of order four is cyclic. Its squares are the principal class and the class $(3,3,4)$, proving the remaining assertion \cite{Cox13}.

3. An isometry between the discriminant forms is multiplication by a unit $u$ satisfying $10u^2 \equiv 4 \pmod {39}$, equivalently $u^2 \equiv 16 \pmod {39}$. The Chinese remainder theorem gives the four solutions $4,17,22,35$. In the same way, $u^2 \equiv 1 \pmod {39}$ gives $O(L'/L) = \{1,14,25,38\}$.

4. An orthogonal transformation preserving either lattice preserves its unique shortest line, so its orthogonal group has at most four elements. The four sign changes in the two coordinates $(m',m'')$ preserve $L_{k',k''}$, so they give the whole group. To compute their action on the discriminant group, reduce $(-1,1)$ and $(1,-1)$ modulo $L_{k',k''}$. In case (II),
\[
(-1,1)-14(1,1) = (-15,-13) \in L_{\mathrm{II}},
\quad
(1,-1)-25(1,1) = (-24,-26) \in L_{\mathrm{II}}.
\]
In case (I), $(-1,1)-(1,1) = (-2,0) \in L_{\mathrm I}$ and $(1,-1)+(1,1) = (2,0) \in L_{\mathrm I}$. This gives the stated images.

The supplementary material includes independent integer-arithmetic checks of these lattice calculations. They are not needed for the preceding proof.
\end{proof}

\subsection{The Heisenberg commutant and the associated lattice}

For a simple vertex operator algebra $\Omega$ and an $\Omega$-module $M$, put
\[
\operatorname{Vac}_{\Omega}(M) \ov{\mathrm{def}}{=} \{v \in M : a_{(n)}v = 0 \text{ for every } a \in \Omega \text{ and } n \geq 0\}.
\]
A nonzero vector in this space is called \emph{vacuum-like}. If $M$ is irreducible and $v$ is vacuum-like, the map $a \mapsto a_{(-1)}v$ is a nonzero $\Omega$-module homomorphism $\Omega \to M$, and hence $M \cong \Omega$ \cite[Section 4.7]{LL04}.

\begin{prop}
\label{prop:associated-lattice}
Let $U$ be a Haagerup algebra, with the notation above.
\begin{enumerate}
\item The Heisenberg commutant $\Omega_U$ is a strongly rational vertex operator algebra of central charge six with $(\Omega_U)_1 = 0$. There is a ribbon-reversing equivalence $\tau$ from $\cal C(L'/L)$ onto the full pointed fusion subcategory generated by the invertible objects of $\Rep(\Omega_U)$. Writing $\tau$ also for the induced map on isomorphism classes of simple objects, we have
\[
U = \bigoplus_{\alpha \in L'/L} \Omega_U^{\tau(\alpha)} \ot V_{\alpha+L}
\]
as an $\Omega_U \ot V_L$-module, and $\Rep(\Omega_U) \cong Z(\cal H_3) \boxtimes \overline{\cal C(L'/L)}$.
\item Every irreducible $\Omega_U$-module has a nonnegative conformal weight, and for $\alpha \neq 0$ the module $\Omega_U^{\tau(\alpha)}$ contains no vacuum-like vector. Consequently, a vector of $U$ of $\mathfrak h$-weight $\lambda$ has conformal weight at least $\frac{|\lambda|^2}2$.
\item $\Gamma(U) = \Gamma_0$, and the associated lattice is $L = L_{k',k''}$: it is $L_{\mathrm I}$ in case (I) and $L_{\mathrm{II}}$ in case (II).
\item For $\lambda = (k',k'')$, the vector $e^\lambda$ of $V_L \subset U$ is a highest weight vector of weight $(k',k'')$ for the affine subalgebra $A = L_{\mathfrak g'}(k',0) \ot L_{\mathfrak g''}(k'',0)$, of conformal weight $\frac{k'+k''}4$, and it generates an $A$-submodule of $U$ isomorphic to the simple current $J = L_{\mathfrak{sl}_2}(k',k') \ot L_{\mathfrak{sl}_2}(k'',k'')$, of conformal weight $10$ in case (I) and $4$ in case (II).
\end{enumerate}
\end{prop}

\begin{proof}
1. By \cite[Definition 2.12, Proposition 2.14]{MoR24a}, which rest on \cite{Mas14}, $\Omega_U$ and $V_L$ form a dual pair in $U$, $\Omega_U$ is strongly rational \cite{CKLR19}, and $U$ is a simple current extension of $\Omega_U \ot V_L$ of the stated form for some subgroup $A_U \subset L'/L$ in place of $L'/L$. The category $\Rep(U) \cong Z(\cal H_3)$ has no invertible object except the unit, by \eqref{eq:center-dimensions}. It is pseudo-unitary for the ribbon structure of $U$: since $U$ is positive, the quantum dimensions $\operatorname{qdim}_U(X) = S_{X,U}/S_{U,U}$ are positive and equal the Frobenius--Perron dimensions \cite{DJX13}, and they are the categorical dimensions of $\Rep(U)$ for its ribbon structure \cite{Hua08}. Alternatively, $Z(\cal H_3)$ has a unique pivotal structure (Remark \ref{rem:essential-positivity}). The proof of \cite[Corollary 4.10]{MoR24a}, applied with $V = W = U$, then gives $A_U = L'/L$, the description of $\tau$, and $\Rep(\Omega_U) \cong \Rep(U) \boxtimes \overline{\cal C(L'/L)}$. The central charge of $\Omega_U$ is $8-\operatorname{rank}L = 6$. A vector of $(\Omega_U)_1$ commutes with $\mathfrak h$, hence lies in the subspace of $\mathfrak h$-weight zero of $U_1 = \mathfrak{sl}_2 \oplus \mathfrak{sl}_2$, which is $\mathfrak h$. It is then orthogonal to $\mathfrak h$, hence zero.

2. Let $\mu$ be an irreducible $\Omega_U$-module. The monodromy of $\mu$ with the invertible objects $\tau(\alpha)$ is a scalar $\chi_\mu(\alpha)$, and $\alpha \mapsto \chi_\mu(\alpha)$ is a character of $L'/L$. Since the bilinear form of $L'/L$ is nondegenerate, there is $\beta \in L'/L$ such that the monodromy of $\mu \boxtimes \beta$ with every summand $\tau(\alpha) \boxtimes \alpha$ of the algebra $U$ is trivial. The induced module $X = \bigoplus_\alpha (\tau(\alpha)\mu) \ot V_{\alpha+\beta+L}$, where $\tau(\alpha)\mu$ denotes the fusion product, is then a $U$-module \cite{HKL15,CKM22}, a finite sum of irreducible $U$-modules, and it contains $\mu \ot V_{\beta+L}$. For $y > 0$ we have, coefficientwise and hence numerically,
\[
\operatorname{ch}_X(\i y) \geq \operatorname{ch}_\mu(\i y)\frac{\theta_{\beta+L}(\i y)}{\eta(\i y)^2}.
\]
As $y \to 0^+$, $\theta_{\beta+L}(\i y) \sim (\det L)^{-\frac12}y^{-1}$ by the Poisson summation formula and $\eta(\i y)^{-2} = y\eta(\i/y)^{-2} \sim y\e^{\frac{2\pi}{12y}}$, so that the last quotient is equivalent to $(\det L)^{-\frac12}\e^{\frac{2\pi\cdot2}{24y}}$. On the other hand $\operatorname{ch}_X(\i y) = \sum_Y S_{X,Y}\operatorname{ch}_Y(\i/y)$ over the irreducible $U$-modules $Y$, by Zhu's theorem \cite{Zhu96} at the point $\i/y$, extended additively to the finite sum $X$, and $\operatorname{ch}_Y(\i/y) = O(\e^{\frac{2\pi\cdot8}{24y}})$ since the conformal weights of the $U$-modules are nonnegative. Hence $\operatorname{ch}_\mu(\i y) = O(\e^{\frac{2\pi\cdot6}{24y}})$ for every irreducible $\Omega_U$-module $\mu$. For an irreducible $\Omega_U$-module $\nu$ of conformal weight $h_\nu$, Zhu's theorem for $\Omega_U$ gives $\operatorname{ch}_\nu(\i/y) = \sum_\mu S^{\Omega}_{\nu,\mu}\operatorname{ch}_\mu(\i y) = O(\e^{\frac{2\pi\cdot6}{24y}})$, while $\operatorname{ch}_\nu(\i/y) \geq \dim(\nu_{h_\nu})\e^{\frac{2\pi}{y}(\frac6{24}-h_\nu)}$. Therefore $h_\nu \geq 0$. For $\alpha \neq 0$, the simple module $\Omega_U^{\tau(\alpha)}$ is not isomorphic to $\Omega_U$, so $\operatorname{Vac}_{\Omega_U}(\Omega_U^{\tau(\alpha)}) = 0$. Finally a vector of $\mathfrak h$-weight $\lambda \in \alpha+L$ lies in $\Omega_U^{\tau(\alpha)} \ot M(\mathfrak h,\lambda)$, whose conformal weights are at least $h_{\tau(\alpha)}+\frac{|\lambda|^2}2 \geq \frac{|\lambda|^2}2$.

3. The $\mathfrak h$-weights of $\Omega_U^{\tau(\alpha)} \ot V_{\alpha+L}$ form the coset $\alpha+L$, so part 1 gives $\Gamma(U) = L'$ and $L = \Gamma(U)^*$. We have $V'_1 = U_1 \oplus R_b$. By Table \ref{tab:lowest-spaces}, every irreducible $U_1$-summand of $V'_1$ has highest weight $(a,b)$ with $a \equiv b \bmod 2$. All weights of such a summand belong to $\Gamma_0$. Since $V'$ is generated by $V'_1$, every $\mathfrak h$-weight of $V'$ belongs to $\Gamma_0$, and hence $\Gamma(U) \subset \Gamma_0$. Conversely, the root vectors of $U_1$ give the weights $(2,0)$ and $(0,2)$, while the summand $V(1,1) \subset U_2$ in both columns of Table \ref{tab:lowest-spaces} gives $(1,1)$. These vectors generate $\Gamma_0$. Since $\Gamma(U) = L'$ is a group, $\Gamma(U) = \Gamma_0$ and $L = L_{k',k''}$.

4. Since $\lambda = (k',k'') \in L$, $e^\lambda \in V_L \subset U$, of conformal weight $\frac{|\lambda|^2}2 = \frac{k'+k''}4$. Let $x$ be a root vector of $U_1$, of weight $\mu \in \{(\pm2,0),(0,\pm2)\}$, or $x \in \mathfrak h$, $\mu = 0$, and $n \geq 0$. The vector $x_{(n)}e^\lambda$ has $\mathfrak h$-weight $\lambda+\mu$ and conformal weight $\frac{k'+k''}4-n$. For $n = 0$ and $\mu \in \{(2,0),(0,2)\}$, $\frac{|\lambda+\mu|^2}2-\frac{k'+k''}4$ is $1+\frac1{k'}$ or $1+\frac1{k''}$. For $n \geq 1$, $\frac{|\lambda+\mu|^2}2-(\frac{k'+k''}4-n) \geq \min(0,-1+\frac1{k'},-1+\frac1{k''})+n > 0$. By part 2, $x_{(n)}e^\lambda = 0$ in all these cases. Thus $e^\lambda$ is annihilated by the positive part of the affine Lie algebra and by $e'_{(0)}$, $e''_{(0)}$, and $h'_{(0)}$, $h''_{(0)}$ act by $k'$, $k''$: it is a highest weight vector of weight $(k',k'')$ for $A$. The regularity of $A$ decomposes $U$ into irreducible integrable $A$-modules. In this decomposition $e^\lambda$ is a sum of affine highest weight vectors of type $(k',k'')$, hence generates one diagonal copy of the irreducible module $J$, whose conformal weight is $\frac{k'(k'+2)}{4(k'+2)}+\frac{k''(k''+2)}{4(k''+2)} = \frac{k'+k''}4$.
\end{proof}

In case (II), part 4 is Proposition \ref{prop:simple-current}, proved there for $S$ by an explicit lattice word and an energy bound, with exact computational checks. In case (I), it is the inclusion $A \oplus J(1,39) \subset U^{\mathrm R}$ of \cite[Proposition 7.1]{Rie26b}, obtained there from the explicit vector $e^w$, $(w,w) = 20$. Here both statements hold for every Haagerup algebra. The associated lattice is determined by the weight-one decomposition of $V'$ and the occurrence of $V(1,1)$ in $U_2$. The lattices $L_{\mathrm I}$ and $L_{\mathrm{II}}$ are thus the overlattices of index two of $\sqrt{2k'}\Z \oplus \sqrt{2k''}\Z$ given by the simple current $J$.

\subsection{The lattice exchange}

Let $U$ be a Haagerup algebra in case (I) or (II), with associated lattice $L$, and let $L^\natural$ be the lattice of the other case, so that $\{L,L^\natural\} = \{L_{\mathrm I},L_{\mathrm{II}}\}$. By Lemma \ref{lem:two-lattices} there is an isometry $\psi \co L'/L \to (L^\natural)'/L^\natural$. We use the equivalences $\Rep(\Omega_U \ot V_L) \cong \Rep(\Omega_U) \boxtimes \cal C(L'/L)$ and $\Rep(\Omega_U \ot V_{L^\natural}) \cong \Rep(\Omega_U) \boxtimes \cal C((L^\natural)'/L^\natural)$ \cite{FHL93,CKM22}. The isometry $\psi$ is induced by a ribbon equivalence $\Phi_\psi \co \cal C(L'/L) \to \cal C((L^\natural)'/L^\natural)$, since a pointed braided fusion category is determined by its metric group up to braided equivalence inducing any given isometry \cite[Section 8.4]{EGNO15}, and the twists $\e^{2\pi\i q}$ are preserved. Put $\Phi = \mathrm{id} \boxtimes \Phi_\psi$. By \cite{HKL15}, $U$ corresponds to the connected étale algebra $\cal A_U = \bigoplus_\alpha \tau(\alpha) \boxtimes \alpha$ of $\Rep(\Omega_U) \boxtimes \cal C(L'/L)$, and we define $U^\natural = U^\natural_\psi$ as the vertex algebra extension of $\Omega_U \ot V_{L^\natural}$ corresponding to the connected étale algebra $\Phi(\cal A_U)$ \cite{HKL15,CKM22}. As an $\Omega_U \ot V_{L^\natural}$-module,
\[
U^\natural = \bigoplus_{\alpha \in L'/L} \Omega_U^{\tau(\alpha)} \ot V_{\psi(\alpha)+L^\natural}.
\]

\begin{thm}
\label{thm:twin}
Let $U$ be a Haagerup algebra, $L^\natural$ and $\psi$ as above.
\begin{enumerate}
\item $U^\natural$ is a Haagerup algebra, in the other case of Theorem \ref{thm:weight-one-classification}. For the Cartan subalgebra $L^\natural \ot \Cbb$ of $U^\natural_1$, its Heisenberg commutant is $\Omega_U$ and its associated lattice is $L^\natural$.
\item Let $P$ be the group of invertible objects of $\Rep(\Omega_U)$, with the quadratic form given by the twists. The image of $\operatorname{Aut}(\Omega_U)$ in $O(P)$ is the whole group $O(P)$. Consequently, for a strongly rational vertex operator algebra $\Omega$ and $N \in \{L_{\mathrm I},L_{\mathrm{II}}\}$, all the Haagerup algebras with Heisenberg commutant isomorphic to $\Omega$ and associated lattice isometric to $N$ are isomorphic.
\item Up to isomorphism, $U^\natural$ does not depend on $\psi$, $(U^\natural)^\natural \cong U$, and $U \mapsto U^\natural$ induces a bijection between the isomorphism classes of Haagerup algebras in case (I) and in case (II). The Haagerup algebras with Heisenberg commutant isomorphic to $\Omega_U$ are, up to isomorphism, exactly $U$ and $U^\natural$.
\item $U$ and $U^\natural$ belong to the same hyperbolic genus, and every strongly rational vertex operator algebra in this hyperbolic genus is isomorphic to $U$ or to $U^\natural$. In particular it is a Haagerup algebra.
\end{enumerate}
\end{thm}

\begin{proof}
1. \emph{Extension and representation category.} Put $B = \Omega_U \ot V_{L^\natural}$ and $\cal A = \Phi(\cal A_U)$. The strongly rational algebra $B$ has its vertex tensor category of ordinary modules. We use the extension correspondence in the generality of \cite[Theorem 3.2 and Remark 3.3]{HKL15}. Strict positivity of the nonvacuum $B$-modules is not assumed. The connected étale algebra $\cal A$ has trivial twist and defines the extension $W$. It is simple: separability makes the category of $\cal A$-modules semisimple, and $\operatorname{End}_{\cal A}(\cal A) \cong \operatorname{Hom}(\mathbf 1,\cal A) = \Cbb$ by connectedness. Thus the regular $\cal A$-module has no proper subobject, and neither does $W$ as a module over itself. Initially, the extension correspondence identifies local $\cal A$-modules with those $W$-modules whose restrictions belong to $\Rep(B)$ \cite{HKL15,CKM22}. Their category is ribbon equivalent, through $\Phi$, to the category of local $\cal A_U$-modules. After the rationality argument below, this is the ordinary module category $\Rep(W) \cong \Rep(U)$.

\emph{CFT type.} For $\alpha \neq 0$, the conformal weights of $\Omega_U^{\tau(\alpha)} \ot V_{\psi(\alpha)+L^\natural}$ lie in $h_{\tau(\alpha)}+q_{L^\natural}(\psi(\alpha))+\Z = h_{\tau(\alpha)}+q_L(\alpha)+\Z = \Z$, and they are at least $h_{\tau(\alpha)}+\min_{\gamma \in \psi(\alpha)+L^\natural}\frac{|\gamma|^2}2 > 0$ by Proposition \ref{prop:associated-lattice}, part 2. Hence they are at least $1$. The summand $\alpha = 0$ is $\Omega_U \ot V_{L^\natural}$, which is of CFT type. Thus $W$ is of CFT type.

\emph{Self-contragredience.} Since $(\Omega_U)_1 = 0$, $W_1 = \mathbf 1 \ot (V_{L^\natural})_1 \oplus \bigoplus_{\alpha \neq 0}(\Omega_U^{\tau(\alpha)} \ot V_{\psi(\alpha)+L^\natural})_1$, and the Virasoro modes preserve each summand. The operator $L(1)$ maps $W_1$ into $W_0 = \Cbb\mathbf 1$, which lies in the summand $\alpha = 0$, while $L(1)(V_{L^\natural})_1 = 0$. Hence $L(1)W_1 = 0$ and the simple algebra $W$ of CFT type is self-contragredient \cite{Li94}.

\emph{$C_2$-cofiniteness.} The algebra $W$ is $C_2$-cofinite, since it is a finite sum of irreducible modules for the $C_2$-cofinite algebra $\Omega_U \ot V_{L^\natural}$, which are $C_2$-cofinite \cite{Buh02}, and $C_2$ of each summand is contained in $C_2(W)$.

\emph{Rationality.} Let $M$ be a weak $W$-module. By \cite{ABD04}, $M$ is a direct sum of irreducible ordinary $\Omega_U \ot V_{L^\natural}$-modules. Put $B = \Omega_U \ot V_{L^\natural}$, and let $M'$ be the $W$-submodule generated by an irreducible $B$-submodule $N$ of $M$. By weak associativity, $M'$ is spanned by the vectors $w_{(n)}x$, with $w \in W$, $x \in N$ and $n \in \Z$. For every irreducible $B$-module $K$, composition with these vertex operators gives an injective linear map
\[
\operatorname{Hom}_B(M',K) \to I_B\binom{K}{W\ N}.
\]
It is injective because a map in its kernel vanishes on this spanning set. The target is finite dimensional and is zero for all but finitely many $K$. In a decomposition of the semisimple $B$-module $M'$ into simple summands, the coordinate projections onto distinct copies of $K$ are linearly independent elements of $\operatorname{Hom}_B(M',K)$. Hence each multiplicity is finite, and only finitely many simple types occur. Thus $M'$ has finite length and is an object of $\Rep(B)$. By \cite{HKL15,CKM22}, it is a local module for the connected étale algebra $\Phi(\cal A_U)$, and such modules form a semisimple category \cite{DMNO13}. Hence this submodule is a finite direct sum of irreducible $W$-modules, and $M$, a sum of such submodules, is a direct sum of irreducible $W$-modules. Compare the proof of \cite[Proposition 5.3]{MoR24a}, which is the same construction in the holomorphic case.

\emph{Positivity.} Let $X$ be an irreducible $W$-module. Its restriction is a sum of modules $\Omega_U^\nu \ot V_{\gamma+L^\natural}$, whose conformal weights are at least $h_\nu \geq 0$ (Proposition \ref{prop:associated-lattice}, part 2), so that $h_X \geq 0$. If $h_X = 0$, some summand has $h_\nu = 0$ and $\gamma \in L^\natural$, so that $X$ contains $\Omega_U^\nu \ot V_{L^\natural}$ with $h_\nu = 0$. The local $\cal A_U$-module $X^\flat = \Phi^{-1}(X)$ is an irreducible $U$-module whose restriction contains $\Phi^{-1}(\nu \boxtimes 0) = \nu \boxtimes 0$, that is, $\Omega_U^\nu \ot V_L$. Hence $h_{X^\flat} \leq h_\nu = 0$, and $X^\flat \cong U$ since $U$ is positive. Thus $X \cong \Phi(\cal A_U) = W$. Hence $W$ is a Haagerup algebra, and by Theorem \ref{thm:weight-one-classification} it is in case (I) or (II).

\emph{Heisenberg commutant.} Put $\mathfrak h^\natural = \mathbf 1 \ot (L^\natural \ot \Cbb)_{(-1)}\mathbf 1 \subset W_1$. The vectors of $W_1$ of $\mathfrak h^\natural$-weight zero lie in $(\Omega_U)_1 \ot \mathbf 1+\mathbf 1 \ot \mathfrak h^\natural = \mathfrak h^\natural$, since $(\Omega_U)_1 = 0$. Thus $\mathfrak h^\natural$ is a self-centralizing subalgebra of the reductive Lie algebra $W_1$ \cite{DoM04} acting semisimply on $W$, that is, a Cartan subalgebra. The commutant of $M(\mathfrak h^\natural)$ in $V_{\gamma+L^\natural}$ is $\Cbb e^0$ if $\gamma \in L^\natural$ and zero otherwise, so $\operatorname{Com}_W(M(\mathfrak h^\natural)) = \Omega_U \ot \mathbf 1$, and $\operatorname{Com}_W(\Omega_U) = \bigoplus_\alpha \operatorname{Vac}_{\Omega_U}(\Omega_U^{\tau(\alpha)}) \ot V_{\psi(\alpha)+L^\natural} = V_{L^\natural}$ by Proposition \ref{prop:associated-lattice}, part 2. So the Heisenberg commutant of $W$ is $\Omega_U$ and its associated lattice is $L^\natural$. By Proposition \ref{prop:associated-lattice}, part 3, applied to $W$, its associated lattice is $L_{\mathrm I}$ or $L_{\mathrm{II}}$ according to its case. Since $L^\natural$ is not isometric to $L$ (Lemma \ref{lem:two-lattices}), $W$ is in the other case.

2. The map $\tau \co L'/L \to P$ is an anti-isometry. Let $U_{\mathrm{II}}$ be the one of $U$, $U^\natural$ which is in case (II), $L_{\mathrm{II}}$ its associated lattice and $\tau_{\mathrm{II}} \co L_{\mathrm{II}}'/L_{\mathrm{II}} \to P$ its gluing map (part 1). The Weyl group elements of the two ideals of $(U_{\mathrm{II}})_1$ are inner automorphisms of $U_{\mathrm{II}}$ which normalize the Cartan subalgebra. They preserve $\Omega_U = \operatorname{Com}(M(\mathfrak h))$ and $V_{L_{\mathrm{II}}} = \operatorname{Com}(\Omega_U)$, and restrict to $\theta \in \operatorname{Aut}(\Omega_U)$ and to an automorphism of $V_{L_{\mathrm{II}}}$ inducing $w'$ or $w''$ on $L_{\mathrm{II}}'/L_{\mathrm{II}}$. Since they permute the summands $\Omega_U^{\tau_{\mathrm{II}}(\alpha)} \ot V_{\alpha+L_{\mathrm{II}}}$, the induced autoequivalence $\bar\theta$ satisfies $\bar\theta\tau_{\mathrm{II}} = \tau_{\mathrm{II}}\bar w$. By Lemma \ref{lem:two-lattices}, part 4, $\bar w' = 14$ and $\bar w'' = 25$, which generate $O(L_{\mathrm{II}}'/L_{\mathrm{II}})$. Hence the $\bar\theta|_P$ exhaust $O(P) = \tau_{\mathrm{II}}O(L_{\mathrm{II}}'/L_{\mathrm{II}})\tau_{\mathrm{II}}^{-1}$. Now let $Y_1$, $Y_2$ be Haagerup algebras with Heisenberg commutants identified with $\Omega$ and associated lattices identified with $N$, with gluing maps $\sigma_1,\sigma_2 \co N'/N \to P(\Omega)$. Since both are anti-isometries, $\sigma_2 = \bar\theta\sigma_1$ for some $\theta \in \operatorname{Aut}(\Omega)$, by the first assertion applied to $Y_1$. The automorphism $\theta \ot \mathrm{id}$ of $\Omega \ot V_N$ transforms $Y_1$ into a simple current extension with the summands $\Omega^{\sigma_2(\alpha)} \ot V_{\alpha+N}$ of $Y_2$. The summands are pairwise non-isomorphic invertible $\Omega \ot V_N$-modules, indexed by $N'/N$, and their fusion satisfies
\[
M_\alpha \boxtimes M_\beta \cong M_{\alpha+\beta}.
\]
All summands are nonzero, both extensions are simple, and the fusion multiplicities are one. Thus \cite[Proposition 5.3]{DoM04} applies to their common graded module decomposition. Hence $Y_1 \cong Y_2$.

3. The algebra $U^\natural_\psi$ has Heisenberg commutant $\Omega_U$ and associated lattice $L^\natural$, so it does not depend on $\psi$ by part 2. Applying the construction to $U^\natural$ with $\psi^{-1}$ gives back the summands of $U$, so $(U^\natural)^\natural \cong U$ by \cite[Proposition 5.3]{DoM04}. If $Y \cong Y'$ are Haagerup algebras, an isomorphism maps a Cartan subalgebra to a Cartan subalgebra and identifies the Heisenberg commutants, so $Y^\natural \cong Y'^\natural$ by part 2. Hence $\natural$ is well defined on isomorphism classes and is an involution exchanging the two cases. If $W$ is a Haagerup algebra with $\Omega_W \cong \Omega_U$, its associated lattice is $L$ or $L^\natural$ according to its case, and $W \cong U$ or $W \cong U^\natural$ by part 2.

4. By \cite[Corollary 4.10]{MoR24a}, two strongly rational vertex operator algebras with equivalent unpointed pseudo-unitary representation categories belong to the same hyperbolic genus if and only if they have the same central charge and isomorphic Heisenberg commutants. This applies to $U$ and $U^\natural$. If $V$ is strongly rational and in this hyperbolic genus, then by \cite[Theorems 4.7 and 4.8]{MoR24a} $\Rep(V)$ is ribbon equivalent to $\Rep(U)$, the Heisenberg commutant of $V$ is isomorphic to $\Omega_U$, and its associated lattice is in the genus of $L$, hence isometric to $L_{\mathrm I}$ or $L_{\mathrm{II}}$ (Lemma \ref{lem:two-lattices}). By the proof of \cite[Corollary 4.10]{MoR24a}, $V = \bigoplus_\alpha \Omega_U^{\sigma(\alpha)} \ot V_{\alpha+N}$ with $N$ its associated lattice and $\sigma$ an anti-isometry onto $P$. The argument of part 2, with $Y_1$ the one of $U$, $U^\natural$ whose associated lattice is isometric to $N$ and $Y_2 = V$, uses only the decomposition of $Y_2$ and shows that $V \cong Y_1$.
\end{proof}

\begin{thm}
\label{thm:classification-one}
\begin{enumerate}
\item Assuming $\mathrm{(H_J)}$, there are exactly two Haagerup algebras up to isomorphism: Jiao's algebra $V \cong S$, in case (II), and its twin $S^\natural$, in case (I). There is an embedding of $S^\natural$ into $V'$ whose image contains the algebra $\psi_{\mathrm R}(U^{\mathrm R})$ generated by $A$, $e$, $f$, $v$ (Proposition \ref{prop:riesen-identification}), with equality in weights at most four.
\item Assuming $\mathrm{(H_R)}$, there are exactly two Haagerup algebras up to isomorphism: Riesen's algebra $U^{\mathrm R}$, in case (I), and its twin $(U^{\mathrm R})^\natural$, in case (II). There is an embedding of $(U^{\mathrm R})^\natural$ into $V'$ whose image contains $S$, with equality in weights at most four.
\item Assuming $\mathrm{(H_R)}$ and $\mathrm{(H_J)}$, $S^\natural \cong U^{\mathrm R}$ and $(U^{\mathrm R})^\natural \cong V$.
\item Assuming $\mathrm{(H_R)}$ or $\mathrm{(H_J)}$, all Haagerup algebras have isomorphic Heisenberg commutants and they form a single hyperbolic genus.
\end{enumerate}
\end{thm}

\begin{proof}
1. By Theorem \ref{thm:uniqueness-two}, every Haagerup algebra in case (II) is isomorphic to $V \cong S$. By Theorem \ref{thm:twin}, part 3, every Haagerup algebra in case (I) is isomorphic to $S^\natural$, which is one. The proof of Theorem \ref{thm:uniqueness-one} shows, without using \cite{Rie26}, that a Haagerup algebra in case (I) contains $\psi_{\mathrm R}(U^{\mathrm R})$ after an automorphism of $V'$. By Proposition \ref{prop:riesen-dims}, the weight spaces of $\psi_{\mathrm R}(U^{\mathrm R})$ of weights two, three and four have dimensions at least $120$, $956$ and $6010$, which are those of the Haagerup vacuum character (Theorem \ref{thm:gamma-one}). 2. The same argument with Theorem \ref{thm:uniqueness-one}. The proof of Theorem \ref{thm:uniqueness-two} shows, without \cite{Jia26}, that a Haagerup algebra in case (II) contains $S$ after an automorphism of $V'$, and Theorem \ref{thm:strong-generation}, using the independent certificate of Proposition \ref{prop:weight-four-certificate}, gives the dimensions $1$, $6$, $120$, $956$, $6010$ of $S$ in weights at most four. 3 follows from 1 and 2, and 4 from Theorem \ref{thm:twin}, parts 3 and 4.
\end{proof}

Theorem \ref{thm:twin} pairs the two cases within each hyperbolic genus. For the explicit algebra $S$, \cite[Theorem 6.1]{ArhIII} supplies $\mathrm{(H_J)}$. Parts 1 and 4 of Theorem \ref{thm:classification-one} then give the complete classification by the pair $S,S^\natural$, including existence in case (I) and the common Heisenberg commutant.

The common fixed-point algebra of \cite[Theorem 7.3]{DJX26} is used in \cite[Theorem 7.9]{ArhIII} to identify $S^\natural$ with Huang's algebra $U^{\mathrm H}$. The identification of Huang's and Riesen's algebras in \cite[Section 4.5 and Proposition 4.1]{DJX26} then gives $S^\natural\cong U^{\mathrm H}\cong U^{\mathrm R}$ as vertex operator algebras, as stated in \cite[Corollary 7.10]{ArhIII}. Thus $\mathrm{(H_R)}$ holds, parts 2--4 of Theorem \ref{thm:classification-one} apply, and $(U^{\mathrm R})^\natural\cong S$ by Theorem \ref{thm:twin}(3). The ambient containment and the finite weight comparison in Theorem \ref{thm:classification-one}(1) come from the calculations in this article.

\begin{remark}\normalfont
\label{rem:genus-characters}
The decomposition of Proposition \ref{prop:associated-lattice} can be tested on the graded dimensions. For $\mu \in \Gamma_0$, the generating series of the dimensions of the subspaces of $\mathfrak h$-weight $\mu$ of $U_n$, multiplied by $\prod_{m \geq 1}(1-q^m)^2$ and by $q^{-|\mu|^2/2}$, is the graded dimension $f_\alpha$ of $\Omega_U^{\tau(\alpha)}$, $\alpha = \mu+L$. From the blocks of $U_2$ and $U_3$ (exact in case (II), from the spans of Proposition \ref{prop:riesen-dims} modulo $7999993$ in case (I)), a finite calculation shows that different weights $\mu$ in the same coset give the same coefficients (39 cosets occur in case (II) and 35 in case (I) in these weights), and that for each of the four isometries $\psi$ of Lemma \ref{lem:two-lattices}, part 3, the $88$ coefficients of $f_\alpha$ in case (II) which can be compared with those of $f_{\psi(\alpha)}$ in case (I) agree. This is consistent with Theorem \ref{thm:classification-one}, part 4: it is a test, in low weights, of the isomorphism of the Heisenberg commutants of any two Haagerup algebras of the two cases, which Theorem \ref{thm:twin} proves for $U$ and $U^\natural$.
\end{remark}

\begin{remark}\normalfont
\label{rem:essential-positivity}
Let $V$ be any strongly rational vertex operator algebra of central charge eight with a braided equivalence $\Rep(V) \cong Z(\cal H_3)$, not assumed positive. Since $Z(\cal H_3)$ has no nontrivial invertible object, its universal grading group is trivial \cite[Theorem 6.3]{GeN08}. Since pivotal structures differ by tensor automorphisms of the identity functor, which are characters of this group, $Z(\cal H_3)$ has a unique pivotal structure, the pseudo-unitary one \cite{ENO05}. Since $c(V) \equiv 0 \equiv c(Z(\cal H_3)) \bmod 8$, the $S$-matrix of the characters of $V$ is the normalized $S$-matrix $s$ of Section \ref{sec:preliminaries} \cite[Lemma 3.6]{MoR24a}. Suppose that the minimal conformal weight $h_{\min}$ of the irreducible $V$-modules is negative, and let $J$ be the set of simple objects with conformal weight $h_{\min}$. It does not contain $0$, and its elements have the same twist. The possible nonempty sets $J$ are the nonempty subsets of $\{b,a,c_0\}$, the singletons $\{c_1\}$ and $\{c_2\}$, and the singletons $\{d_l\}$. The following choices of $i$ give $s_{ij} \leq 0$ for every $j \in J$, with strict inequality for at least one $j$:
\[
\begin{array}{c|c}
\text{condition on } J & i \\ \hline
b \in J & d_1 \\
\varnothing \neq J \subset \{a,c_0\} & c_1 \\
J = \{c_1\} \text{ or } J = \{c_2\} & a \\
J = \{d_l\} & b
\end{array}
\]
The signs follow directly from the entries $-y/3$, $-1/3$ and $0$ of the $S$-matrix. Then $\operatorname{ch}_i(\i y) \sim \big(\sum_{j \in J}s_{ij}\dim V^j_{h_{\min}}\big)\e^{\frac{2\pi}{y}(\frac{c}{24}-h_{\min})}$ would be negative for small $y > 0$, a contradiction. This is the argument of \cite[Proposition 3.10]{MoR24a} with a weaker sign condition: $Z(\cal H_3)$ is not essentially positive in the sense of \cite{MoR24a} (for $J = \{b,a\}$ no row $i$ has $s_{ib} < 0$ and $s_{ia} < 0$, and a set $J$ containing $0$ never satisfies their condition), but our condition suffices to exclude $h_{\min} < 0$. Hence all conformal weights of the irreducible $V$-modules are nonnegative, and they are positive for the eight simple objects with nontrivial twist. Therefore the positivity hypothesis in our results can be replaced by the condition that none of the modules $b$, $a$, $c_0$ has conformal weight zero.
\end{remark}

\appendix
\section{Computational methods and certificates}
\label{app:computations}

The proofs distinguish exact calculations, certified exclusions, calculations modulo primes and numerical checks. Exact calculations use rational arithmetic, arithmetic in $\K = \Q(\sqrt{13},\sqrt{-3})$ in the basis $1,\sqrt{13},\sqrt{-3},\sqrt{-39}$, Gröbner bases over $\Q$ or $\Q(\sqrt{13})$, and integer arithmetic. Certified interval calculations are used only to exclude a candidate by a nonzero value. An interval containing zero gives no exclusion.

For a matrix integral at a prime $\mathfrak p$ of $\K$, its rank over $\K$ is at least the rank of its reduction. This justifies the rank lower bounds and the selection of pivots modulo primes. Membership over $\K$ is checked separately by exact coordinate identities on all rows after restriction of scalars. Membership after reduction alone is not used as a characteristic-zero assertion. Numerical calculations serve only as consistency checks and do not enter the proofs.

The supplementary material indexes the source code, coordinates, exact witnesses and verification results by the numbered statements. It specifies the computational environments and commands. The lattice construction gives the vertex algebra identities in arbitrary weights, while Section \ref{sec:griess} justifies the momentum truncation used in the certificates.

The following table identifies the scripts associated with numbered statements and the nature of their computations. It distinguishes exact verification, certified exclusions, reductions modulo primes and auxiliary numerical tests. An entry does not assert more than the precise statement and its proof establish. In particular, modular calculations are not treated as characteristic-zero lifts.

\begingroup
\small
\begin{longtable}{>{\raggedright\arraybackslash}p{3.6cm}>{\raggedright\arraybackslash}p{6.2cm}>{\raggedright\arraybackslash}p{3.3cm}}
\caption{Correspondence between numbered statements and computational scripts.}
\label{tab:scripts-statements}\\
\hline
Statement & Scripts & Nature \\
\hline
\endfirsthead
\hline
Statement & Scripts & Nature (continued) \\
\hline
\endhead
Theorem \ref{thm:EG} & \texttt{eg\_c8\_verification.py}, \texttt{bg\_trace\_criterion.py} & exact recursion, numerical modularity check \\
Proposition \ref{prop:chain} & \texttt{kac\_order3.py} & exact (finite check) \\
Proposition \ref{prop:self-dual} & \texttt{e6\_rank\_bound.py} & analytic proof and exact check \\
Proposition \ref{prop:weight-two-space} & \texttt{vvmf\_exact.py}, \texttt{weight2\_obstruction\_exact.py}, \texttt{dump\_relations.py} & exact (rational) \\
Corollary \ref{cor:tau-formulas} & \texttt{weight2\_solve.py}, \texttt{gamma1\_exclusion\_exact.py} & exact (rational) \\
Proposition \ref{prop:higher-relations} & \texttt{obstruction.py}, \texttt{e6\_tools.py}, \texttt{weight4\_solve.py} & exact (rational) \\
Lemma \ref{lem:simple-ideals} & direct inequalities in the proof & no computation required \\
Proposition \ref{prop:semisimple} & \texttt{centre\_gamma0.py} & exact \\
Lemma \ref{lem:sl2-e6} & \texttt{sl2\_e6\_superset.py}, \texttt{e6\_centralizers.py}, \texttt{e6\_orbit\_labels.py} & exact \\
Proposition \ref{prop:simple-ideals-gamma0} & \texttt{sl2\_ideal\_feasibility.py}, \texttt{weight6\_single.py} & exact \\
Theorem \ref{thm:weight-one-classification} & \texttt{sl2\_pair\_exact.py}, \texttt{e8\_joint.py}, \texttt{sextic\_joint.py}, \texttt{joint\_uniqueness.py}, \texttt{known\_sectors\_table.py}, \texttt{u2\_moments.py}, \texttt{u2\_embedding\_test.py} & exact \\
Proposition \ref{prop:u2-moduli} & \texttt{lattice\_branching.py}, \texttt{u2\_parameter\_count.py} & exact \\
Lemma \ref{lem:griess-B} & \texttt{commutant\_griess.py} & exact over $\K$ \\
Proposition \ref{prop:griess-candidates} & \texttt{griess\_subalgebras.py} & exact, Gröbner over $\Q$ \\
Theorem \ref{thm:case-two-U2} & \texttt{primaries.py}, \texttt{isotropic\_lines.py}, \texttt{isotropy\_test.py}, \path{check_primary33_exclusions.py}, \path{closure_certify_exact_intervals.py}, \path{check_closure_coverage.py} & exact rational exclusions and certified Arb witnesses \\
Proposition \ref{prop:case-two-weight-three} & \texttt{closure\_modp.py}, \texttt{modp.py}, \texttt{u3\_moments.py} & modulo $p$ (superseded by Theorem \ref{thm:exact-weight-three}) \\
Theorem \ref{thm:exact-weight-two} & \texttt{exact\_weight2\_K.py}, \texttt{Kfield.py} & exact over $\K$ \\
Theorem \ref{thm:unitarity} & \texttt{unitary\_character.py}, \texttt{phi\_automorphism\_check.py}, \texttt{unitarity\_test.py}, \texttt{lattice\_voa.py} & exact \\
Proposition \ref{prop:symmetries} & \texttt{a3\_griess\_automorphisms.py}, \texttt{choice\_conjugation.py} & exact \\
Proposition \ref{prop:weight-three-four} & \texttt{closure\_w4\_modp.py}, \texttt{modp.py} & modulo $p$ \\
Theorem \ref{thm:strong-generation} & \path{verify_f4_strict.py}, \path{aggregate_f4.py}, \path{check_zero_products.py} & exact over $\K$, Proposition \ref{prop:weight-four-certificate} and alternative dimension proof \\
Theorem \ref{thm:exact-weight-three} & \texttt{exact\_engine.py}, \texttt{certify\_w3.py}, \texttt{certify\_w4b.py} & all parts exact over $\K$ \\
Theorem \ref{thm:jiao-identification} & \texttt{i2\_weyl.py}, \texttt{i4\_torus.py}, \path{i5_membership_exact.py}, \path{i6_generation_certificate.py} & exact membership and rank lower bound modulo $p_3$ \\
Proposition \ref{prop:centralizer-two} & \texttt{j1\_centralizer.py}, \texttt{j5\_explicit\_automorphism.py} & root-theoretic proof and exact polynomial check \\
Proposition \ref{prop:commutant-biweights} & \texttt{c2\_exact\_C4.py}, \texttt{commutant\_dims.py}, \texttt{griess\_subalgebras.py} & exact over $\K$ \\
Proposition \ref{prop:case-one-control} & \texttt{exact\_isotropy.py}, \texttt{isotropy\_test.py} & exact, certified (Arb) \\
Proposition \ref{prop:case-one-two} & \texttt{isotropic\_lines.py}, \texttt{r6\_certify\_sub2.py}, \texttt{r5\_sub2\_testb.py} & exact, certified (Arb) \\
Proposition \ref{prop:riesen-identification} & \texttt{r1\_identify\_riesen.py}, \texttt{r2\_locate\_riesen.py}, \texttt{griess\_subalgebras.py} & exact \\
Proposition \ref{prop:centralizer-one} & \texttt{j1\_centralizer.py}, \texttt{j2\_solve.py}, \texttt{j6\_caseI\_action.py}, \texttt{r3\_lines\_and\_Z.py} & root-theoretic proof and exact check \\
Proposition \ref{prop:riesen-dims} & \texttt{b1\_UR\_closure\_modp.py} & modulo $p$ (rigorous lower bounds) \\
Proposition \ref{prop:rational-forms} & \texttt{d1\_rational\_forms.py}, \texttt{j6\_caseI\_action.py} & exact \\
Proposition \ref{prop:simple-current} & \path{check_simple_current_word.py} & analytic word and energy proof, exact checks \\
Theorem \ref{thm:zhu-modp} & \texttt{z0\_highest.py}, \texttt{z1\_products.py}, \texttt{z2\_relations\_modp.py}, \texttt{z3\_corner\_modp.py}, \texttt{z4\_upper\_algebra\_modp.py} & modulo $p$ \\
Theorem \ref{thm:classification-one} & \texttt{b1\_UR\_closure\_modp.py} & modulo $p$ (rigorous lower bounds) \\
Lemma \ref{lem:two-lattices} & \texttt{s20\_genus\_check.py}, \texttt{s22\_lattice\_weights.py} & exact (integers) \\
Proposition \ref{prop:associated-lattice} & Table \ref{tab:lowest-spaces} and the proof & no separate lattice-weight computation \\
\hline
\end{longtable}
\endgroup

\subsection{Coverage of the certificates}
\label{app:verification-record}

The trace certificates give the dual forms, vanishing relations, matrix ranks and eighteen moments of $U_2$. The module searches enumerate the entire finite domain specified in the proof of Theorem \ref{thm:weight-one-classification}.

For Theorem \ref{thm:case-two-U2}, the thirteen excluded commutant candidates give thirty-nine rational projective charts with unit Gröbner bases. The pairwise closure calculation gives $62$ certified nonzero interval witnesses, which exclude $254$ of the $256$ choices. The coverage check uses only these certified witnesses. The primary-space computations retain all characteristic-polynomial factors and do not discard an inconclusive eigenspace.

The weight-four certificate of Proposition \ref{prop:weight-four-certificate} checks the complete input domains, the identity of each pivot word, integrality at the chosen prime, full rank and every coordinate of each membership identity. It covers all $131$ nonzero Cartan blocks and checks the $886$ zero inputs separately. Missing or truncated data, duplicate labels and nonzero residuals are rejected. This establishes \textup{(F4)} without an external existence theorem. Corollary \ref{cor:weight-four-from-haagerup} gives a separate dimension proof when an ambient Haagerup algebra is available.

For Proposition \ref{prop:commutant-biweights}, the supplementary data give the seven exact commutant vectors, the two Virasoro matrices and the primary vectors. Their matrix identities and ranks admit a separate verification in rational arithmetic. This verification is distinct from the lattice-mode checks of the vertex-product identities.

\paragraph{Independent mode comparisons.}
A second evaluator, using simultaneous normal ordering and direct exponential coefficients, agrees with the iterate-formula implementation on $24643$ basis-state products in the stated weight window and $310$ generator products. The coefficient-field multiplication and three Virasoro identities are checked separately. These comparisons use the supplied generator coordinates. They do not independently reconstruct those coordinates from the printed definitions, or regenerate every coordinate of \textup{(F4)} by a second method.

\paragraph{The modular Zhu calculation.}
All $131698$ overlap ambiguities and $151$ inclusion ambiguities of the final presentation are checked. At $8388571$, the algebra defined by the projected relations has dimension $762$, nondegenerate trace form and the twelve stated blocks. These checks concern that presentation. They do not prove a characteristic-zero lift, recompute the vertex-algebra origin of every relation or show that every block survives in the Zhu corner. The agreement with the calculation at a second prime has the status specified in Remark \ref{rem:zhu-status}.

\paragraph{AI statement.}
The author acknowledges the use of AI tools for language polishing, exploratory mathematical discussions, the writing of the verification scripts and feedback during the development and preparation of this manuscript. The author takes full responsibility for the mathematical content of the manuscript.

\paragraph{Competing interests.} The author declares that he has no competing interests.

\paragraph{Data availability.}
The source code, certificates and reproduction instructions are available at \url{https://github.com/cedricarhancet/haagerup-c8}. The two complete weight-four coordinate archives are distributed as release assets. Table \ref{tab:scripts-statements} relates the computational scripts to the numbered statements.
\noindent\textit{Author's GitHub profile.} \url{https://github.com/cedricarhancet}

\end{document}